\documentclass[12pt]{article}
\usepackage{rotating}
\input tcilatex
\input xy
\xyoption{all}
\usepackage{amsfonts,epsf}
\newcommand{\psdiag}[3]{\hspace{1mm}\raisebox{-#1mm}{\epsfysize#2mm
\epsffile{#3.eps}}\hspace{1mm}}

\begin{document}
\author{Rui Pedro Carpentier\\
\\
{\small\it Departamento de Matem\'{a}tica, Centro de
Matem\'{a}tica e Aplica\c c\~oes}\\ {\small\it and Centro de
An\'alise Matem\'{a}tica, Geometria e Sistemas Din\^amicos}\\
{\small\it Instituto Superior T\'{e}cnico}\\ {\small\it Avenida
Rovisco Pais, 1049-001 Lisboa}\\ {\small\it Portugal}}
\title{Representations of non-singular
planar tangles by operators}
\date{12rd April, 2005}
\maketitle
\begin{abstract}
In this paper we study how to distinguish two embeddings of a
finite collection of disjoint circles into the plane up to planar
isotopy. We adopt the spirit of the approach by V. Turaev,
Operator Invariants of Tangles, Math. USSR-Izv. 35  (1990),
411--444, by considering a category of planar tangles and
representing it in an ``algebraic'' category. From this we can
extract a numerical invariant for embeddings of a finite
collection of disjoint circles and this invariant is, up to
certain choices, complete.
\end{abstract}

\section{Introduction}

Let ${\bf{PT}}$ be the {\it category of non-singular planar
tangles} whose objects are finite sets of points in the real line
identified up to 1-dimensional isotopies\footnote{these objects
can be regarded as finite ordinal numbers $\emptyset , \{0\},
\{0,1\}, ...$}, and whose morphisms between two objects $O_1$ and
$O_2$ are piecewise regular 1-dimensional manifolds, with boundary
$O_1 \times\{1\} \cup O_2 \times\{0\}$, embedded in
$\mathbb{R}\times [0,1]$ identified up to planar isotopies:

$$\psdiag{10}{20}{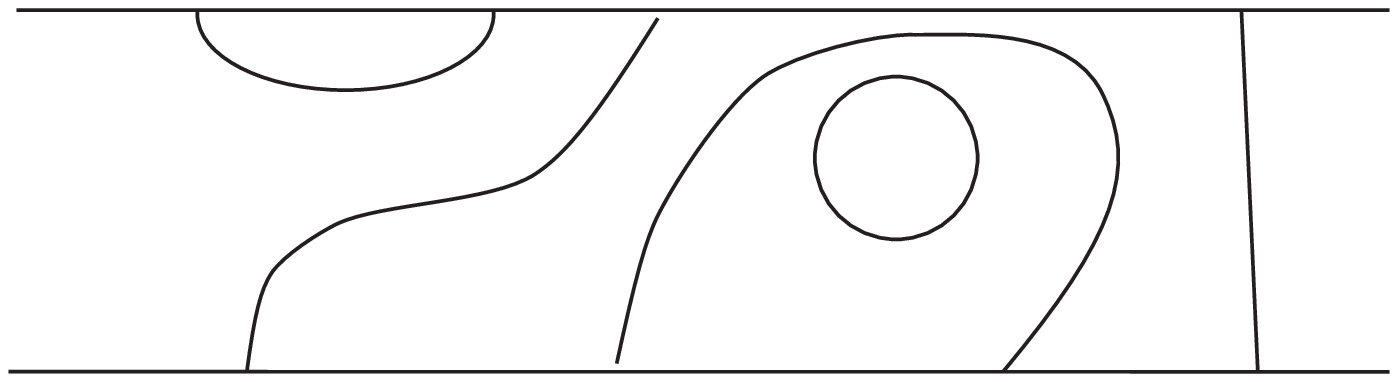}$$

The composition of two morphisms $t_1$ and $t_2$ is defined by
$$t_2\circ t_1 := g(f(t_1)\cup t_2)$$ where $f(x,y)=(x,y+1)$ and
$g(x,y)=(x,y/2)$.

$$t_1 = \psdiag{8}{16}{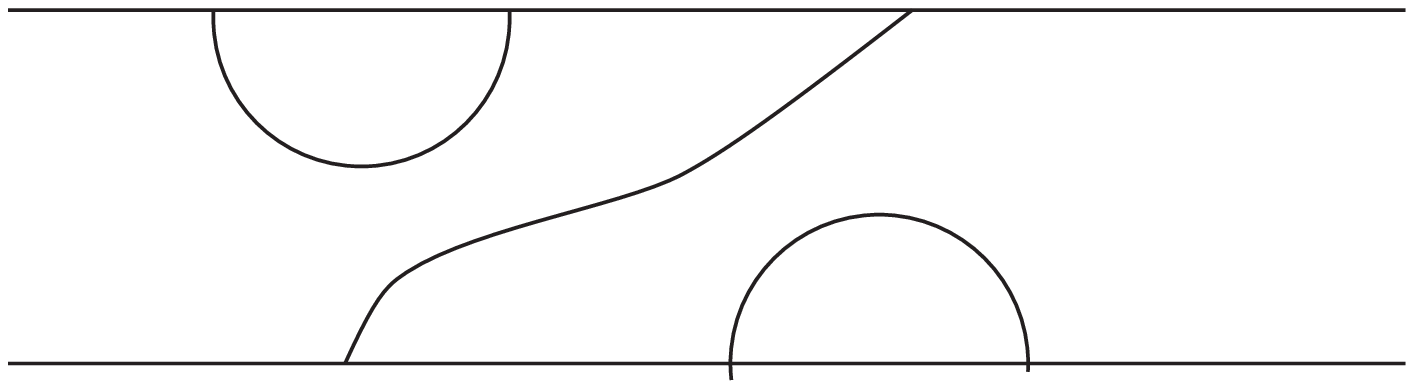} \qquad t_2 =
\psdiag{8}{16}{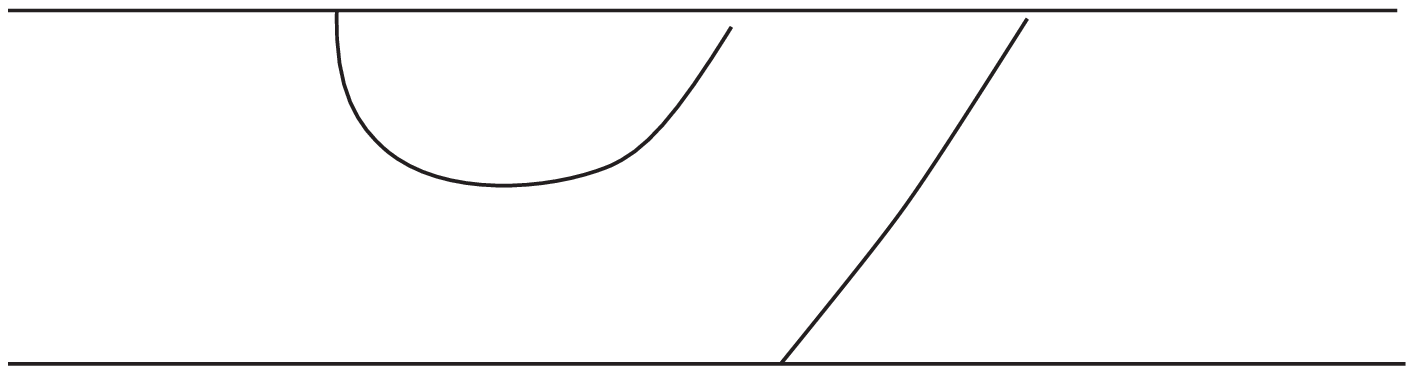}$$

$$ \psdiag{8}{32}{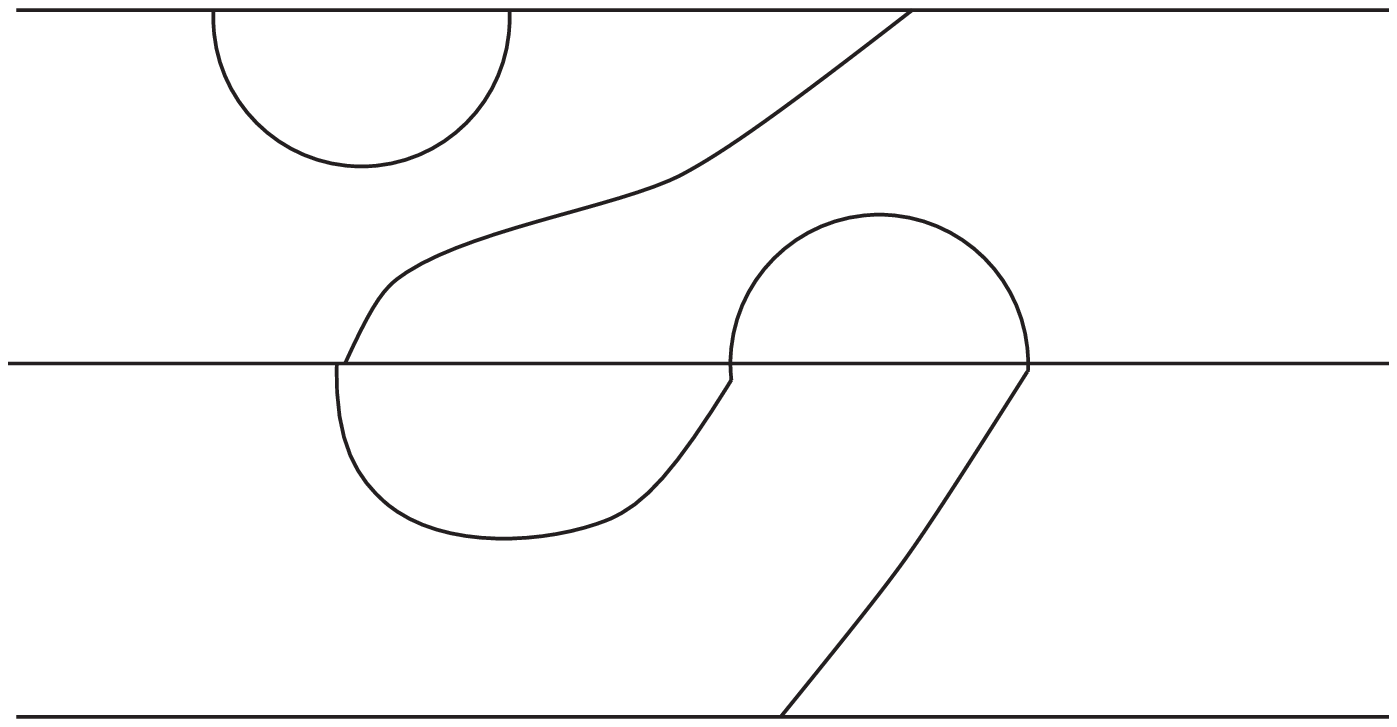} \longrightarrow t_2\circ t_1 =
\psdiag{8}{16}{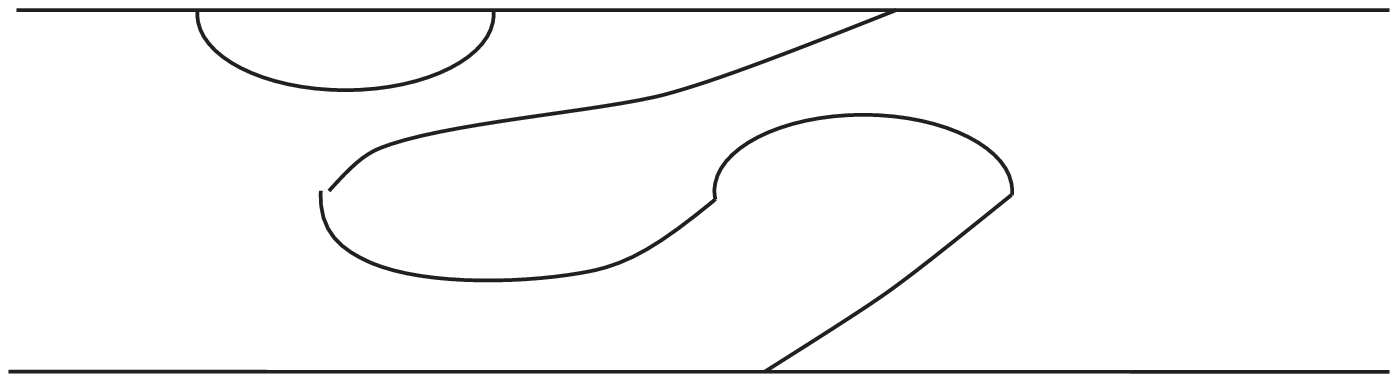}$$

In this paper, the downward direction composition is used, some
authors use the opposite direction.

 This category has the following presentation:

The generators are morphisms $\hat{t}_{n,k} \in
\hom(\{1,...,n-1\},\{1,...,n+1\})$ connecting $(i,1)$ to $(i,0)$
if $i\leq k-2$, $(i,1)$ to $(i+2,0)$ if $i\geq k-1$ and $(k-1,0)$
to $(k,0)$: $$\hat{t}_{n,k}=\psdiag{10}{20}{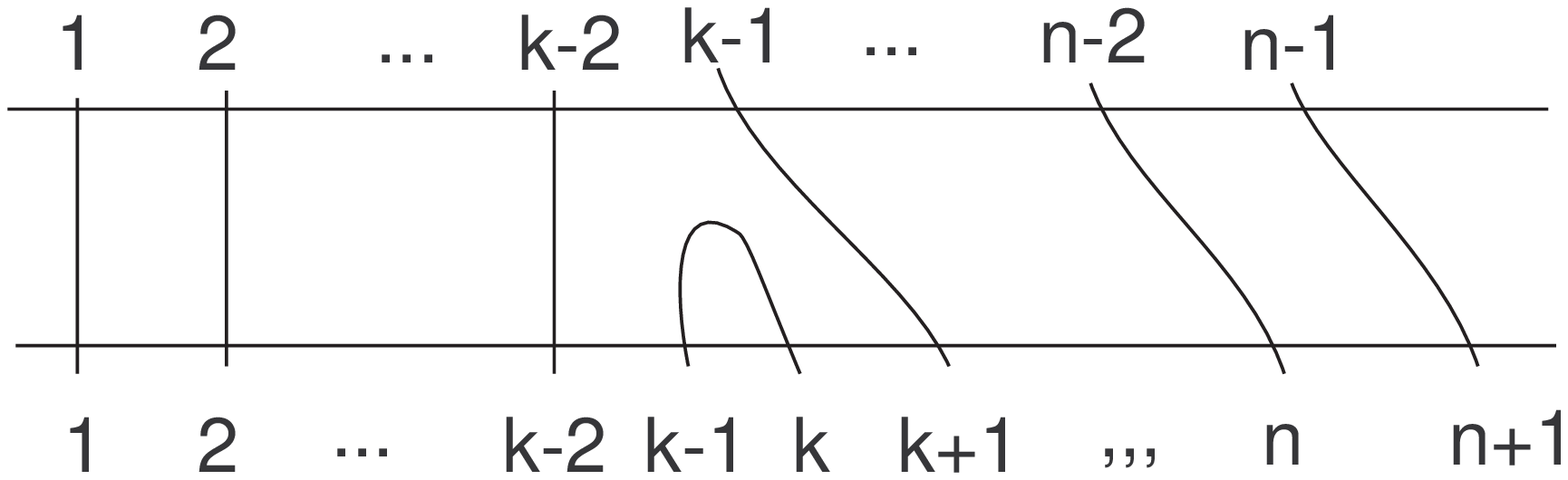}$$ and
$\check{t}_{n,k} \in \hom(\{1,...,n+1\},\{1,...,n-1\})$ connecting
$(i,1)$ to $(i,0)$ if $i\leq k-2$, $(i+2,1)$ to $(i,0)$ if $i\geq
k-1$ and $(k-1,1)$ to $(k,1)$:
$$\check{t}_{n,k}=\psdiag{10}{20}{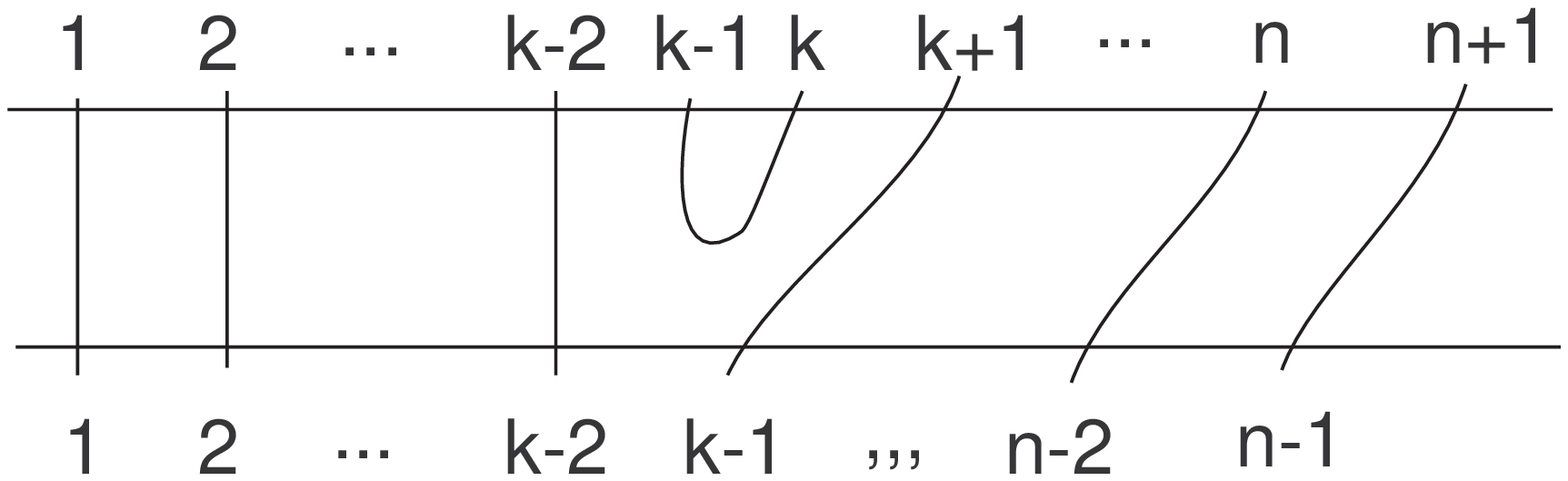}$$ for any $k,n\in
\mathbb{N}$ with $2 \leq k \leq n+1$.

For the rest of this paper it is better to number the intervals
instead of the points:

$$\hat{t}_{n,k}=\psdiag{9}{18}{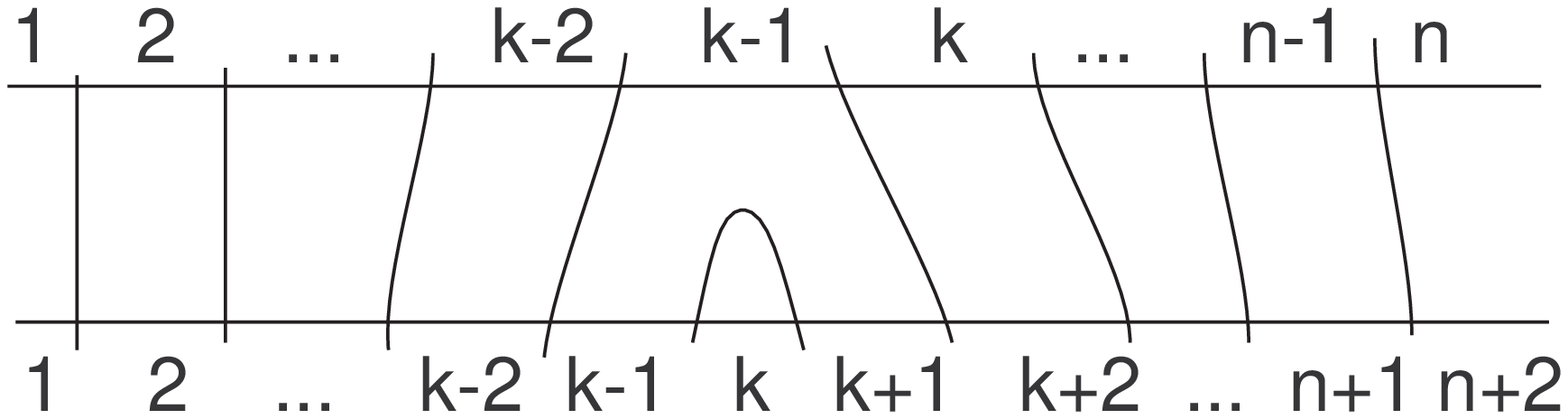}$$
$$\check{t}_{n,k}=\psdiag{9}{18}{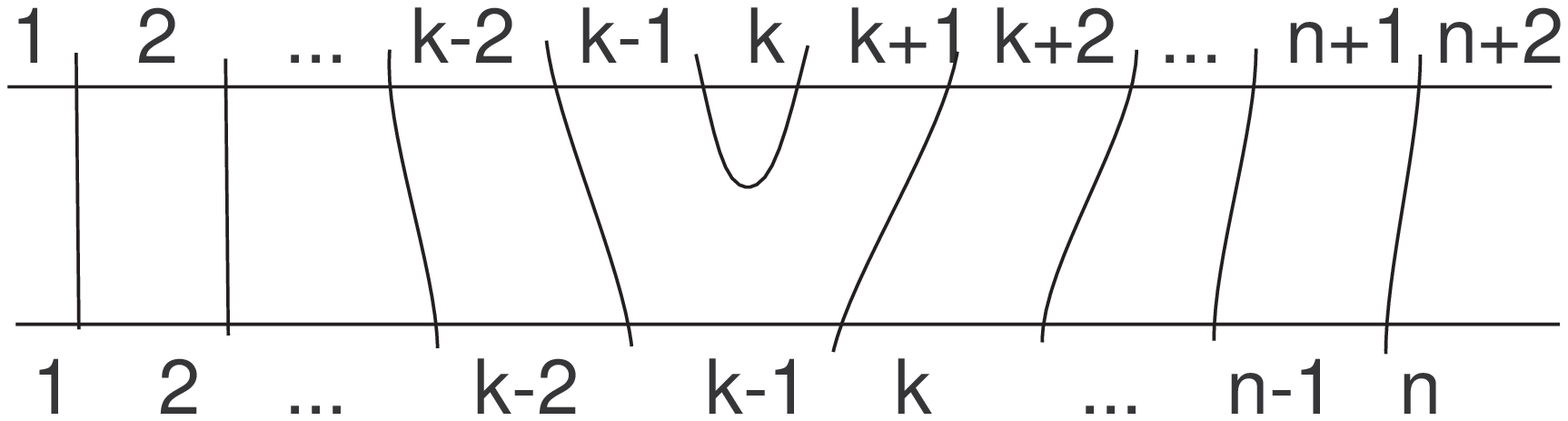}$$

These generators satisfy the following relations:

$$\psdiag{15}{30}{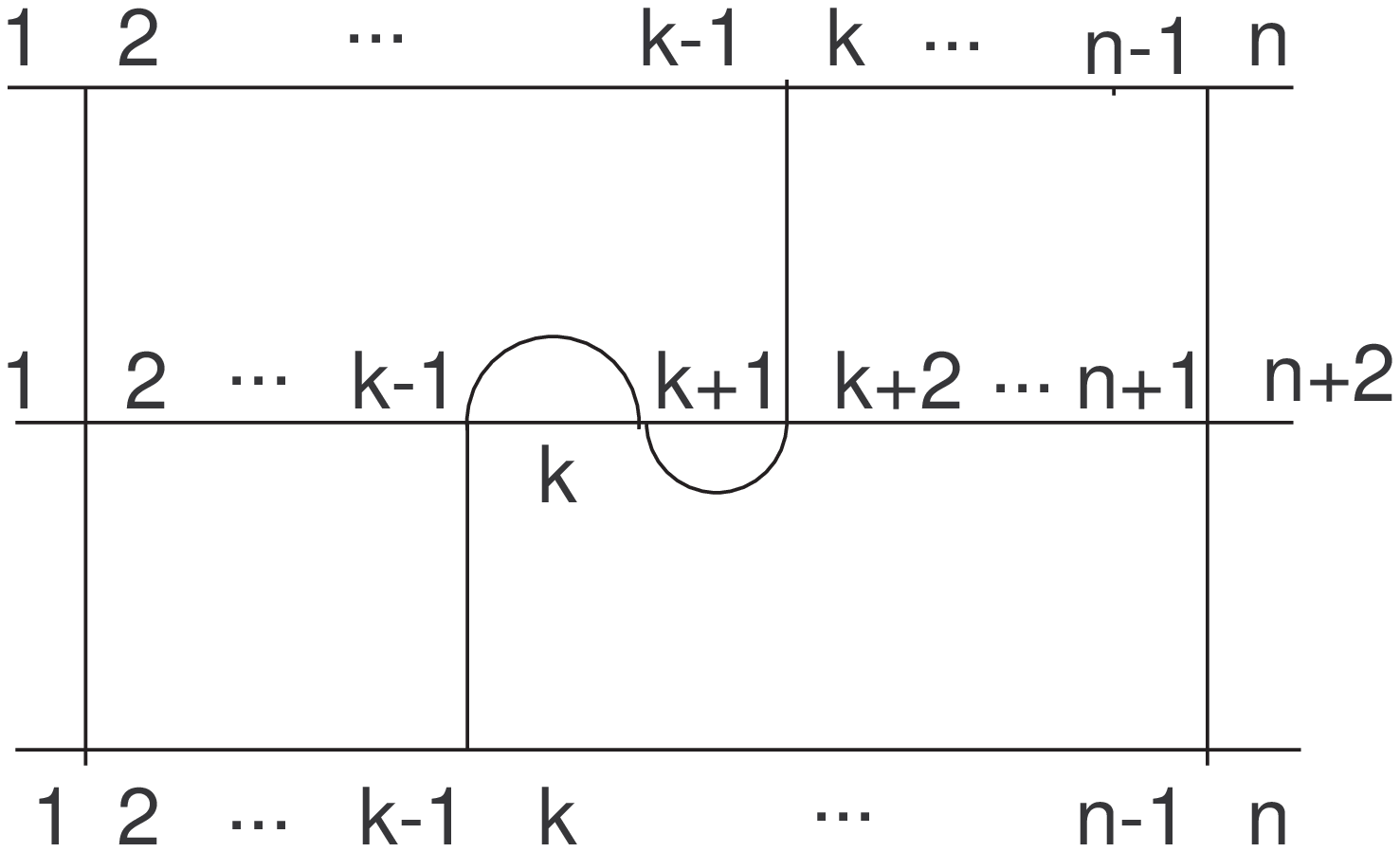}=\psdiag{15}{30}{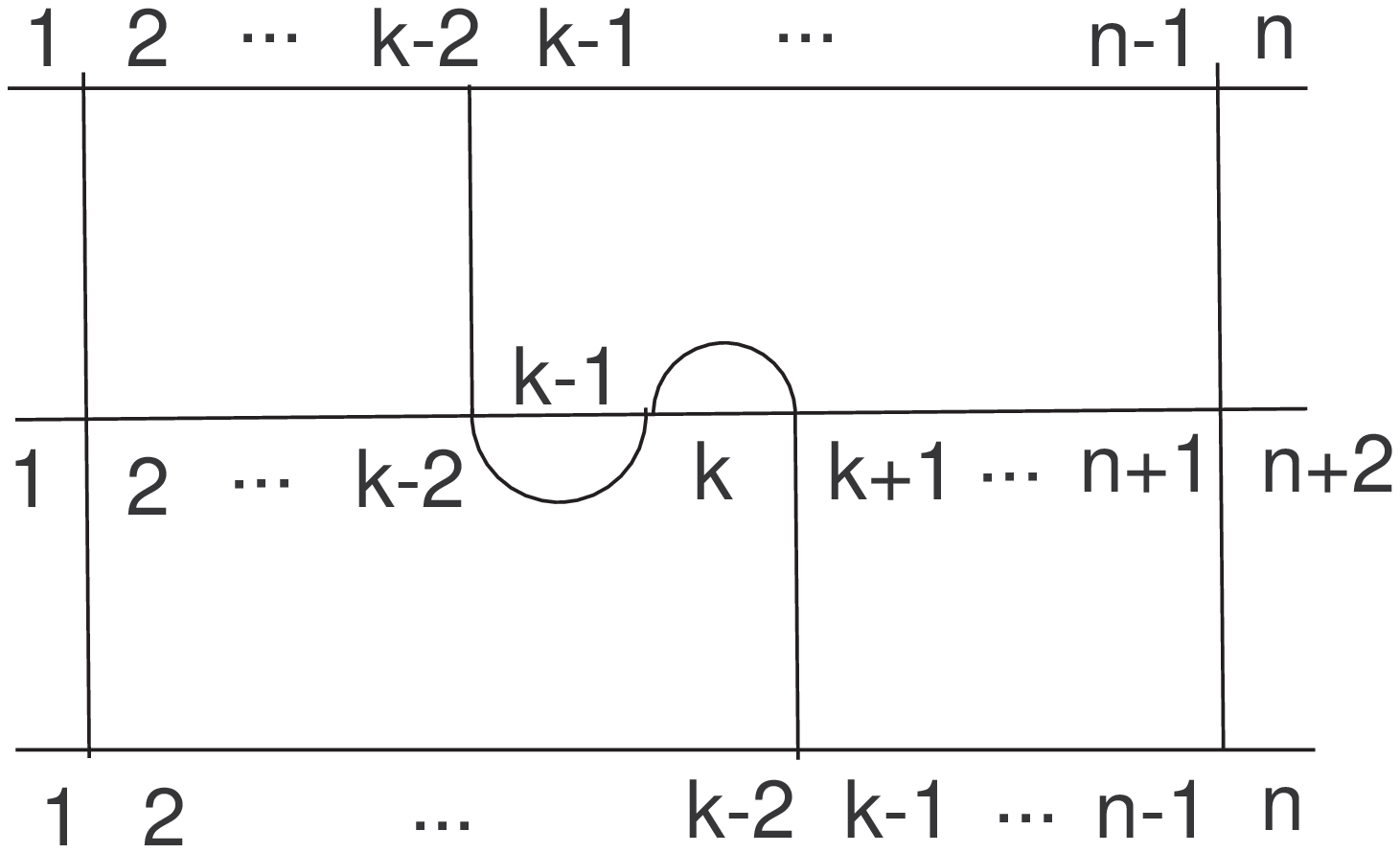}=\psdiag{15}{30}{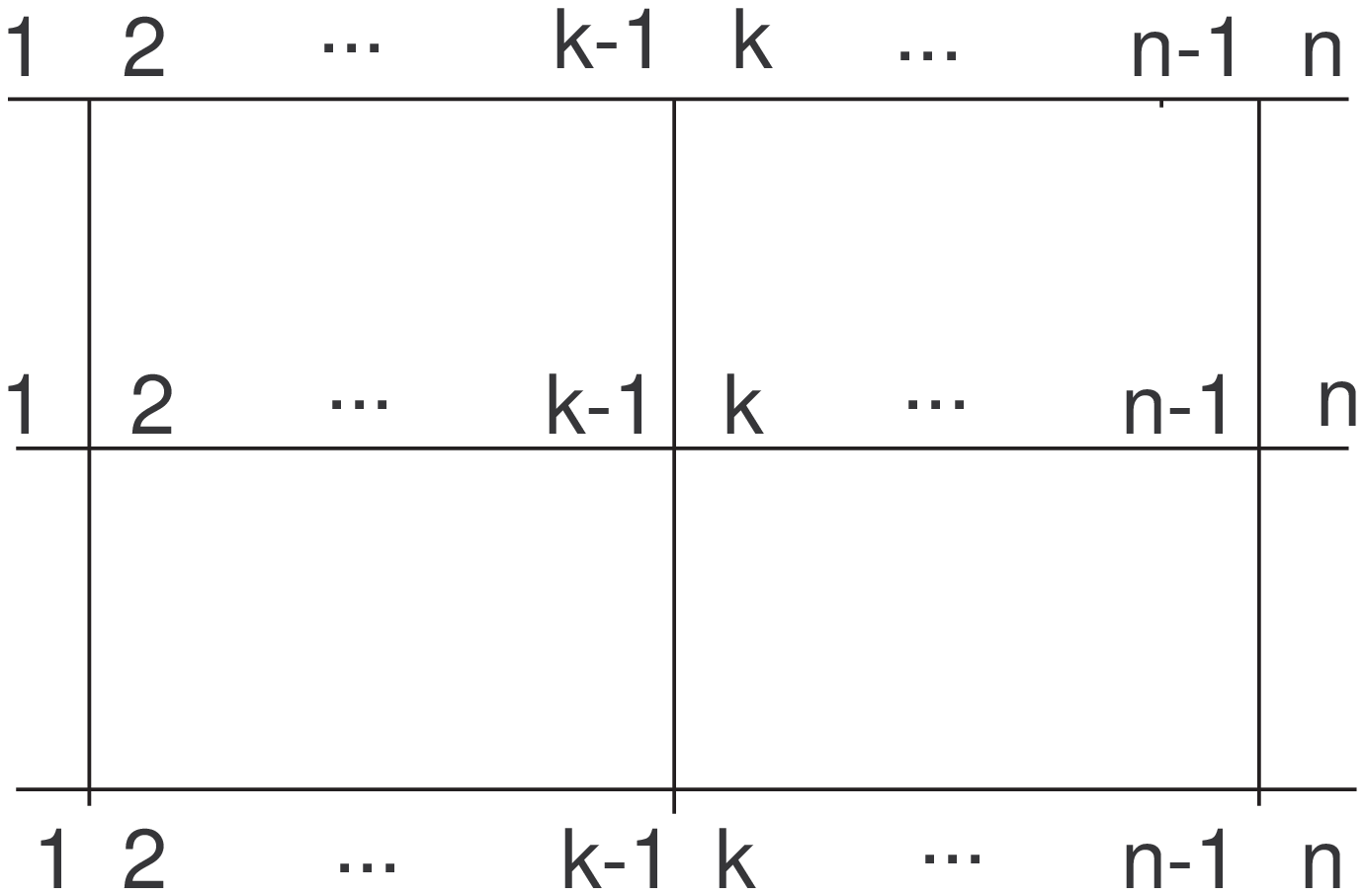}$$

$$\check{t}_{n,k+1}\circ\hat{t}_{n,k}=
\check{t}_{n,k-1}\circ\hat{t}_{n,k}= id_n$$ where $id_n := \{1,
...,n\}\times [0,1]$ is the identity morphism on $\{1, ...,n\}$;

$$\psdiag{15}{30}{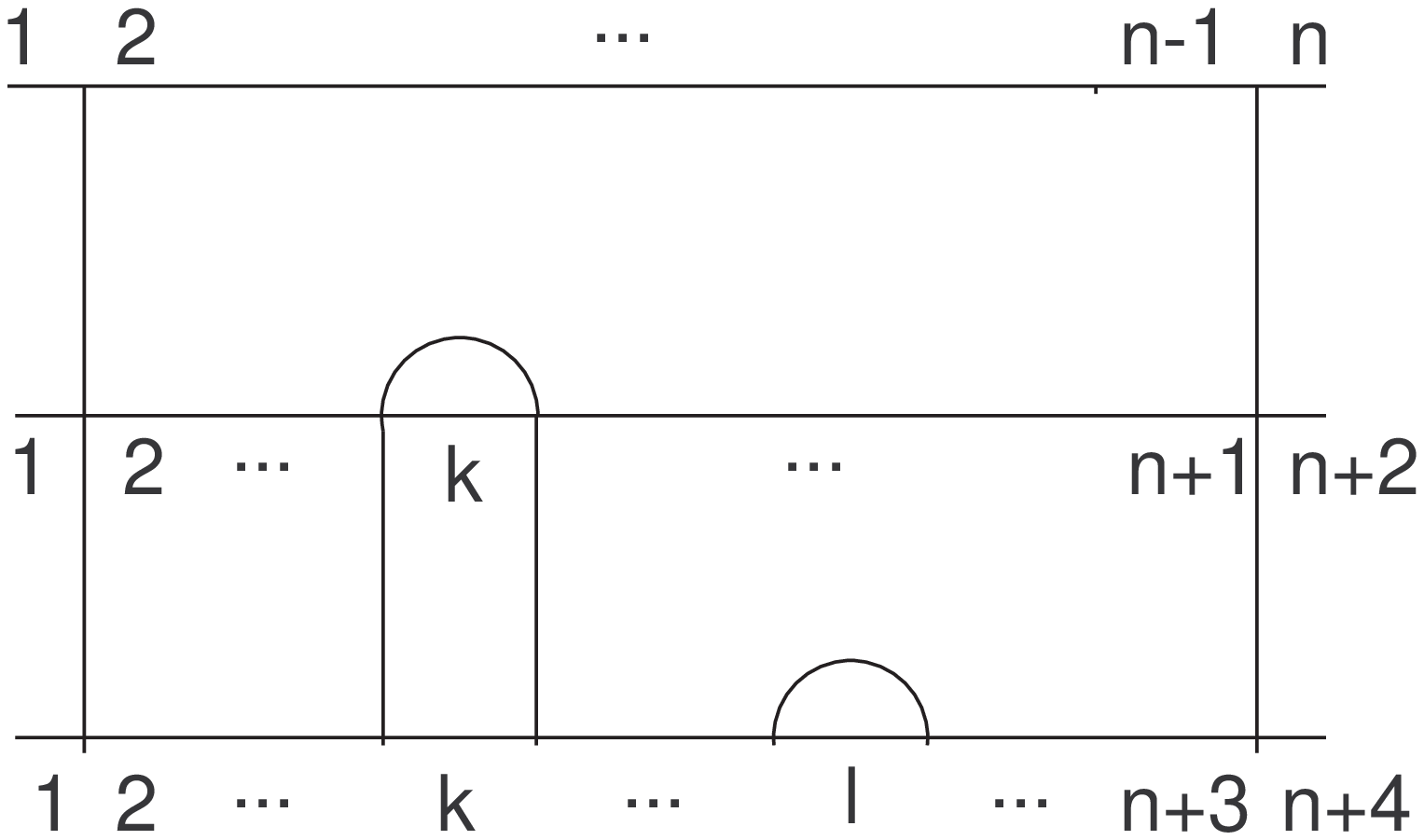}=\psdiag{15}{30}{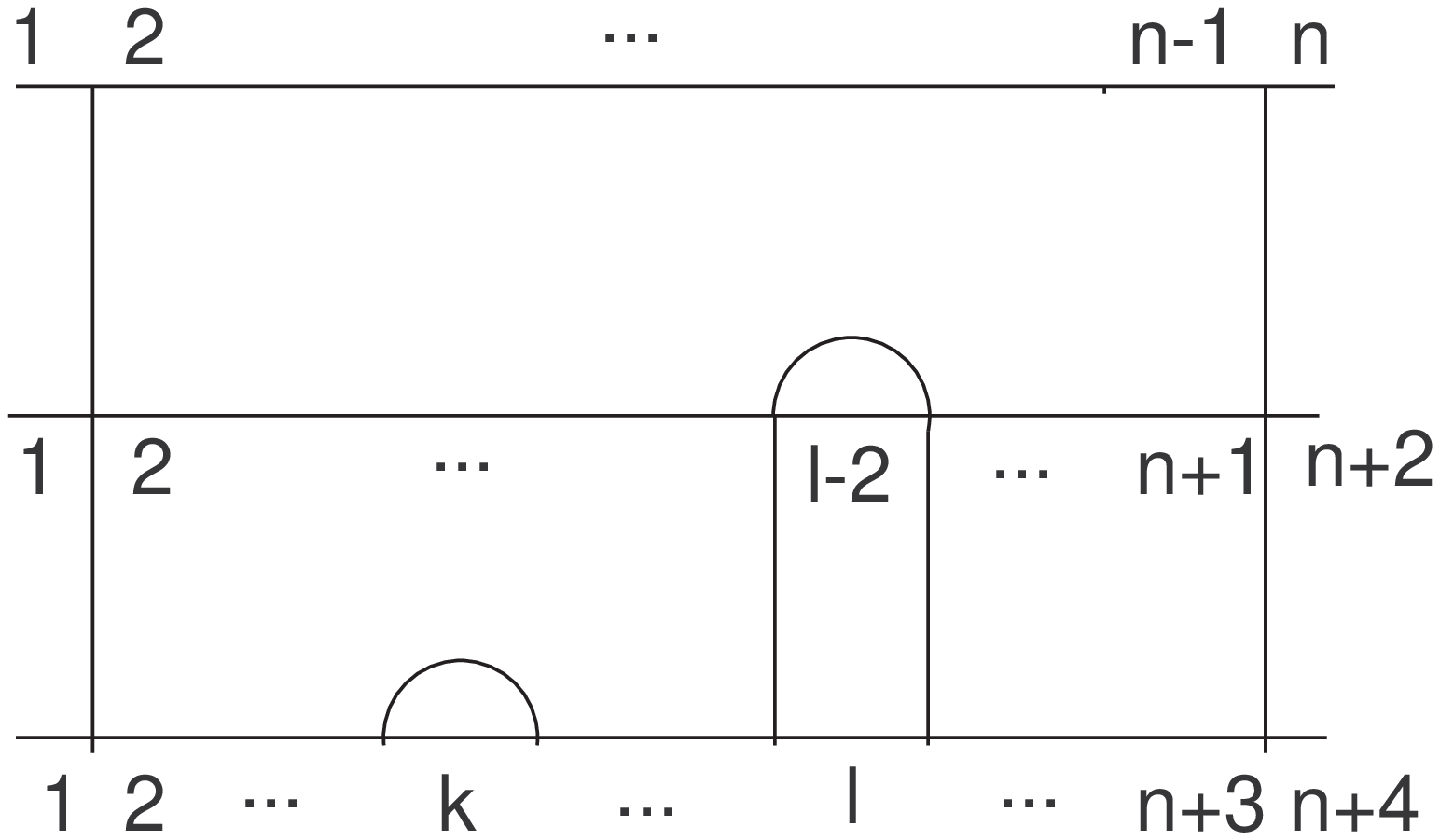}$$

$$\hat{t}_{n+2,l}\circ\hat{t}_{n,k}=
\hat{t}_{n+2,k}\circ\hat{t}_{n,l-2}$$ for $l\geq k+2$;

$$\psdiag{15}{30}{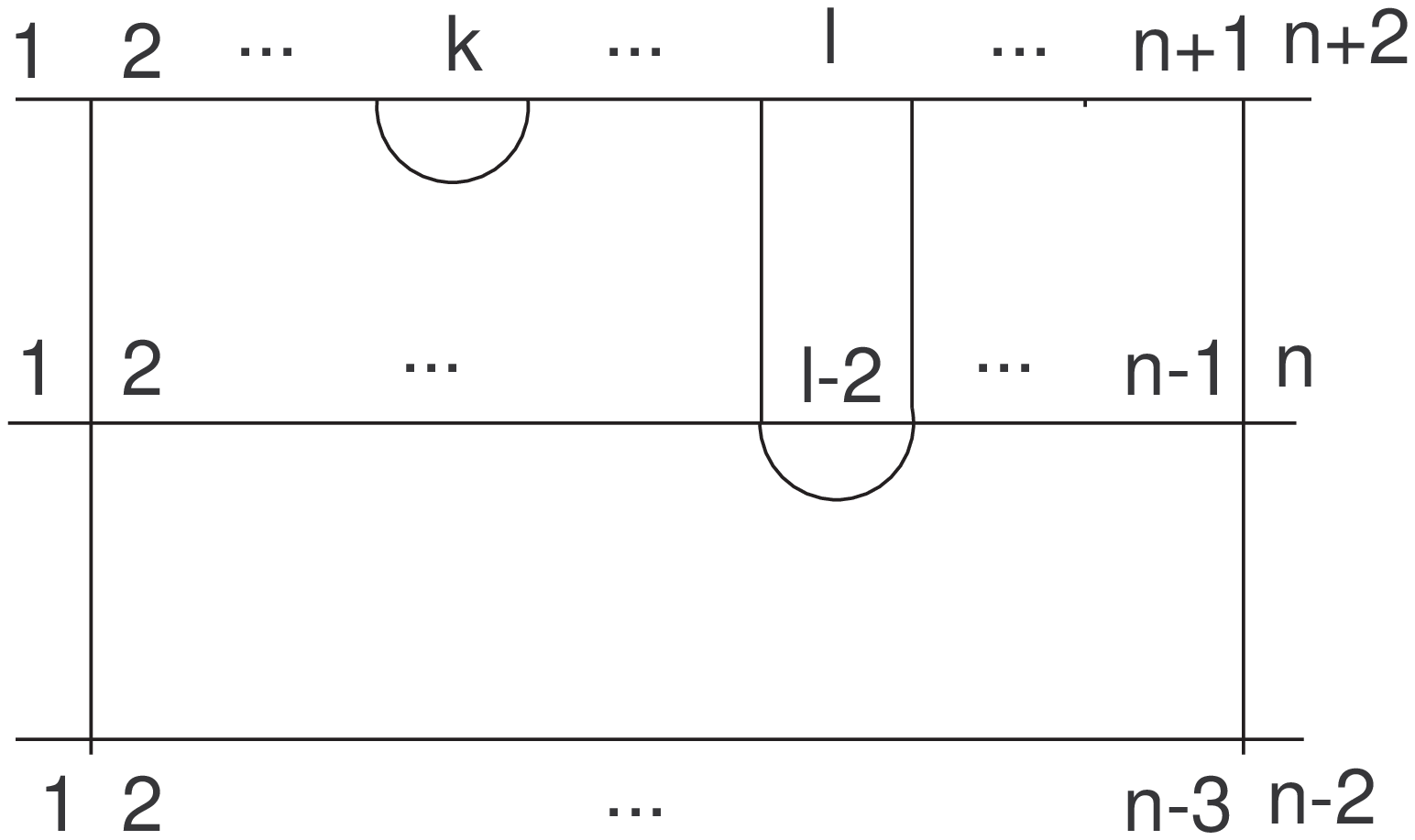}=\psdiag{15}{30}{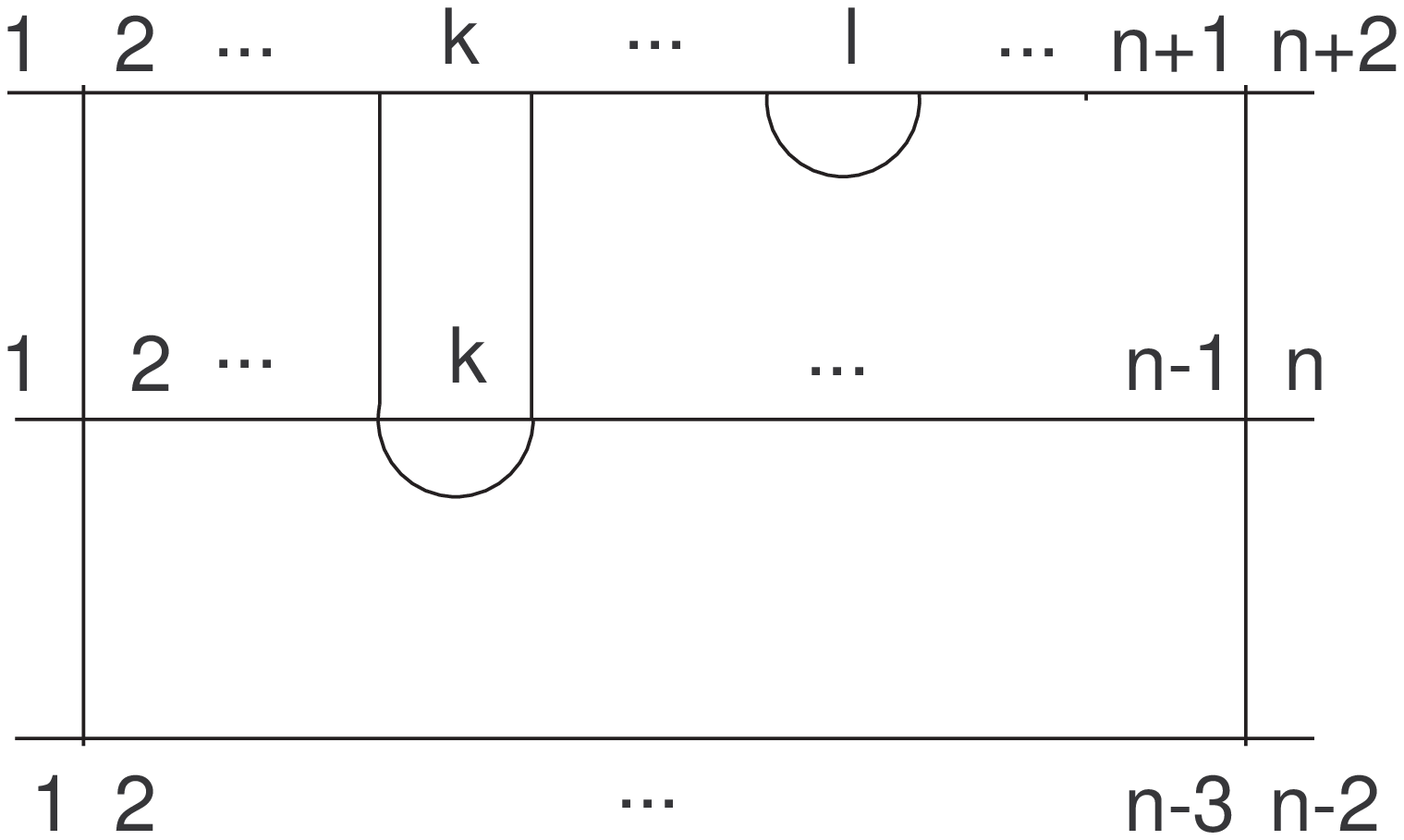}$$

$$\check{t}_{n-2,l-2}\circ\check{t}_{n,k}=
\check{t}_{n-2,k}\circ\check{t}_{n,l}$$ for $l\geq k+2$;

$$\psdiag{15}{30}{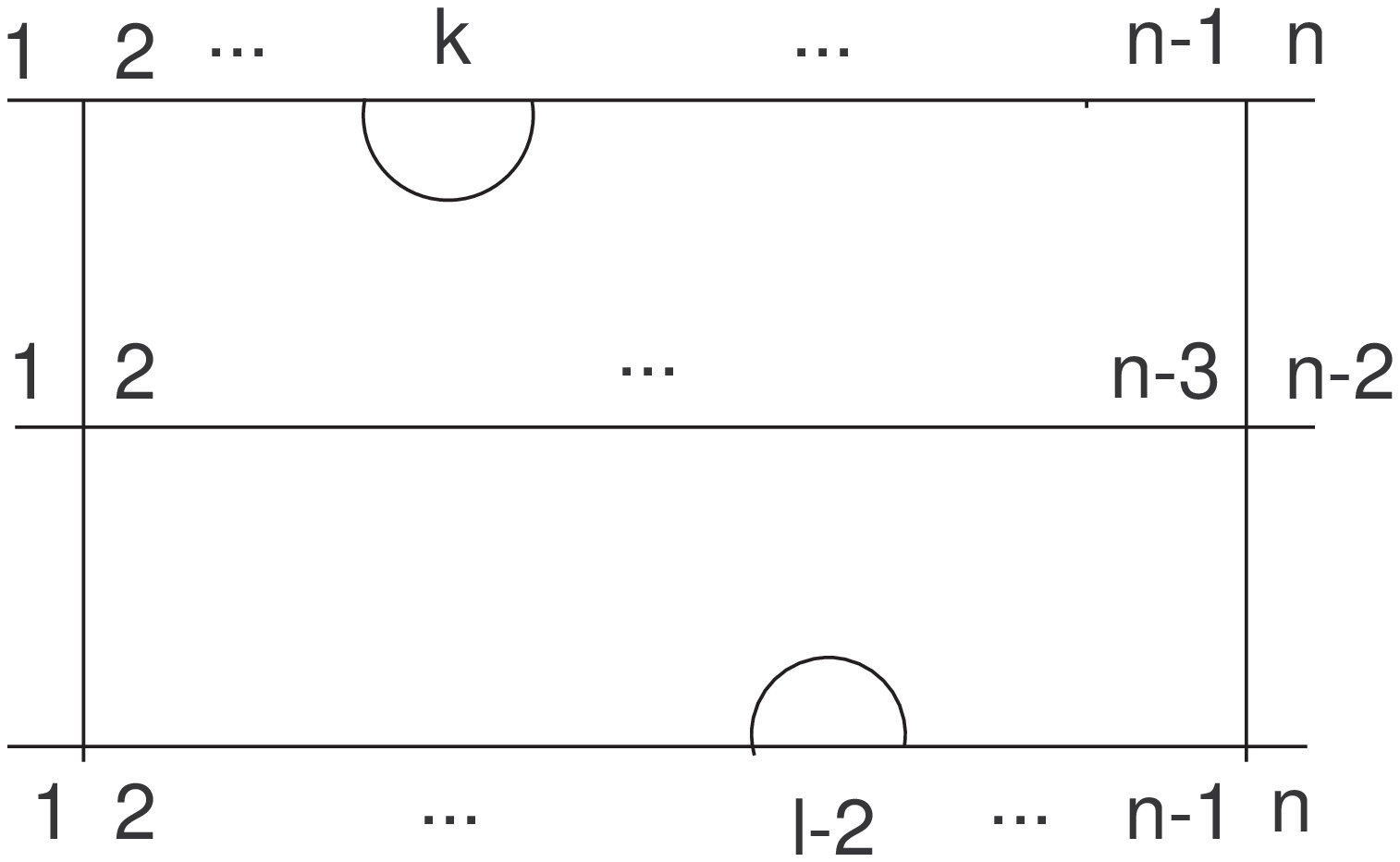}=\psdiag{15}{30}{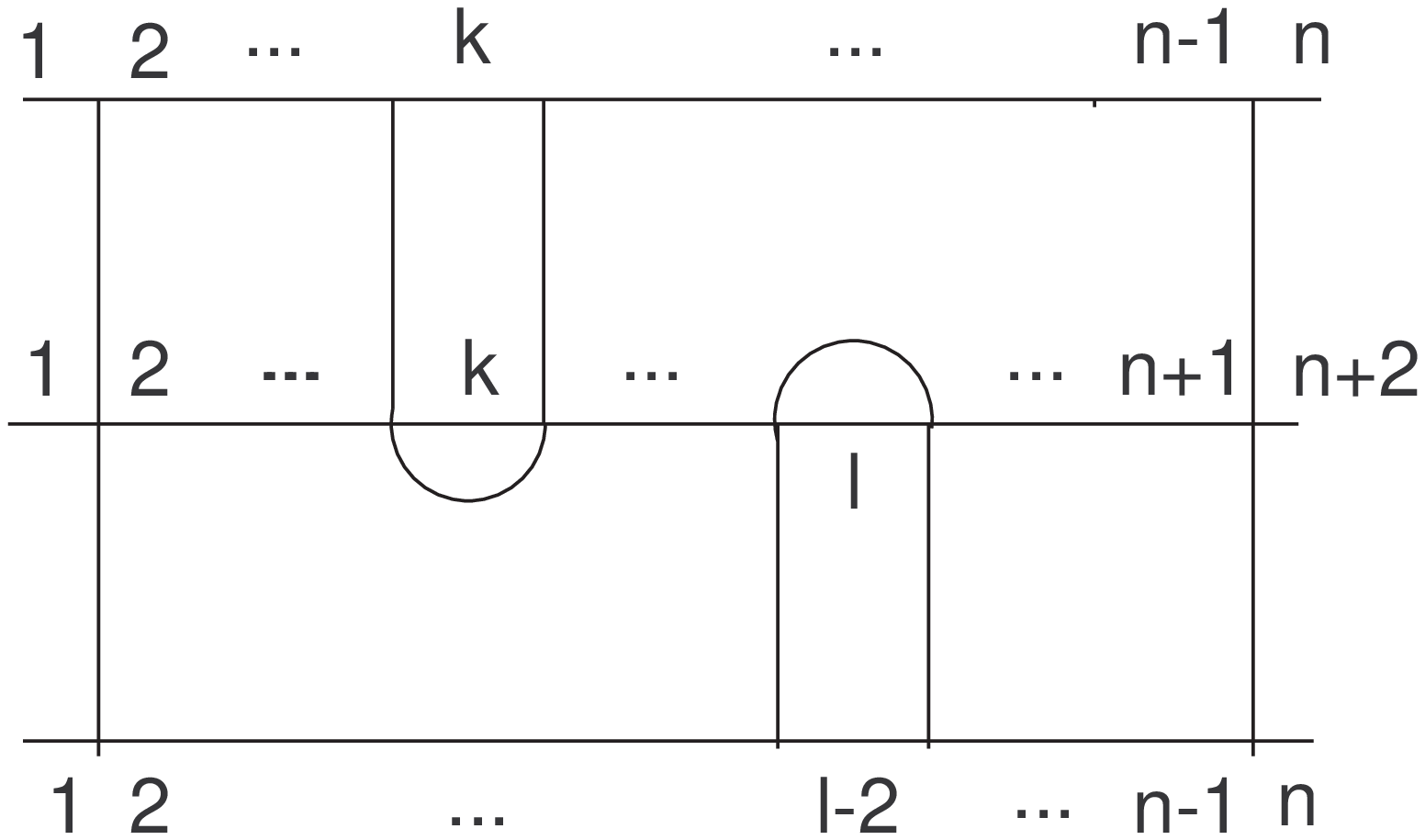}$$

$$\hat{t}_{n-2,l-2}\circ\check{t}_{n-2,k}=
\check{t}_{n,k}\circ\hat{t}_{n,l}$$ for $l\geq k+2$ and

$$\psdiag{15}{30}{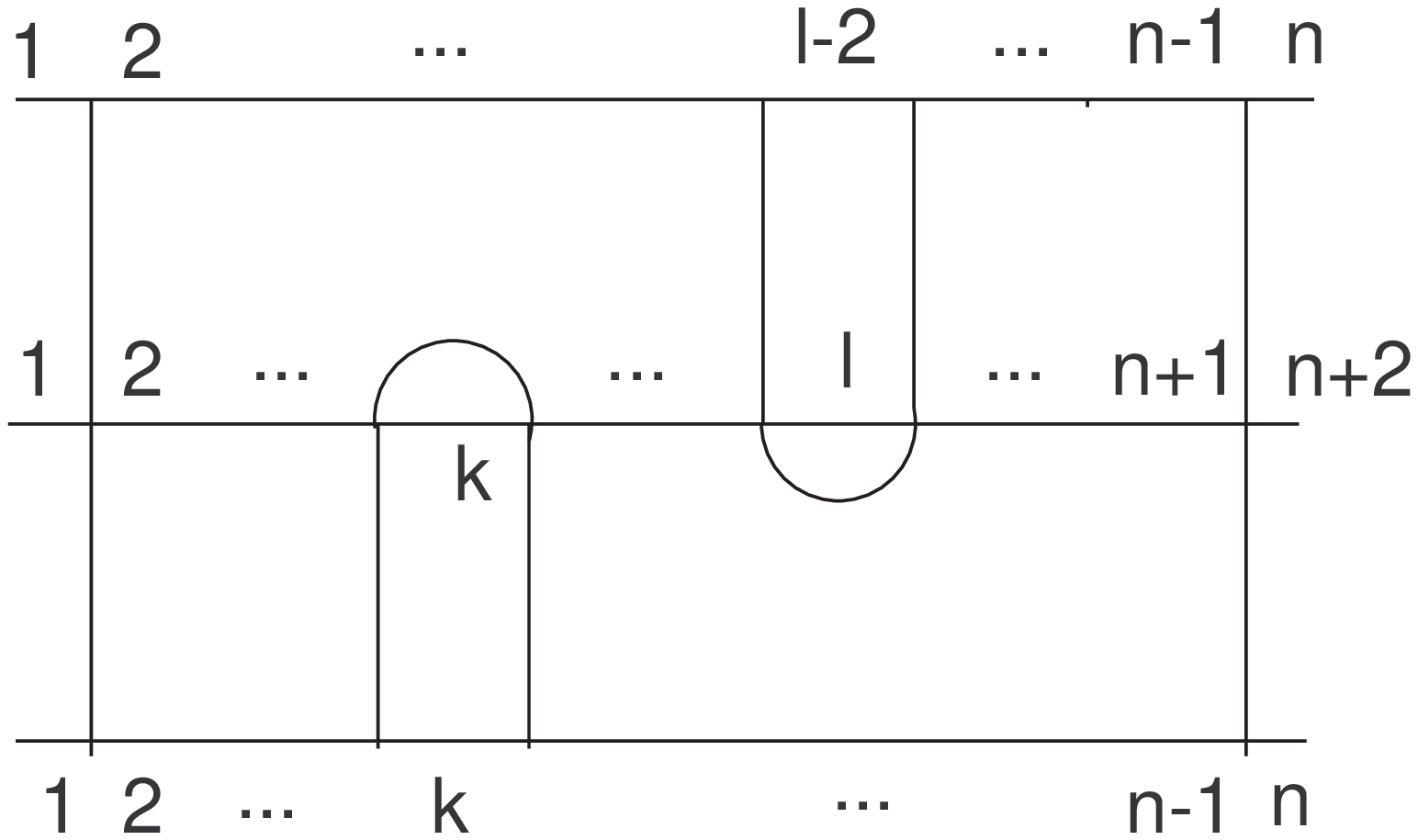}=\psdiag{15}{30}{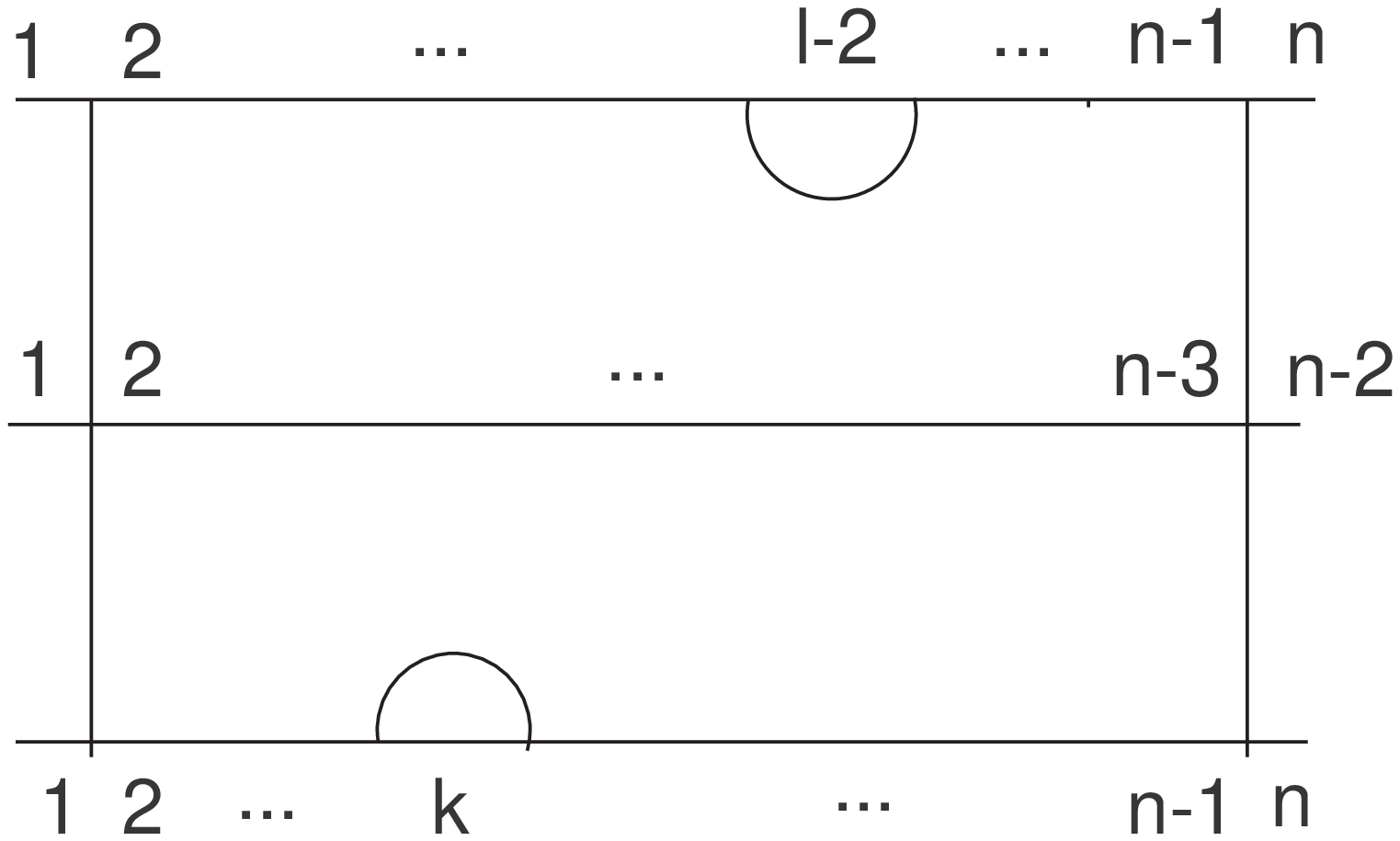}$$

$$\check{t}_{n,l}\circ\hat{t}_{n,k}=
\hat{t}_{n-2,k}\circ\check{t}_{n-2,l-2}$$ for $l\geq k+2$.

Now we want to represent non-singular planar tangles by functions.
That is, we want to find functions $\hat{T}_{n,k}$ and
$\check{T}_{n,k}$ satisfying the following relations:

$$\check{T}_{n,k+1}\circ\hat{T}_{n,k}=
\check{T}_{n,k-1}\circ\hat{T}_{n,k}= id$$

$$\hat{T}_{n+2,l}\circ\hat{T}_{n,k}=
\hat{T}_{n+2,k}\circ\hat{T}_{n,l-2}$$ for $l\geq k+2$;

$$\check{T}_{n-2,l-2}\circ\check{T}_{n,k}=
\check{T}_{n-2,k}\circ\check{T}_{n,l}$$ for $l\geq k+2$;

$$\hat{T}_{n-2,l-2}\circ\check{T}_{n-2,k}=
\check{T}_{n,k}\circ\hat{T}_{n,l}$$ for $l\geq k+2$ and

$$\check{T}_{n,l}\circ\hat{T}_{n,k}=
\hat{T}_{n-2,k}\circ\check{T}_{n-2,l-2}$$ for $l\geq k+2$.

Our proposal is to represent ${\bf{PT}}$ by the following category
${\bf{PI}}_{\mathbb{M}}$ whose objects are ${\cal O}_{1}, {\cal
O}_2, {\cal O}_3, ...$ where ${\cal O}_n$ is the set of pairs
$(R,\overrightarrow{v})$ such that $R=[r_{i,j}]$ is a $n\times n$
Boolean matrix satisfying the following properties:

\begin{description}
\item[E1.] $R\geq I$ (where $I$ is the identity matrix);
\item[E2.] $R^t =R$ (the matrix is symmetric);
\item[E3.] $R^2 =R$ (the matrix is idempotent);
\item[T1.] If $r_{i,j}=1$ then $|i-j|$ is even;
\item[T2.] For any $\alpha \leq\beta \leq\gamma \leq\delta$,
$r_{\alpha,\gamma}r_{\beta,\delta} \leq
r_{\alpha,\beta}r_{\beta,\gamma}r_{\gamma,\delta}$;
\item[T3.] For any $\alpha < \beta$ if $r_{\alpha,\beta}=1$
then either $r_{\alpha +1,\beta -1}=1$ or there exists $\gamma$
between $\alpha$ and $\beta $ such that $r_{\alpha,\gamma}=1$.
\end{description}
and $\overrightarrow{v}$ is an array of $n$ entries with values in
a chosen lattice ordered monoid $\mathbb{M}$ such that it is fixed
by the action induced by $R$ which we will define later:

\begin{description}
\item[EC.] $R*\overrightarrow{v}=\overrightarrow{v}$.
\end{description}

Note: The properties E1, E2 and E3 represent an equivalence
relation, and the properties T1, T2 and T3 have topological
motivations (see the explanation in section 3.1).

A morphism between $ \mathcal{O}_m$ and $\mathcal{O}_n$ is just a
set function between the sets $ \mathcal{O}_m$ and
$\mathcal{O}_n$.

\section{Algebraic interlude}

\begin{definition}
The canonical Boolean algebra $\mathcal{B}$ is the set $\{0,1\}$
with two binary operations: the sum $+$ and the multiplication
$\cdot$, and a unary operation the negation $\neg$ such that
$(\{0,1\},+,\cdot)$ is the (unique) semi-ring with $1+1=1$, $\neg
0=1$ and $\neg 1=0$.
\end{definition}

\begin{definition}
A Boolean matrix is a matrix with values in the canonical Boolean
algebra.
\end{definition}

We define the operations sum, multiplication and transpose in the
same way as on real matrices:
$$\mbox{Sum: }
[a_{i,j}]+[b_{i,j}]:=[c_{i,j}] \mbox{ where }
c_{i,j}=a_{i,j}+b_{i,j}$$
$$\mbox{Multiplication: }
[a_{i,j}][b_{i,j}]:=[c_{i,j}] \mbox{ where }
c_{i,j}=\sum_{k=1}^{n} a_{i,k}b_{k,j}$$
$$\mbox{Transpose: }
[a_{i,j}]^t :=[a_{j,i}] $$

There is a natural partial order relation on these matrices given
in the following way: $$[a_{i,j}]\leq [b_{i,j}] \mbox{ iff }
a_{i,j} \leq b_{i,j} \forall_{i,j}$$

These matrices have many of the properties of real matrices.

\begin{proposition}
Let $A$, $B$ and $C$ be Boolean matrices with appropriate
dimensions. We have:
\begin{description}
\item[1.] (commutativity of the sum) $A+B=B+A$;
\item[2.] (associativity) $(A+B)+C=A+(B+C)$ and $(AB)C=A(BC)$;
\item[3.] (distributivity) $A(B+C)=AB+AC$ and $(A+B)C=AC+BC$;
\item[4.] (existence of the zero matrix) $A+O=A$ , $AO=O$ and
$OA=O$ where $O$ is the matrix with all entries equal to zero;
\item[5.] (existence of the identity matrix) $AI=A$ and
$IA=A$ where $I=[\delta_{i,j}]$ with
$\delta_{i,j}=1\Leftrightarrow i=j$;
\item[7.] (idempotency of the sum) $A+A=A$;
\item[8.] $A\leq B \Leftrightarrow A+B=B$;
\item[9.] $(AB)^t=B^t A^t$ and $(A+B)^t=A^t +B^t$;
\item[10.] $A\leq B \Rightarrow A+C\leq B+C$ and $CA\leq CB$
 and $AC\leq BC$ and $A^t \leq B^t$.
\end{description}
\end{proposition}

We can regard a square Boolean matrix $R=[r_{i,j}]$ of dimension
$n$ as a binary relation $\sim_R$ on the set $\{1,...,n\}$:
$$i\sim_R j \Leftrightarrow r_{i,j}=1$$ Then:

\begin{proposition}
The binary relation $\sim_R$ represented by the matrix $R$ is:
\begin{description}
 \item[i.] reflexive iff $I\leq R$;
 \item[ii.] symmetric iff $R^t=R$;
 \item[iii.] transitive iff $R^2\leq R$.
 \end{description}
\end{proposition}

Thus we can transpose the notions of reflexivity, symmetry and
transitivity from the binary relations to square Boolean matrices.
Notice that reflexivity and transitivity imply idempotency of the
product for Boolean matrices.

\begin{proposition} (definition)
Let $A$ be a square Boolean matrix and let
$$\overline{A}:=\sum_{n=1}^{\infty}A^n =A+A^2+...$$ Then:
\begin{description}
 \item[i.] $\overline{A}$ is transitive;
 \item[ii.] For any transitive matrix $B$, $A\leq B \Rightarrow
 A\leq \overline{A}\leq B$;
 \item[iii.] $A\leq B \Rightarrow \overline{A}\leq \overline{B}$;
 \item[iv.] $\overline{\overline{A}}=\overline{A}$;
 \item[v.] If $A\geq I$ then $\overline{A}=A^n$ for some natural
 $n$.\footnote{This is not true for matrices with infinite dimension.}
 \end{description}
$\overline{A}$ is called the transitive closure of $A$.
\end{proposition}

Next we will define a {\it lattice ordered additive monoid} to be
a commutative monoid $(\mathbb{M},\oplus,\emptyset)$, where
$\oplus$ is the binary operation of the monoid and $\emptyset$ is
the zero element, with a partial order relation $\leq$ such that
$(\mathbb{M},\leq)$ is a distributive lattice with minimum
$\emptyset$ and where the sum $\oplus$ is distributive over the
operations meet $\wedge$ and join $\vee$. Formally, it is a set
$\mathbb{M}$ with three binary operations $\oplus$, $\vee$ and
$\wedge$ and an element $\emptyset$ such that for any $a,b,c \in
\mathbb{M}$ :

\begin{description}
\item[M.] $(\mathbb{M},\oplus)$ is a commutative monoid:
\begin{description}
\item[M1.] (commutativity) $a\oplus b=b\oplus a$;
\item[M2.] (associativity) $(a\oplus b)\oplus c=a\oplus (b\oplus c)$;
\item[M3.] (existence of the zero element) $\emptyset\oplus a=a$;
\end{description}
\item[L.] $(\mathbb{M},\vee,\wedge)$ is a distributive lattice:
\begin{description}
\item[L1.] (idempotency) $a\vee a=a$ and
$a\wedge a=a$;
\item[L2.] (commutativity) $a\vee b=b\vee a$ and
$a\wedge b=b\wedge a$;
\item[L3.] (associativity) $(a\vee b)\vee c=a\vee (b\vee c)$ and
$(a\wedge b)\wedge c=a\wedge (b\wedge c)$;
\item[L4.] (absorption) $a\wedge (a\vee b)=a$ and $a\vee (a\wedge b)=a$;
\item[L5.] (distributivity) $a\wedge (b\vee c)=(a\wedge b)\vee (a\wedge c)$
 and $a\vee (b\wedge c)=(a\vee b)\wedge (a\vee c)$;
\end{description}
\item[C.] The lattice and monoid structures of $\mathbb{M}$ are
compatible by the following axioms:
\begin{description}
\item[C1.] $\emptyset\vee a=a$ and $\emptyset\wedge a=\emptyset$;
\item[C2.] $a\oplus (b\vee c)=(a\oplus b)\vee (a\oplus c)$
 and $a\oplus (b\wedge c)=(a\oplus b)\wedge (a\oplus c)$.
\end{description}
\end{description}

Remember that by definition $a\vee b =\sup\{a,b\}$ and $a\wedge b
=\inf\{a,b\}$. Also we have $a\leq b \Leftrightarrow a\vee b=b
\Leftrightarrow a\wedge b=a$. Using the axioms of such
monoids\footnote{In this paper we will only consider monoids of
this type and will refer to them simply as monoids.} we have the
following properties:

\begin{description}
\item[P1.] $a\leq a\oplus b$;
\item[P2.] $(a\vee b)\oplus (a\wedge b)=a\oplus b$.
\end{description}

The property P1 is very easy to prove and the proof of
 P2 follows from the following inequalities:
$$\begin{array}{ccc}
  (a\vee b)\oplus (a\wedge b) & =
  & [a\oplus (a\wedge b)]\vee [b\oplus (a\wedge b)] \\
   & = & [(a\oplus a)\wedge (a\oplus b))]\vee
   [(b\oplus a)\wedge (b\oplus b)] \\
   & \leq & (a\oplus b)\vee (b\oplus a)\\
   & = & a\oplus b
\end{array}$$

$$\begin{array}{ccc}
  (a\vee b)\oplus (a\wedge b) & =
  & [(a\vee b)\oplus a]\wedge [(a\vee b)\oplus b] \\
   & = & [(a\oplus a)\vee (b\oplus a))]\wedge
   [(a\oplus b)\vee (b\oplus b)] \\
   & \geq & (b\oplus a)\wedge (a\oplus b)\\
   & = & a\oplus b
\end{array}$$

{\bf Examples:}
 \begin{description}
 \item[1.] $\mathbb{M}:=\mathbb{N}_0 = \{0,1,2,...\}$, $\emptyset:=0$,
 $a\oplus b:=a+b$, $a\vee b:=\max\{a,b\}$ and $a\wedge b:=\min\{a,b\}$;
 \item[2.] $\mathbb{M}:=\mathbb{N}_1 = \{1,2,...\}$, $\emptyset:=1$,
 $a\oplus b:=ab$, $a\vee b:=l.c.m.\{a,b\}$ and $a\wedge b:=g.c.d.\{a,b\}$;
 \item[3.] $\mathbb{M}$ a distributive lattice with minimum and
 $\oplus:=\vee$.
 \end{description}

Now we consider the following action of the canonical Boolean
algebra $\mathcal{B}$ on a monoid $\mathbb{M}$:
$$\xymatrix{ \mathcal{B}\times\mathbb{M} \ar[r] 
&\mathbb{M} \\ (v,m)\ar @{|->}[r] &v*m}$$ where $$v*m:=\left\{
\begin{array}{c}
 m \mbox{ if }v=1  \\ \emptyset \mbox{ if }v=0
\end{array}\right. $$

Then we have:
 \begin{description}
 \item[i.] $(v_1 v_2)*m=v_1 *(v_2 *m)$
 $\forall v_1 ,v_2 \in \mathcal{B}; m\in \mathbb{M}$;
 \item[ii.] $(v_1 + v_2)*m=(v_1 *m)\vee (v_2 *m)$
 $\forall v_1 ,v_2 \in \mathcal{B}; m\in \mathbb{M}$;
 \item[iii.] $v*(m_1 \vee m_2)=(v *m_1)\vee (v *m_2)$
 $\forall v\in \mathcal{B}; m_1 ,m_2 \in \mathbb{M}$;
 \item[iv.] $v*(m_1 \oplus m_2)=(v *m_1)\oplus(v *m_2)$
 $\forall v\in \mathcal{B}; m_1 ,m_2 \in \mathbb{M}$.
 \end{description}

Now we can define an action of Boolean matrices on arrays with
values in the monoid $\mathbb{M}$.

\begin{definition}
Let $[v_{i,j}]_{m\times n}$ be a Boolean matrix and
$(a_j)_{j=1,...,n}$ be an array in $ \mathbb{M}^n$. We define
$$[v_{i,j}]*(a_j):=(b_i)_{i=1,...,m} \mbox{ where } b_i
=\bigvee_{j=1}^{n}v_{i,j}*a_j$$
\end{definition}

\begin{proposition} For any Boolean matrices $A$ and $B$ and any
arrays $\overrightarrow{x}$ and $\overrightarrow{y}$ with values
in $\mathbb{M}$, we have:
 \begin{description}
 \item[1.] $(AB)*\overrightarrow{x}=A*(B*\overrightarrow{x})$;
 \item[2.] $(A+B)*\overrightarrow{x}=
 (A*\overrightarrow{x})\vee (B*\overrightarrow{x})$;
 \item[3.] $I*\overrightarrow{x}=\overrightarrow{x}$ and
 $O*\overrightarrow{x}=\overrightarrow{\emptyset}$;
 \item[4.] $A*(\overrightarrow{x}\vee\overrightarrow{y})=
 (A*\overrightarrow{x})\vee(A*\overrightarrow{y})$;
 \item[5.] $A*(\overrightarrow{x}\oplus\overrightarrow{y})\leq
 (A*\overrightarrow{x})\oplus(A*\overrightarrow{y})$.
 \end{description}
where $I$ is the identity matrix, $O$ is the zero matrix, the
operators $\vee$ and $\oplus$ are defined coordinate by coordinate
in $\mathbb{M}^n$ and $\overrightarrow{\emptyset}=
(\emptyset,...,\emptyset)$.
\end{proposition}

\section{Representation of the category
$\mathbf{PT}$ on $\mathbf{PI}_{\mathbb{M}}$}

To each object of $\mathbf{PT}$ with cardinality $n$ we associate
the object $\mathcal{O}_{n+1}$ of $\mathbf{PI}_{\mathbb{M}}$.

The motivation is the following. An object $O$ of $\mathbf{PT}$
gives a decomposition of the real line into intervals, and each
planar tangle that ends on $O$ decomposes the strip
$\mathbb{R}\times [0,1]$ into regions whose boundaries contain
these intervals. Ordering the intervals in the natural way we will
store in a Boolean matrix the information about which intervals
are in the same region, that is, the intervals $i$ and $j$ are in
the same region if and only if the $(i,j)$ entry of the Boolean
matrix is $1$. Also to each interval we associate a value (in the
given monoid $\mathbb{M}$) which is specific for the region to
which the interval belongs. Thus in this way, intervals in the
same region have the same value and therefore the array of values
is fixed by the action of the matrix.

This should make clear the reason for the properties that the
matrices in $\mathcal{O}_n$ have to satisfy. Indeed, the author
conjectures that any matrix with the properties E1, E2, E3, T1, T2
and T3 has a geometric realization in this form.

To obtain a functor from the category $\mathbf{PT}$ to the
category $\mathbf{PI}_{\mathbb{M}}$ we need to associate to each
elementary tangle $\hat{t}_{n,k}$ and $\check{t}_{n,k}$ functions
$\hat{T}_{n,k} : \mathcal{O}_n \longrightarrow \mathcal{O}_{n+2}$
and $\check{T}_{n,k}: \mathcal{O}_{n+2} \longrightarrow
\mathcal{O}_n$ that satisfy the same relations as $\hat{t}_{n,k}$
and $\check{t}_{n,k}$.

We want these functions to preserve the motivation for the
definition of $\mathcal{O}_n$. Specifically if
$(R,\overrightarrow{v})$ is an element of $\mathcal{O}_n$, and $R$
is the matrix of connectivity of the intervals for a specific
tangle that ends on the object associated with $\mathcal{O}_n$
then the image $(R',\overrightarrow{v}')$ of
$(R,\overrightarrow{v})$ by $\hat{T}_{n,k}$ (or $\check{T}_{n,k}$)
has $R'$ as the matrix of connectivity of the intervals which
terminate the composition of the tangle $\hat{t}_{n,k}$ (or
$\check{t}_{n-2,k}$) with the specific tangle. Furthermore, if
$\overrightarrow{v}$ gives the values assigned to the intervals,
then $\overrightarrow{v}'$ gives the values assigned to the
intervals after composition with the tangle $\hat{t}_{n,k}$ (or
$\check{t}_{n-2,k}$).

We will define $\hat{T}_{n,k}$ and $\check{T}_{n,k}$ as follows
$$\hat{T}_{n,k}(R,\overrightarrow{v})=(R',\overrightarrow{v}')$$
where $$R'=B_{n,k}RB^t_{n,k}+D_{n+2,k}$$ and
$$\overrightarrow{v}'=B_{n,k}*\overrightarrow{v}$$ $B_{n,k}$ is a
Boolean matrix with $n+2$ rows and $n$ columns defined by
$$B_{n,k}:=[b_{i,j}] \mbox{ with } b_{i,j}=1 \mbox{ iff } i=j<k
\mbox{ or } i=j+2>k.$$ $D_{n,k}$ is the diagonal square Boolean
matrix of dimension $n$ defined by $$D_{n,k}:=[d_{i,j}] \mbox{
with } d_{i,j}=1 \mbox{ iff } i=j=k.$$

We can regard the matrix $B_{n,k}$ as the connectivity relation
between the upper and lower intervals of the tangle
$\hat{t}_{n,k}$, that is, $b_{i,j}=1$ iff the upper interval $j$
and the lower interval $i$ are in the same region for the tangle
$\hat{t}_{n,k}$ (or equivalently, iff the upper interval $i$ and
the lower interval $j$ are in the same region for the tangle
$\check{t}_{n,k}$).

In this sense the formula $R'=B_{n,k}RB^t_{n,k}+D_{n+2,k}$ means
that two distinct intervals $i$ and $j$ are in the same region
after the composition with $\hat{t}_{n,k}$ if $i\not=k$ and
$j\not=k$ and the intervals $\hat{k}(i)$ and $\hat{k}(j)$
($\hat{k}(i)=i$ if $i<k$, $\hat{k}(i)=i+2$ if $i>k$) are in the
same region before the composition by $\hat{t}_{n,k}$. In other
words two intervals which not $k$ are in the same region if the
respective intervals above them in the tangle $\hat{t}_{n,k}$ are
in the same region before the composition with $\hat{t}_{n,k}$.

The formula $\overrightarrow{v}'=B_{n,k}*\overrightarrow{v}$ means
that the extended regions (after the composition with the tangle
$\hat{t}_{n,k}$) preserve the old values and the new region
created over the interval $k$ receives the value $\emptyset$.


Now we define $\check{T}_{n,k}$:

$$\check{T}_{n,k}(R,\overrightarrow{v})=(R',\overrightarrow{v}')$$
where $$R'=(B^t_{n,k}RB_{n,k})^2$$ and $$\overrightarrow{v}'=
R'*[(B^t_{n,k}*\overrightarrow{v})\oplus(e_{n,k-1}*x_k)]$$ where
$e_{n,k-1}$ is a 1-column Boolean matrix of dimension $n$ defined
by $e_{n,k-1}:=[\epsilon_i]$ with $\epsilon_i=1$ iff $i=k-1$ and
 $x_k$ is a monoid value which depends
 on $R$ and$\overrightarrow{v}$ (despite this, we use the
 symbol $x_k$ instead of $x_k(R,\overrightarrow{v})$ to simplify
 the notation),
 given by the following formula:
$$x_k=[r_{k-1,k+1}*\varphi(v_k)]\oplus[(\neg
r_{k-1,k+1})*(v_{k-1}\wedge v_{k+1})]=\left\{\begin{array}{ccc}
 \varphi(v_k) &\mbox{ if }&r_{k-1,k+1}=1  \\
 v_{k-1}\wedge v_{k+1} &\mbox{ if }&r_{k-1,k+1}=0
\end{array}\right.$$
where $r_{k-1,k+1}=e^t_{n,k-1}Re_{n,k+1}$ (the $(k-1,k+1)$ entry
of $R$) and $\varphi :\mathbb{M}\longrightarrow \mathbb{M}$ is a
fixed function (without structure) independent\footnote{The
representation depends on the choice of the function $\varphi$
i.e. a different function $\varphi$ gives a different
representation.} of $(R,\overrightarrow{v})$. 


The idea behind the formula $R'=(B^t_{n,k}RB_{n,k})^2$ is the same
as before. The matrix $B^t_{n,k}RB_{n,k}$ transfers the relation
between two intervals of belonging to the same region from the top
of the tangle $\check{t}_{n,k}$ to the bottom, and we need to take
the square power because the matrix $B^t_{n,k}RB_{n,k}$ may not be
transitive, since $\check{t}_{n,k}$ joins the regions associated
to the intervals $k-1$ and $k+1$ (which may or may not be the
same).

The formula $\overrightarrow{v}'=
R'*[(B^t_{n,k}*\overrightarrow{v})\oplus(e_{n,k-1}*x_k)]$ plays a
crucial rule in the construction and needs a more careful
explanation. What it says is that the interval $k-1$ receives the
values of the old intervals $k-1$ and $k+1$ and if these intervals
are in distinct regions then we sum them by the operation $\oplus$
(since $(v_{k-1}\vee v_{k+1})\oplus(v_{k-1}\wedge v_{k+1})
=(v_{k-1}\oplus v_{k+1})$ by P2). $$\psdiag{20}{40}{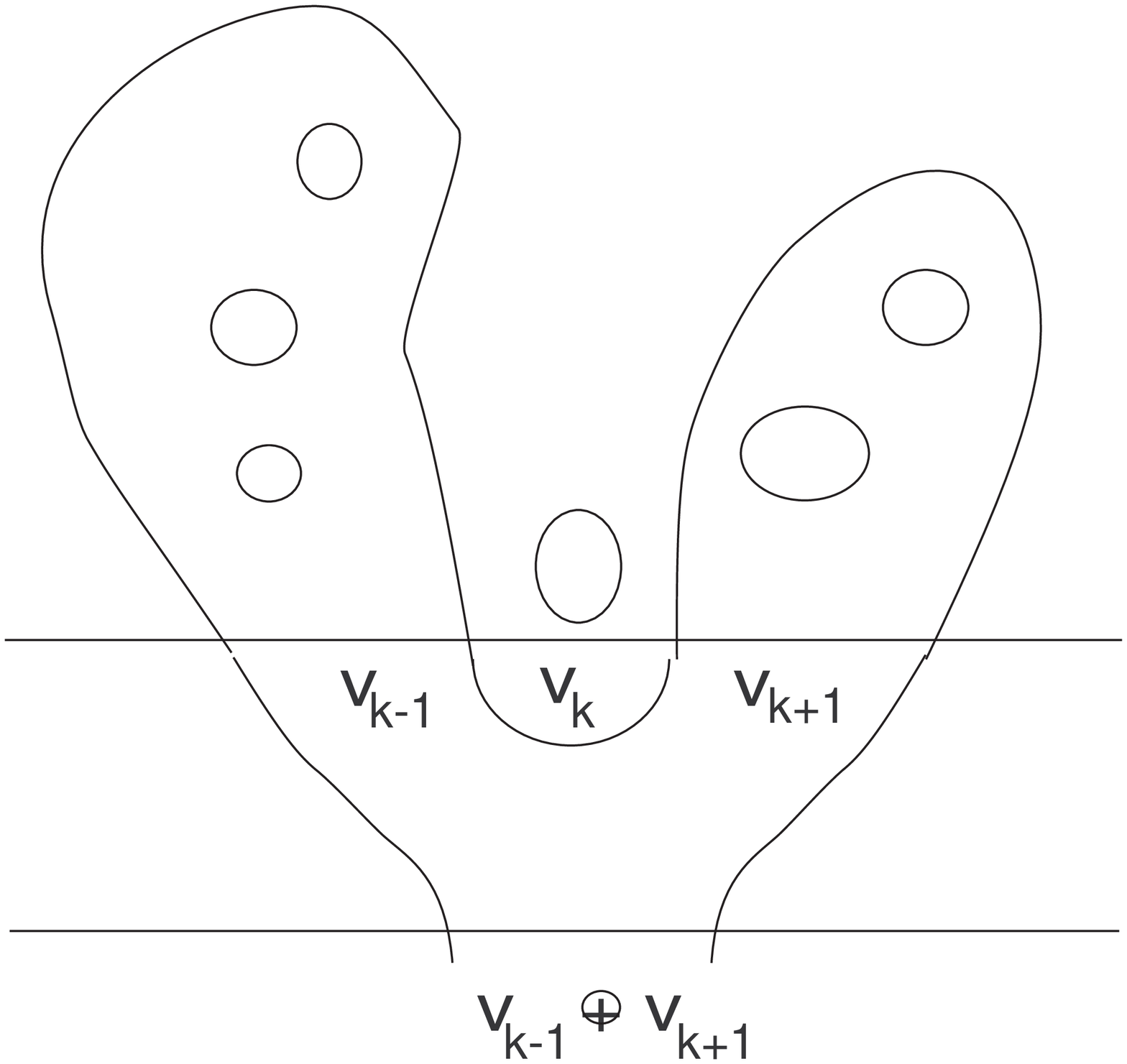}$$
 If they are in the same region then we take their common value
and sum to it some modification (given by the function $\varphi$)
of the value of the interval $k$ corresponding to a region which
is closed after the composition with $\check{t}_{n,k}$.
$$\psdiag{20}{40}{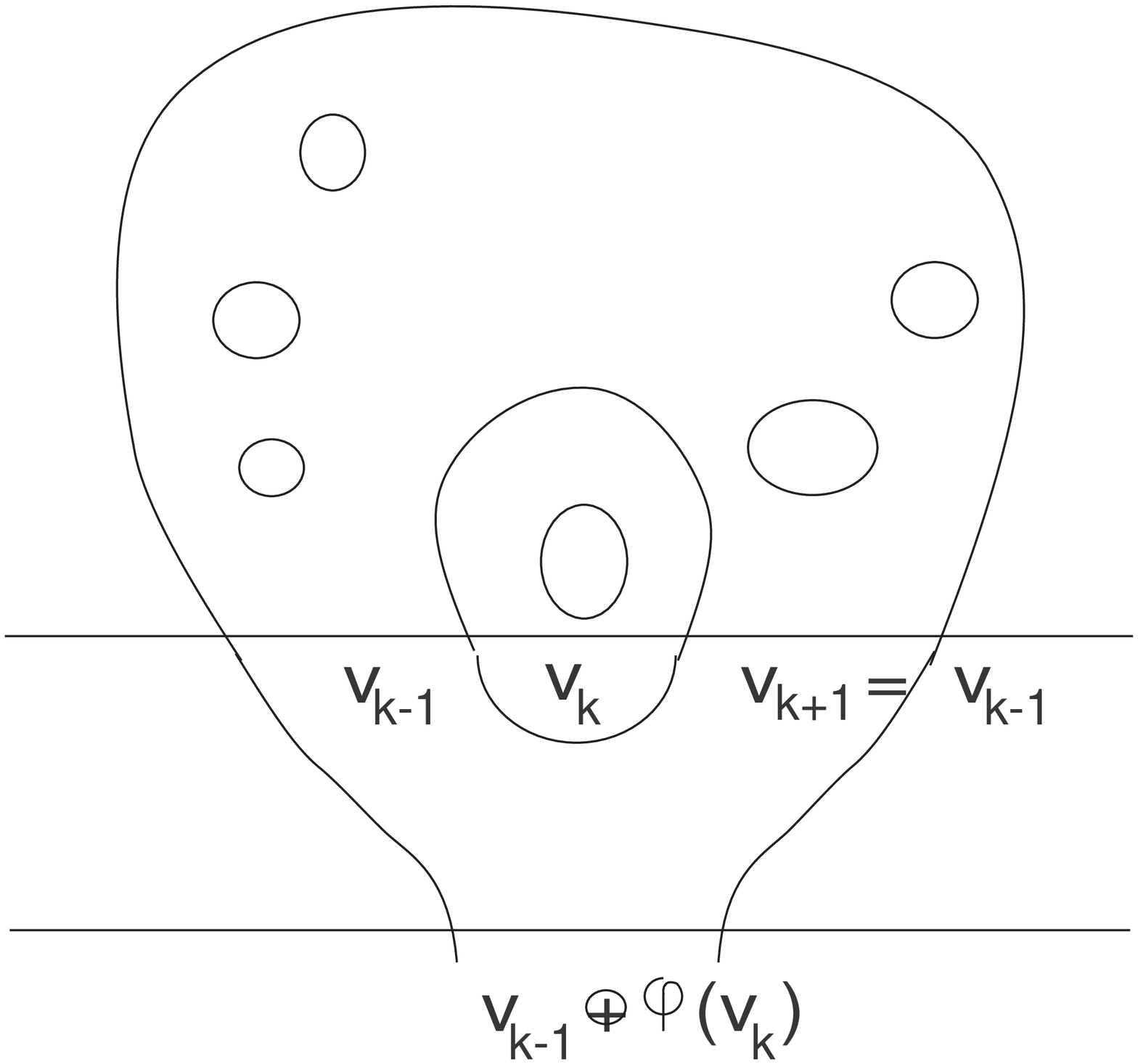}$$ The other intervals receive their former
value if they are not in the region associated to the interval
$k-1$ or receive the new value of the interval $k-1$ if they are
in the same region as that interval. This is why we take the
action of the matrix $R'$ on the array
$(B^t_{n,k}*\overrightarrow{v})\oplus(e_{n,k-1}*x_k)$ so as to
transfer the value of the interval $k-1$ to others connected with
it.

A better way of thinking about this may be that the array of
values describes the histories of the regions associated to each
interval, which keep track of the histories of any closed region
inside them by means of the function $\varphi$. The Boolean matrix
essentially plays the role of an assistant storing the information
about which intervals are in the same region. For example, in the
case of closed planar curves, which are morphisms from the empty
set to itself, we get in the end a one-dimensional square matrix
(which is unique by the condition E1) and a one-dimensional array
(or simply a monoid value). So in this case the Boolean matrix
doesn't matter at all and the only significant content is the
monoid value.


So as to simplify the notation we will always substitute
$\hat{t}_{n,k}$, $\check{t}_{n,k}$, $\hat{T}_{n,k}$,
$\check{T}_{n,k}$, $B_{n,k}$, $D_{n,k}$ and $e_{n,k}$ by
$\hat{t}_k$, $\check{t}_k$, $\hat{T}_k$, $\check{T}_k$, $B_k$,
$D_k$ and $e_k$ when $n$ is implicit.

\subsection{The well-definedness of the functions
$\hat{T}_{n,k}$ and $\check{T}_{n,k}$} $$\begin{array}{cccc}
\hat{T}_{n,k}: & \mathcal{O}_n & \longrightarrow &
\mathcal{O}_{n+2}\\ & (R,\overrightarrow{v}) & \longmapsto
&(R',\overrightarrow{v}') \end{array}$$ with $R'=B_k R B_k^t +D_k$
and $\overrightarrow{v}'= B_k *\overrightarrow{v}$.

We need to check that if $R=[r_{i,j}]$ is an $n$-dimensional
matrix that satisfies the conditions:
\begin{description}
\item[E1.] $R\geq I$;
\item[E2.] $R^t =R$;
\item[E3.] $R^2 =R$;
\item[T1.] If $r_{i,j}=1$ then $|i-j|$ is even;
\item[T2.] For any $\alpha \leq\beta \leq\gamma \leq\delta$,
$r_{\alpha,\gamma}r_{\beta,\delta} \leq
r_{\alpha,\beta}r_{\beta,\gamma}r_{\gamma,\delta}$;
\item[T3.] For any $\alpha < \beta$ if $r_{\alpha,\beta}=1$
then $r_{\alpha +1,\beta -1}=1$ or there exists $\gamma$ between
$\alpha$ and $\beta $ such that $r_{\alpha,\gamma}=1$.
\end{description}
then the matrix $R'=B_k R B_k^t +D_k$ is an $n+2$-dimensional
matrix that satisfies the same conditions.


Also we need to prove that if $\overrightarrow{v}$ is fixed by the
action of $R$:
\begin{description}
\item[EC.] $R*\overrightarrow{v}=\overrightarrow{v}$
\end{description}
then $\overrightarrow{v}'= B_k *\overrightarrow{v}$ is likewise
fixed by the action of $R'$.


First we will check that $R'$ satisfies the conditions E1, E2, E3,
T1, T2 and T3.

\begin{description}
\item[E1:] $$R\geq I \Rightarrow R'=B_k R B_k^t +D_k \geq B_k I B_k^t
+D_k \geq (I-D_k)+D_k =I$$ The operation {\it minus} on Boolean
matrices is defined in the following way:
$$[a_{i,j}]-[b_{i,j}]=[c_{i,j}] \mbox{ where } c_{i,j}=1 \mbox{
iff } a_{i,j}>b_{i,j} (\mbox{i.e. }a_{i,j}=1 \mbox{ and }
b_{i,j}=0)$$ It easy to see that, for any matrices $A$ and $B$,
$(A-B)+B\geq A$ and $A\geq B \Rightarrow (A-B)+B= A$. We leave it
to the reader to check that $ B_k B_k^t \geq I-D_k$.
\item[E2:] $$R=R^t \Rightarrow R'^t=(B_k R B_k^t +D_k)^t = B_k R^t B_k^t
+D_k^t = B_k R B_k^t +D_k =R'$$
\item[E3:] $$\begin{array}{rcl}
R^2=R \Rightarrow R'^2 &=& (B_k R B_k^t +D_k)^2 \\
 &=& B_k R B_k^t B_k R B_k^t + B_k R B_k^t D_k +D_k B_k R B_k^t +D_k^2 \\
 &=& B_k R I R B_k^t +B_k R O +O R B_k +D_k \\
 &=& B_k R B_k^t +D_k \\
 &=& R'
 \end{array}$$

We leave it to the reader to check that $B_k^t B_k =I$, $B_k^t D_k
=O$ and $D_k B_k =O$.

\item[T1:] The condition
$$r_{i,j}=1 \Rightarrow i-j\in 2\mathbb{Z}$$ is equivalent to the
condition $$R\leq C_{n\times n}$$ where $C_{m\times n}$ is the
{\it chess board matrix of dimension $m\times n$} defined in the
following way: $$C_{m\times
n}=[c_{i,j}]_{{i=1,...,m}\atop{j=1,...,n}} \mbox{ with } c_{i,j}=1
\mbox{ iff } i-j\in 2\mathbb{Z}$$

It is easy to see that $C_{l\times m}C_{m\times n}\leq C_{l\times
n}$ (in fact, this is an equality unless $m=1$), and also we have
$B_{n,k}\leq C_{(n+2)\times n}$ and $D_{n,k}\leq C_{(n+2)\times
(n+2)}$. Thus $R\leq C_{n\times n} \Rightarrow R'=B_k R B_k^t +D_k
\leq C_{(n+2)\times n}C_{n\times n}C_{n\times (n+2)}
+C_{(n+2)\times (n+2)}\leq C_{(n+2)\times (n+2)}$.

\item[T2:] Let $[r_{i,j}]=R$ and $[r'_{i,j}]=R'=B_k R B_k^t +D_k$.

Suppose that:
\begin{quote}
Hypothesis: $\forall_{\alpha\leq\beta\leq\gamma\leq\delta} \qquad
r_{\alpha,\gamma}r_{\beta,\delta} \leq
r_{\alpha,\beta}r_{\beta,\gamma}r_{\gamma,\delta}$.
\end{quote}

We want to prove that:
\begin{quote}
Thesis: $\forall_{\alpha\leq\beta\leq\gamma\leq\delta} \qquad
r'_{\alpha,\gamma}r'_{\beta,\delta} \leq
r'_{\alpha,\beta}r'_{\beta,\gamma}r'_{\gamma,\delta}$.
\end{quote}

In the case $\alpha=\beta$ or $\beta=\gamma$ or $\gamma=\delta$
this assertion is true if $[r'_{i,j}]=R'$ satisfies the conditions
of an equivalence relation (E1, E2 and E3):
\begin{quote}
 (1) $r'_{i,i}=1$ for any $i$;
\end{quote}
\begin{quote}
 (2) $r'_{i,j}=r'_{j,i}$ for any $i$ and $j$;
\end{quote}
\begin{quote}
 (3) $r'_{i,j}r'_{j,k}\leq r'_{i,k}$ for any $i$, $j$ and $k$.
\end{quote}
which we have already seen to be true.

In fact, if $\alpha=\beta$ we have
$$r'_{\alpha,\gamma}r'_{\beta,\delta}
=r'_{\alpha,\gamma}r'_{\alpha,\delta}
=r'_{\gamma,\alpha}r'_{\alpha,\delta} \leq r'_{\gamma,\delta}$$
and $$r'_{\alpha,\gamma}r'_{\beta,\delta} \leq
r'_{\alpha,\gamma}=r'_{\beta,\gamma}$$ and
$$r'_{\alpha,\gamma}r'_{\beta,\delta} \leq
1=r'_{\alpha,\alpha}=r'_{\alpha,\beta}$$ thus
$$r'_{\alpha,\gamma}r'_{\beta,\delta} \leq
r'_{\alpha,\beta}r'_{\beta,\gamma}r'_{\gamma,\delta}$$ if
$\beta=\gamma$ we have $$r'_{\alpha,\gamma}r'_{\beta,\delta}
=r'_{\alpha,\beta}r'_{\beta,\delta}
=r'_{\alpha,\beta}r'_{\beta,\beta}r'_{\beta,\delta} =
r'_{\alpha,\beta}r'_{\beta,\gamma}r'_{\gamma,\delta}$$ and if
$\gamma=\delta$ we have $$r'_{\alpha,\gamma}r'_{\beta,\delta}
=r'_{\alpha,\gamma}r'_{\beta,\gamma}
=r'_{\alpha,\gamma}r'_{\gamma,\beta} \leq r'_{\alpha,\beta}$$ and
$$r'_{\alpha,\gamma}r'_{\beta,\delta} \leq
r'_{\beta,\delta}=r'_{\beta,\gamma}$$ and
$$r'_{\alpha,\gamma}r'_{\beta,\delta} \leq
1=r'_{\gamma,\gamma}=r'_{\gamma,\delta}$$ thus
$$r'_{\alpha,\gamma}r'_{\beta,\delta} \leq
r'_{\alpha,\beta}r'_{\beta,\gamma}r'_{\gamma,\delta}$$ So we are
left with the case $\alpha<\beta<\gamma<\delta$. It easy to see
that $$ r'_{i,j}=\left\{\begin{array}{ccl}
r_{\hat{k}(i),\hat{k}(j)}&\mbox{if}&i,j\not=k\\
\delta_{k,j}&\mbox{if}&i=k\\ \delta_{i,k}&\mbox{if}&j=k
\end{array}\right.$$
$$ \mbox{ where } \hat{k}(i)=\left\{
\begin{array}{ccl}
i&\mbox{if}&i<k\\ i-2&\mbox{if}&i>k
\end{array}\right.\mbox{ and } \delta_{i,k}=\left\{
\begin{array}{ccl}
1&\mbox{if}&i=k\\ 0&\mbox{if}&i\not=k
\end{array}\right.$$
if $k\in\{\alpha,\beta,\gamma,\delta\}$ (with
$\alpha<\beta<\gamma<\delta$) then
$$r'_{\alpha,\gamma}r'_{\beta,\delta}=0\leq
r'_{\alpha,\beta}r'_{\beta,\gamma}r'_{\gamma,\delta}$$
if $k\not\in\{\alpha,\beta,\gamma,\delta\}$ we have
$\hat{k}(\alpha)\leq \hat{k}(\beta)\leq \hat{k}(\gamma)\leq
\hat{k}(\delta)$ and then $$r'_{\alpha,\gamma}r'_{\beta,\delta}=
r'_{\hat{k}(\alpha),\hat{k}(\gamma)}
r'_{\hat{k}(\beta),\hat{k}(\delta)} \leq
r'_{\hat{k}(\alpha),\hat{k}(\beta)}
r'_{\hat{k}(\beta),\hat{k}(\gamma)}
r'_{\hat{k}(\gamma),\hat{k}(\delta)}=
r'_{\alpha,\beta}r'_{\beta,\gamma}r'_{\gamma,\delta}$$

\item[T3:]

Suppose by hypothesis that:
\begin{quote}
Hypothesis: $\forall_{\alpha<\beta} \qquad r_{\alpha,\beta}=1
\Rightarrow r_{\alpha +1,\beta -1 } =1 \mbox{ or }
\exists_{\alpha<\gamma<\beta} : \quad r_{\alpha,\gamma}=1$.
\end{quote}
We want to prove
\begin{quote}
Thesis: $\forall_{\alpha<\beta} \qquad r'_{\alpha,\beta}=1
\Rightarrow r'_{\alpha +1,\beta -1 } =1 \mbox{ or }
\exists_{\alpha<\gamma<\beta} : \quad r'_{\alpha,\gamma}=1$.
\end{quote}
where $[r_{i,j}]=R$ and $[r'_{i,j}]=R'=B_k R B_k^t +D_k$.
If $\beta =\alpha +1$ then the thesis is true by the condition E2
or by the condition T1.

If $\beta =\alpha +2$ then the thesis is true by the condition E1.

Now, we consider $\beta\geq\alpha +3$. If $k=\alpha$ or $k=\beta$
we have $r'_{\alpha,\beta}=0$, so the thesis is true. If $k=\alpha
+1$ then $r'_{\alpha,k+1}=r'_{k-1,k+1}=
r_{\hat{k}(k-1),\hat{k}(k+1)}= r_{k-1,k-1}=1$ and we have
$\alpha=k-1<k+1<\beta$, so the thesis is true. If $k=\beta -1$
then $r'_{\alpha,k-1}= r'_{\alpha,\beta}r'_{\beta,k-1}=
r'_{\alpha,\beta}r'_{k+1,k-1}= r'_{\alpha,\beta}$ and we have
$\alpha<k-1<k+1=\beta$, so the thesis is true. Now suppose that
$k\not\in \{\alpha,\alpha+1,\beta-1,\beta\}$. If
$r'_{\alpha,\beta}= r_{\hat{k}(\alpha),\hat{k}(\beta)}=1$ then, by
hypothesis, $r_{\hat{k}(\alpha)+1,\hat{k}(\beta)-1}=1$ or there
exists $\hat{k}(\alpha)<\gamma'<\hat{k}(\beta)$ s.t.
$r_{\hat{k}(\alpha),\gamma'}=1$. Since $k\not\in
\{\alpha,\alpha+1,\beta-1,\beta\}$ we have that
$\hat{k}(\alpha)+1=\hat{k}(\alpha+1)$ and
$\hat{k}(\beta)-1=\hat{k}(\beta-1)$ ($k\not\in \{\alpha,\alpha+1\}
\Rightarrow \alpha+1<k \mbox{ or } k<\alpha \Rightarrow
\hat{k}(\alpha)+1=\alpha+1=\hat{k}(\alpha+1) \mbox{ or }
\hat{k}(\alpha)+1=\alpha-2+1=\hat{k}(\alpha+1)$, and the same
argument for $\hat{k}(\beta)-1=\hat{k}(\beta-1)$). Thus
$r_{\hat{k}(\alpha)+1,\hat{k}(\beta)-1}=r'_{\alpha+1,\beta-1}$.
Since $\hat{k}:\mathbb{N}\setminus\{k\} \longrightarrow\mathbb{N}$
is surjective and monotone, for any $\gamma'$ between
$\hat{k}(\alpha)$ and $\hat{k}(\beta)$ there exists $\gamma$
between $\alpha$ and $\beta$ such that $\gamma'=\hat{k}(\gamma)$.
And then $r_{\hat{k}(\alpha),\gamma'}=r'_{\alpha,\gamma}$. Thus
the thesis is true.
\end{description}


Now we only have to check the extra condition:
 \begin{description}
 \item[EC.] $R'*\overrightarrow{v}'=\overrightarrow{v}'$
 \end{description}
 assuming that $\overrightarrow{v}$ is fixed by the action of
$R$.

$$\begin{array}{rcl} R'*\overrightarrow{v}' &=& (B_k R B_k^t
+D_k)* B_{k}*\overrightarrow{v}\\
 &=& [(B_k R B_k^t
+D_k) B_{k}]*\overrightarrow{v} \\
 &=& (B_k R B_k^t B_{k}+D_k B_{k})*\overrightarrow{v}  \\
 &=& (B_k R)*\overrightarrow{v} \\
 &=& B_k *(R*\overrightarrow{v}) \\
 &=& B_k *\overrightarrow{v}\\
 &=& \overrightarrow{v}'
 \end{array}$$

Next we show that the function $\check{T}_{k}$ is well-defined.
Recall that:

$$\begin{array}{cccc} \check{T}_{n,k}: & \mathcal{O}_{n+2} &
\longrightarrow & \mathcal{O}_n\\ & (R,\overrightarrow{v}) &
\longmapsto &(R',\overrightarrow{v}')
\end{array}$$ with $$R'=(B^t_kRB_k)^2$$ and $$\overrightarrow{v}'=
R'*[(B^t_k*\overrightarrow{v})\oplus(e_{k-1}*x_k)]$$ where
$$x_k=[r_{k-1,k+1}*\varphi(v_k)]\oplus[(\neg
r_{k-1,k+1})*(v_{k-1}\wedge v_{k+1})]=\left\{\begin{array}{ccc}
 \varphi(v_k) &\mbox{ if }&r_{k-1,k+1}=1  \\
 v_{k-1}\wedge v_{k+1} &\mbox{ if }&r_{k-1,k+1}=0
\end{array}\right.$$

Let us prove that $R'$ satisfies the six conditions in
$\mathcal{O}_n$ and that $\overrightarrow{v}'$ is fixed by the
action of $R'$.
\begin{description}
\item[E1:] $$R\geq I \Rightarrow R'=(B_k^t R B_k)^2 \geq (B_k^t I
B_k)^2 =I^2 =I$$
\item[E2:] $$R=R^t \Rightarrow R'^t=[(B_k^t R B_k)^2]^t =
[(B_k^t R B_k)^t]^2 = (B_k^t R^t B_k)^2 =(B_k^t R B_k)^2 =R'$$
\item[E3:] To check the transitivity of $R'$ it is sufficient to
prove that $(B_k^t R B_k)^3\leq (B_k^t R B_k)^2$ (assuming that
$R^2=R$).
 $$\begin{array}{rcl}
(B_k^t R B_k)^3 &=& B_k^t R B_k B_k^t R B_k B_k^t R B_k \\
 &\leq & B_k^t R (I+B_k D_{k-1} B_k^t) R B_k B_k^t R B_k \\
 &=& B_k^t RIRB_k B_k^t RB_k +B_k^t RB_k D_{k-1}B_k^t RB_k B_k^t RB_k \\
 &\leq & B_k^t R B_k B_k^t R B_k +
 B_k^t RB_k D_{k-1}B_k^t R(I+B_k D_{k-1} B_k^t) RB_k\\
 &=& B_k^t R B_k B_k^t R B_k +
 B_k^t RB_k D_{k-1}B_k^t R^2 B_k \\ & & +
 B_k^t RB_k D_{k-1}B_k^t R B_k D_{k-1} B_k^t RB_k\\
 &\leq & B_k^t RB_k B_k^t R B_k +B_k^t RB_k D_{k-1}B_k^t R B_k
 +B_k^t RB_k D_{k-1}B_k^t R B_k\\
 &=& B_k^t RB_k B_k^t R B_k\\
 &=& (B_k^t R B_k)^2
 \end{array}$$

We leave it to the reader to check that $B_k B_k^t \leq I+B_k
D_{k-1} B_k^t$ and $D_{k-1}B_k^t R B_k D_{k-1}\leq D_{k-1}$.

\item[T1:] $$\begin{array}{ccl}
R\leq C_{(n+2)\times (n+2)} \Rightarrow R' &=& (B_k^t R B_k)^2
\leq (C_{n\times (n+2)}C_{(n+2)\times (n+2)}C_{(n+2)\times n})^2
\\ &\leq& (C_{n\times n})^2 =C_{n\times n}
\end{array}$$.

\item[T2:] Using the "equality" $r_{i,j}=e_i^tRe_j$ we need to
check that:

\begin{quote}
Thesis: $\forall_{\alpha\leq\beta\leq\gamma\leq\delta} \qquad
e_{\alpha}^tR'e_{\gamma}e_{\beta}^tR'e_{\delta} \leq
e_{\alpha}^tR'e_{\beta}e_{\beta}^tR'e_{\gamma}e_{\gamma}^tR'e_{\delta}$.
\end{quote}
assuming the hypothesis:

\begin{quote}
Hypothesis: $\forall_{\alpha\leq\beta\leq\gamma\leq\delta} \qquad
e_{\alpha}^tRe_{\gamma}e_{\beta}^tRe_{\delta} \leq
e_{\alpha}^tRe_{\beta}e_{\beta}^tRe_{\gamma}e_{\gamma}^tRe_{\delta}$.
\end{quote}

Since $R'$ satisfies the equivalence relation conditions (E1, E2
and E3) the thesis is satisfied for $\alpha=\beta$ or
$\beta=\gamma$ or $\gamma=\delta$. So we can assume that
$\alpha<\beta<\gamma<\delta$.

We will substitute the hypothesis by a more appropriate
hypothesis. But for that we need to introduce a new definition and
same properties. Let $u$ and $v$ be two one-column non-zero
matrices. We define: $$u\prec v \mbox{ iff } \max\{i: u_i=1
\}<\min\{i: v_i=1 \}.$$

\begin{proposition}
\begin{description}
\item[1.] $\prec$ defines a strict order relation on the set of
non-zero one-column matrices;
\item[2.] $u\prec v , e_{\alpha}\leq u \mbox{ and } e_{\beta}\leq
u \Rightarrow \alpha<\beta$;
\item[3.] $(\forall_{\alpha<\beta<\gamma<\delta} \quad
e_{\alpha}^tRe_{\gamma}e_{\beta}^tRe_{\delta} \leq
e_{\alpha}^tRe_{\beta}e_{\beta}^tRe_{\gamma}e_{\gamma}^tRe_{\delta})
\Rightarrow (\forall_{u\prec v\prec w\prec x} \quad u^tR w v^tRx
\leq u^tRvv^tRww^tRx)$;
\item[4.] $\alpha<\beta \Rightarrow B_k e_{\alpha}\prec B_k
e_{\beta}$.
\end{description}
\end{proposition}

Now we take a new (weaker) hypothesis:
\begin{quote}
N.H.: $\forall_{\alpha<\beta<\gamma<\delta} \quad
e_{\alpha}^tB_k^tRB_ke_{\gamma}e_{\beta}^tB_k^tRB_ke_{\delta} \leq
e_{\alpha}^tB_k^tRB_ke_{\beta}e_{\beta}^tB_k^tRB_ke_{\gamma}
e_{\gamma}^tB_k^tRB_ke_{\delta}$.
\end{quote}
$$
\begin{array}{rcl}
e_{\alpha}^tR'e_{\gamma}e_{\beta}^tR'e_{\delta} &=&
e_{\alpha}^tB_k^tRB_k B_k^tRB_k e_{\gamma}e_{\beta}^t B_k^tRB_k
B_k^tRB_ke_{\delta}\\
 &\leq & e_{\alpha}^tB_k^tR(I+A_k)RB_k e_{\gamma}
 e_{\beta}^t B_k^tR(I+A_k)RB_ke_{\delta}\\
 &=& e_{\alpha}^tB_k^tRB_ke_{\gamma}e_{\beta}^tB_k^tRB_ke_{\delta}
 +e_{\alpha}^tB_k^tRA_kRB_ke_{\gamma}e_{\beta}^tB_k^tRB_ke_{\delta}\\
 & &+e_{\alpha}^tB_k^tRB_ke_{\gamma}e_{\beta}^tB_k^tRA_kRB_ke_{\delta}
 +e_{\alpha}^tB_k^tRA_kRB_ke_{\gamma}e_{\beta}^tB_k^tRA_kRB_ke_{\delta}
\end{array}
$$ where $$A_k= B_kD_{k-1}B_k^t = B_ke_{k-1}e_{k-1}^tB_k^t$$

Now, we only need to prove that:
\begin{description}
\item[(i)]
$e_{\alpha}^tB_k^tRB_ke_{\gamma}e_{\beta}^tB_k^tRB_ke_{\delta}
\leq
e_{\alpha}^tR'e_{\beta}e_{\beta}^tR'e_{\gamma}e_{\gamma}^tR'e_{\delta}$;
\item[(ii)]
$e_{\alpha}^tB_k^tRA_kRB_ke_{\gamma}e_{\beta}^tB_k^tRB_ke_{\delta}
\leq
e_{\alpha}^tR'e_{\beta}e_{\beta}^tR'e_{\gamma}e_{\gamma}^tR'e_{\delta}$;
\item[(iii)]
$e_{\alpha}^tB_k^tRB_ke_{\gamma}e_{\beta}^tB_k^tRA_kRB_ke_{\delta}
\leq
e_{\alpha}^tR'e_{\beta}e_{\beta}^tR'e_{\gamma}e_{\gamma}^tR'e_{\delta}$;
\item[(iv)]
$e_{\alpha}^tB_k^tRA_kRB_ke_{\gamma}e_{\beta}^tB_k^tRA_kRB_ke_{\delta}
\leq
e_{\alpha}^tR'e_{\beta}e_{\beta}^tR'e_{\gamma}e_{\gamma}^tR'e_{\delta}$.
\end{description}

For that it is useful to observe that, for arbitrary square
matrices,
\begin{quote}
$e_i^tXe_j e_k^tYe_l = e_k^tYe_l e_i^tXe_j$,
\end{quote}
\begin{quote}
$e_i^tXe_j = e_j^tX^te_i$
\end{quote}
and
\begin{quote}
$(e_i^tXe_j)^2 = e_i^tXe_j$,
\end{quote}
since $e_i^tXe_j$ and $e_k^tYe_l$ are one-dimensional square
matrices.
\begin{description}
\item[(i)]
$$
\begin{array}{rcl}
e_{\alpha}^tB_k^tRB_ke_{\gamma}e_{\beta}^tB_k^tRB_ke_{\delta}
&\leq &
e_{\alpha}^tB_k^tRB_ke_{\beta}e_{\beta}^tB_k^tRB_ke_{\gamma}
e_{\gamma}^tB_k^tRB_ke_{\delta}\\ &\leq &
e_{\alpha}^tR'e_{\beta}e_{\beta}^tR'e_{\gamma}e_{\gamma}^tR'e_{\delta}
\end{array}
$$

Observe that $$B_k^tRB_k\leq (B_k^tRB_k)^2=R' \Rightarrow
e_i^tB_k^tRB_ke_j \leq e_i^tR'e_j$$
\item[(ii)]
$e_{\alpha}^tB_k^tRA_kRB_ke_{\gamma}e_{\beta}^tB_k^tRB_ke_{\delta}=
e_{\alpha}^tB_k^tRB_ke_{k-1}e_{k-1}^tB_k^tRB_ke_{\gamma}
e_{\beta}^tB_k^tRB_ke_{\delta}$.

If $k-1<\beta$ then $$
\begin{array}{l}
e_{\alpha}^tB_k^tRB_ke_{k-1}e_{k-1}^tB_k^tRB_ke_{\gamma}
e_{\beta}^tB_k^tRB_ke_{\delta}\\ \leq
e_{\alpha}^tB_k^tRB_ke_{k-1}e_{k-1}^tB_k^tRB_ke_{\beta}
e_{\beta}^tB_k^tRB_ke_{\gamma} e_{\gamma}^tB_k^tRB_ke_{\delta}\\
\leq
e_{\alpha}^tR'e_{\beta}e_{\beta}^tR'e_{\gamma}e_{\gamma}^tR'e_{\delta}
\end{array}
$$

If $k-1=\beta$ then $$
\begin{array}{l}
e_{\alpha}^tB_k^tRB_ke_{k-1}e_{k-1}^tB_k^tRB_ke_{\gamma}
e_{\beta}^tB_k^tRB_ke_{\delta} \\=
e_{\alpha}^tB_k^tRB_ke_{\beta}e_{\beta}^tB_k^tRB_ke_{\gamma}
e_{\beta}^tB_k^tRB_ke_{\delta}\\ \leq
e_{\alpha}^tR'e_{\beta}e_{\beta}^tR'e_{\gamma}e_{\beta}^tR'e_{\delta}\\
= e_{\alpha}^tR'e_{\beta}e_{\beta}^tR'e_{\gamma}
e_{\beta}^tR'e_{\gamma} e_{\beta}^tR'e_{\delta}\\ =
e_{\alpha}^tR'e_{\beta}e_{\beta}^tR'e_{\gamma}
e_{\gamma}^tR'e_{\beta}e_{\beta}^tR'e_{\delta}\\ \leq
e_{\alpha}^tR'e_{\beta}e_{\beta}^tR'e_{\gamma}e_{\gamma}^tR'e_{\delta}
\end{array}
$$

If $\beta<k-1<\delta$ then $$
\begin{array}{l}
e_{\alpha}^tB_k^tRB_ke_{k-1}e_{k-1}^tB_k^tRB_ke_{\gamma}
e_{\beta}^tB_k^tRB_ke_{\delta} \\=
e_{\alpha}^tB_k^tRB_ke_{k-1}e_{\beta}^tB_k^tRB_ke_{\delta}
e_{k-1}^tB_k^tRB_ke_{\gamma}
 \\ \leq
e_{\alpha}^tB_k^tRB_ke_{\beta}e_{\beta}^tB_k^tRB_ke_{k-1}
e_{k-1}^tB_k^tRB_ke_{\delta} e_{k-1}^tB_k^tRB_ke_{\gamma}\\ \leq
 e_{\alpha}^tR'e_{\beta}e_{\beta}^tR'e_{k-1}
e_{k-1}^tR'e_{\delta}e_{k-1}^tR'e_{\gamma}\\ =
e_{\alpha}^tR'e_{\beta}e_{\beta}^tR'e_{k-1}
e_{k-1}^tR'e_{\gamma}e_{\gamma}^tR'e_{k-1}e_{k-1}^tR'e_{\delta}\\
\leq
e_{\alpha}^tR'e_{\beta}e_{\beta}^tR'e_{\gamma}e_{\gamma}^tR'e_{\delta}
\end{array}
$$

If $k-1=\delta$ then $$
\begin{array}{l}
e_{\alpha}^tB_k^tRB_ke_{k-1}e_{k-1}^tB_k^tRB_ke_{\gamma}
e_{\beta}^tB_k^tRB_ke_{\delta} \\=
e_{\alpha}^tB_k^tRB_ke_{\delta}e_{\delta}^tB_k^tRB_ke_{\gamma}
e_{\beta}^tB_k^tRB_ke_{\delta}\\ =
e_{\alpha}^tB_k^tRB_ke_{\delta}(e_{\delta}^tB_k^tRB_ke_{\beta})^2
e_{\gamma}^tB_k^tRB_ke_{\delta}\\ \leq
e_{\alpha}^tR'e_{\beta}e_{\beta}^tR'e_{\delta}
(e_{\gamma}^tR'e_{\delta})^2 \\ =
e_{\alpha}^tR'e_{\beta}e_{\beta}^tR'e_{\delta}
e_{\delta}^tR'e_{\gamma} e_{\gamma}^tR'e_{\delta}\\  \leq
e_{\alpha}^tR'e_{\beta}e_{\beta}^tR'e_{\gamma}e_{\gamma}^tR'e_{\delta}
\end{array}
$$

If $k-1>\delta$ then $$
\begin{array}{l}
e_{\alpha}^tB_k^tRB_ke_{k-1}e_{k-1}^tB_k^tRB_ke_{\gamma}
e_{\beta}^tB_k^tRB_ke_{\delta} \\=
e_{\alpha}^tB_k^tRB_ke_{k-1}e_{\beta}^tB_k^tRB_ke_{\delta}
e_{\gamma}^tB_k^tRB_ke_{k-1}
 \\ \leq
e_{\alpha}^tB_k^tRB_ke_{k-1}e_{\beta}^tB_k^tRB_ke_{\gamma}
e_{\gamma}^tB_k^tRB_ke_{\delta} e_{\delta}^tB_k^tRB_ke_{k-1}\\
\leq  e_{\alpha}^tR'e_{k-1}e_{\beta}^tR'e_{\gamma}
e_{\gamma}^tR'e_{\delta}e_{\delta}^tR'e_{k-1}\\ =
e_{\alpha}^tR'e_{k-1}e_{k-1}^tR'e_{\delta}e_{\beta}^tR'e_{\gamma}
e_{\gamma}^tR'e_{\delta}\\ \leq
e_{\alpha}^tR'e_{\delta}e_{\beta}^tR'e_{\gamma}e_{\gamma}^tR'e_{\delta}
\\ =
e_{\alpha}^tR'e_{\beta}e_{\beta}^tR'e_{\gamma}e_{\gamma}^tR'e_{\delta}
\end{array}
$$
\item[(iii)]
$e_{\alpha}^tB_k^tRB_ke_{\gamma}e_{\beta}^tB_k^tRA_kRB_ke_{\delta}
=e_{\delta}^tB_k^tRA_kRB_ke_{\beta}e_{\gamma}^tB_k^tRB_ke_{\alpha}$.

Since all arguments in case (ii) are valid reversing the order of
$\alpha,\beta,\gamma$ and $\delta$, we have that:
$e_{\alpha}^tB_k^tRB_ke_{\gamma}e_{\beta}^tB_k^tRA_kRB_ke_{\delta}
\leq
e_{\alpha}^tR'e_{\beta}e_{\beta}^tR'e_{\gamma}e_{\gamma}^tR'e_{\delta}$.
\item[(iv)]
$$
\begin{array}{l}
e_{\alpha}^tB_k^tRA_kRB_ke_{\gamma}e_{\beta}^tB_k^tRA_kRB_ke_{\delta}\\
=e_{\alpha}^tB_k^tRB_ke_{k-1}e_{k-1}^tB_k^tRB_ke_{\gamma}
e_{\beta}^tB_k^tRB_ke_{k-1}e_{k-1}^tB_k^tRB_ke_{\delta} \\ \leq
e_{\alpha}^tR'e_{k-1}e_{k-1}^tR'e_{\gamma}
e_{\beta}^tR'e_{k-1}e_{k-1}^tR'e_{\delta}\\
=e_{\alpha}^tR'e_{k-1}e_{k-1}^tR'e_{\beta}e_{\beta}^tR'e_{k-1}
e_{k-1}^tR'e_{\gamma}e_{\gamma}^tR'e_{k-1}e_{k-1}^tR'e_{\delta}\\
\leq
e_{\alpha}^tR'e_{\beta}e_{\beta}^tR'e_{\gamma}e_{\gamma}^tR'e_{\delta}
\end{array}
$$
\end{description}
\item[T3:]

\begin{lemma}
If $R=[r_{i,j}]$ represents an equivalence relation, then the
following statements are equivalent:
\begin{description}
\item[(i)] $\forall_{\alpha<\beta} \quad r_{\alpha,\beta}=1
\Rightarrow r_{\alpha +1,\beta -1 } =1 \mbox{ or }
\exists_{\alpha<\gamma<\beta} : \quad r_{\alpha,\gamma}=1$;
\item[(ii)] $\forall_{\alpha<\beta} \quad r_{\alpha,\beta}=1
\Rightarrow \exists_{\alpha<\gamma\leq \beta} : \quad
r_{\alpha,\gamma}=r_{\alpha+1,\gamma-1}=1$;
\item[(iii)] $\forall_{\alpha<\beta} \quad r_{\alpha,\beta}=1
\Rightarrow r_{\alpha +1,\beta -1 } =1 \mbox{ or }
\exists_{\alpha<\gamma<\beta} : \quad r_{\gamma,\beta}=1$;
\item[(iv)] $\forall_{\alpha<\beta} \quad r_{\alpha,\beta}=1
\Rightarrow \exists_{\alpha\leq \gamma< \beta} : \quad
r_{\gamma,\beta}=r_{\gamma+1,\beta-1}=1$.
\end{description}
\end{lemma}
\TeXButton{Proof}{\proof} It easy to see that $(ii)\Rightarrow
(i)$ and $(iv)\Rightarrow (iii)$. To see that $(i)\Rightarrow
(ii)$ we take $\gamma=\inf\{\delta\leq \beta :
r_{\alpha,\delta}=1\}$ and to see that $(iii)\Rightarrow (iv)$ we
take $\gamma=\sup\{\delta\geq \alpha : r_{\delta,\beta}=1\}$.
$(i)\Leftrightarrow(iii)$ results from the transitivity of
$[r_{i,j}]$. \TeXButton{End Proof}{\endproof}

Now we want to see that $[r'_{i,j}]=(B_k^tRB_k)^2$ satisfies one
of the statements of the lemma (assuming that $R=[r_{i,j}]$ also
satisfies the same statements and the remaining conditions on
$\mathcal{O}_{n+2}$). Using the properties of $R=[r_{i,j}]$ and
the relation $R'=[r'_{i,j}]=(B_k^tRB_k)^2$ we can set:
$$r'_{\alpha,\beta}=\left\{ \begin{array}{lll}
r_{\alpha,\beta}+r_{\alpha,k+1}r_{k-1,\beta} &\mbox{if}&
k-1>\beta\\ r_{\alpha,k-1}+r_{\alpha,k+1} &\mbox{if}&
k-1=\beta\\r_{\alpha,\beta+2}+r_{\alpha,k-1}r_{k+1,\beta+2}
&\mbox{if}& \alpha<k-1<\beta\\r_{k-1,\beta+2}+r_{k+1,\beta+2}
&\mbox{if}&
k-1=\alpha\\r_{\alpha+2,\beta+2}+r_{\alpha+2,k+1}r_{k-1,\beta+2}
&\mbox{if}& k-1<\alpha
\end{array}\right.
$$

\begin{description}
\item[Case A:] $k-1>\beta$.
$r'_{\alpha,\beta}=r_{\alpha,\beta}+r_{\alpha,k+1}r_{k-1,\beta}=1
\Rightarrow r_{\alpha,\beta}=1 \mbox{ or }
r_{\alpha,k+1}=r_{k-1,\beta}=1$.
\begin{description}
\item[A.1:] $r_{\alpha,\beta}=1 \Rightarrow
\exists_{\alpha<\gamma_1 \leq\beta}: r_{\alpha+1,\gamma_1
-1}=r_{\alpha,\gamma_1}=1 \Rightarrow \exists_{\alpha<\gamma_1
\leq\beta}: r'_{\alpha+1,\gamma_1 -1}=r'_{\alpha,\gamma_1}=1$.
\item[A.2:] $r_{\alpha,k+1}=r_{k-1,\beta}=1 \Rightarrow
\exists_{\alpha<\gamma_2 \leq k+1}: r_{\alpha+1,\gamma_2
-1}=r_{\alpha,\gamma_2}=1$.
\begin{description}
\item[A.2.1:] If $\gamma_2\leq \beta$ then we have $r'_{\alpha+1,\gamma_2
-1}=r'_{\alpha,\gamma_2}=1$ with $\alpha<\gamma_2\leq\beta$.
\item[A.2.2:] If $\beta<\gamma_2\leq k-1$ we can use the
condition T2 to get $r_{\alpha,\beta}=1$ (case A.1) since
$r_{\alpha,\gamma_2}=1$ and $r_{\beta,k-1}=1$.
\item[A.2.3:] If $\gamma_2= k$ we have $r_{\alpha,k}=1$ which
together with $r_{\alpha,k+1}=1$ contradicts the condition T1.
\item[A.2.4:] If $\gamma_2=k+1$ then $r_{\alpha+1,k}=1 \Rightarrow
\exists_{\alpha+1\leq \gamma_3<k}: r_{\gamma_3+1,k-1}=
r_{\gamma_3,k}=1$.
\begin{description}
\item[A.2.4.1:] If $\beta\leq\gamma_3\leq k-1$ we can use the
condition T2 to get $r_{k-1,k}=1$ (which contradicts the condition
T1) since $r_{\beta,k-1}=1$ and $r_{\gamma_3,k}=1$.
\item[A.2.4.2:] If $\alpha+1\leq\gamma_3\leq\beta-1$ we have
$r'_{\gamma_3+1,\alpha}=1$ (since $r_{\gamma_3+1,k-1}=1$ and
$r_{k+1,\alpha}=1$) and $r_{\gamma_3,\alpha+1}=1$ (since
$r_{\gamma_3,k}=1$ and $r_{k,\alpha+1}=1$). Taking
$\gamma_4=\gamma_3+1$ we have $r'_{\alpha+1,\gamma_4
-1}=r'_{\alpha,\gamma_4}=1$ with $\alpha<\gamma_4\leq\beta$.
\end{description}
\end{description}
\end{description}
\item[Case B:] $k-1=\beta$.
$r'_{\alpha,\beta}=r_{\alpha,k-1}+r_{\alpha,k+1}=1 \Rightarrow
r_{\alpha,k-1}=1 \mbox{ or } r_{\alpha,k+1}=1$.
\begin{description}
\item[B.1:] $r_{\alpha,k-1}=1 \Rightarrow
\exists_{\alpha<\gamma_1 \leq k-1}: r_{\alpha+1,\gamma_1
-1}=r_{\alpha,\gamma_1}=1 \Rightarrow \exists_{\alpha<\gamma_1
\leq k-1}: r'_{\alpha+1,\gamma_1 -1}=r'_{\alpha,\gamma_1}=1$.
\item[B.2:] $r_{\alpha,k+1}=1 \Rightarrow
\exists_{\alpha<\gamma_2 \leq k+1}: r_{\alpha+1,\gamma_2
-1}=r_{\alpha,\gamma_2}=1$.
\begin{description}
\item[B.2.1:] If $\gamma_2\leq k-1$ then we have $r'_{\alpha+1,\gamma_2
-1}=r'_{\alpha,\gamma_2}=1$ with $\alpha<\gamma_2\leq k-1=\beta$.
\item[B.2.2:] If $\gamma_2=k$ we have $r_{\alpha,k}=1$ which
together with $r_{\alpha,k+1}=1$ contradicts the condition T1.
\item[B.2.3:] If $\gamma_2=k+1$ then $r_{\alpha+1,k}=1 \Rightarrow
\exists_{\alpha+1\leq \gamma_3<k}: r_{\gamma_3+1,k-1}=
r_{\gamma_3,k}=1$.
\begin{description}
\item[B.2.3.1:] If $\gamma_3= k-1$ then $r_{\gamma_3+1,k-1}=
r_{k,k-1}=1$ which contradicts the condition T1.
\item[B.2.3.2:] If $\alpha+1\leq\gamma_3\leq k-2$ we have
$r'_{\gamma_3+1,\alpha}=1$ (since $r_{\gamma_3+1,k-1}=1$ and
$r_{k+1,\alpha}=1$) and $r_{\gamma_3,\alpha+1}=1$ (since
$r_{\gamma_3,k}=1$ and $r_{k,\alpha+1}=1$). Taking
$\gamma_4=\gamma_3+1$ we have $r'_{\alpha+1,\gamma_4
-1}=r'_{\alpha,\gamma_4}=1$ with $\alpha<\gamma_4\leq k-1$.
\end{description}
\end{description}
\end{description}
\item[Case C:] $\alpha<k-1<\beta$.
$r'_{\alpha,\beta}=r_{\alpha,\beta+2}+r_{\alpha,k+1}r_{k-1,\beta+2}=1
\Rightarrow r_{\alpha,\beta+2}=1 \mbox{ or }
r_{\alpha,k+1}=r_{k-1,\beta+2}=1$.
\begin{description}
\item[C.1:] $r_{\alpha,\beta+2}=1 \Rightarrow
\exists_{\alpha<\gamma_1 \leq\beta+2}: r_{\alpha+1,\gamma_1
-1}=r_{\alpha,\gamma_1}=1 \mbox{ and } \exists_{\alpha\leq
\gamma_2 <\beta+2}: r_{\gamma_2,\beta+2}=
r_{\gamma_2+1,\beta+1}=1$.
\begin{description}
\item[C.1.1:] If $\gamma_1\geq \gamma_2+1$ we can use the
condition T2 to get $r_{\alpha,\beta+1}=1$ since
$r_{\alpha,\gamma_1}=1$ and $r_{\gamma_2+1,\beta+1}=1$. This
together with $r_{\alpha,\beta+2}=1$ contradicts the condition T1.
\item[C.1.2:] If $\gamma_1\leq\gamma_2$ then $\gamma_1<k$ or $\gamma_2>k$
or  $\gamma_1=\gamma_2=k$.
\begin{description}
\item[C.1.2.1:] If $\gamma_1< k$ then $r_{\alpha+1,\gamma_1
-1}=r_{\alpha,\gamma_1}=1 \Rightarrow r'_{\alpha+1,\gamma_1
-1}=r'_{\alpha,\gamma_1}=1$ with $\alpha<\gamma_1<k\leq \beta$.
\item[C.1.2.2:] If $\gamma_2> k$ then $r_{\gamma_2,\beta+2}=
r_{\gamma_2+1,\beta+1}=1 \Rightarrow r'_{\gamma_2,\beta+2}=
r'_{\gamma_2+1,\beta+1}=1$ with $\alpha\leq k-2<\gamma_2-2<\beta$.
\item[C.1.2.3:] If $\gamma_1=\gamma_2=k$ then $r_{\alpha+1,\gamma_1-1}=
r_{\alpha+1,k-1}=1$ and $r_{\gamma_2+1,\beta+1}=
r_{k+1,\beta+1}=1$ which implies $r'_{\alpha+1,\beta-1}=1$.
\end{description}
\end{description}
\item[C.2:] $r_{\alpha,k-1}=r_{k+1,\beta+2}=1 \Rightarrow
r'_{\alpha,k-1}=1$ with $\alpha<k-1<\beta$.
\end{description}
\item[Case D:] $k-1=\alpha$. This case is analogous to case B.
\item[Case E:] $k-1<\alpha$. This case is analogous to case A.
\end{description}
\end{description}


Now we only have to check the condition:
\begin{description}
\item[EC.] $R'*\overrightarrow{v}'=\overrightarrow{v}'$
\end{description}
which is obvious from the definition of $\overrightarrow{v}'$
because $R'$ is idempotent.

\subsection{Checking the relations on $\hat{T}_k$ and
$\check{T}_k$}

Now we are going to prove that $\hat{T}_k$ and $\check{T}_k$
satisfy the following relations:
\begin{description}
\item[1.] $\check{T}_{k+1}\circ\hat{T}_k=
\check{T}_{k-1}\circ\hat{T}_k= id$;
\item[2.] $\hat{T}_l\circ\hat{T}_k=
\hat{T}_k\circ\hat{T}_{l-2}$ for $l\geq k+2$;
\item[3.] $\check{T}_k\circ\hat{T}_l=
\hat{T}_{l-2}\circ\check{T}_k$ and $\check{T}_l\circ\hat{T}_k=
\hat{T}_k\circ\check{T}_{l-2}$ for $l\geq k+2$;
\item[4.] $\check{T}_k\circ\check{T}_l=
\check{T}_{l-2}\circ\check{T}_k$ for $l\geq k+2$.
\end{description}

We begin by checking the first relation.

\begin{description}
\item[1.] $\check{T}_{k+1}\circ\hat{T}_k=
\check{T}_{k-1}\circ\hat{T}_k= id$.

Let $(R_1,\overrightarrow{a})\in \mathcal{O}_n$,
$(R_2,\overrightarrow{b})=\hat{T}_k(R_1,\overrightarrow{a})$ and
$(R_3,\overrightarrow{c})=\check{T}_{k+1}(R_2,\overrightarrow{b})=
\check{T}_{k+1}\circ\hat{T}_k(R_1,\overrightarrow{a})$.

We want to show that
$(R_3,\overrightarrow{c})=(R_1,\overrightarrow{a})$.
$$\left\{\begin{array}{lll} R_2 &=&B_kR_1B_k^t+D_k=B_kR_1B_k^t+I\\
R_3 &=&(B_{k+1}^tR_2B_{k+1})^2 \end{array}\right. \Rightarrow
R_3=[B_{k+1}^t(B_kR_1B_k^t+I)B_{k+1}]^2$$
$$R_3=(B_{k+1}^tB_kR_1B_k^tB_{k+1}+B_{k+1}^tB_{k+1})^2=
(R_1+I)^2=R_1^2=R_1$$ We leave it to the reader to check the
identities $B_{k+1}^tB_k=I$ and $B_{k+1}^tB_{k+1}=I$.
$$\left\{\begin{array}{lll} \overrightarrow{c} &=&
R_3*[(B_{k+1}^t*\overrightarrow{b})\oplus(e_k*x_{k+1})]\\
\overrightarrow{b} &=&B_k*\overrightarrow{a}
\end{array}\right. \Rightarrow
\overrightarrow{c}=
R_3*[(B_{k+1}^t*(B_k*\overrightarrow{a}))\oplus(e_k*x_{k+1})]$$
where $x_{k+1}=[r_{k,k+2}*\varphi(b_{k+1})]\oplus[(\neg
r_{k,k+2})*(b_{k}\wedge b_{k+2})]$ with
$$r_{k,k+2}=e_k^tR_2e_{k+2} =e_k^t(B_kR_1B_k^t+D_k)e_{k+2}
=e_k^tB_kR_1B_k^te_{k+2}+e_k^tD_ke_{k+2} =0$$ thus
$x_{k+1}= b_{k}\wedge b_{k+2}= \emptyset\wedge a_k=\emptyset$.

Then we have $$\overrightarrow{c}=
R_3*[((B_{k+1}^tB_k)*\overrightarrow{a})\oplus(e_k*\emptyset)]
=R_1*(\overrightarrow{a}\oplus\emptyset) =
R_1*\overrightarrow{a}=\overrightarrow{a}$$ To check the identity
$\check{T}_{k-1}\circ\hat{T}_k= id$ we use the same procedure.

\item[2.] $\hat{T}_l\circ\hat{T}_k=
\hat{T}_k\circ\hat{T}_{l-2}$ for $l\geq k+2$.

Let $(R_1,\overrightarrow{a})\in \mathcal{O}_n$,
$(R_2,\overrightarrow{b})=\hat{T}_k(R_1,\overrightarrow{a})$ and
$(R_3,\overrightarrow{c})=\hat{T}_l(R_2,\overrightarrow{b})=
\hat{T}_{l}\circ\hat{T}_k(R_1,\overrightarrow{a})$, and let
$(R'_2,\overrightarrow{b}')=\hat{T}_{l-2}(R_1,\overrightarrow{a})$
and
$(R'_3,\overrightarrow{c}')=\hat{T}_k(R'_2,\overrightarrow{b}')=
\hat{T}_k\circ\hat{T}_{l-2}(R_1,\overrightarrow{a})$. We want to
check that $(R'_3,\overrightarrow{c}')=(R_3,\overrightarrow{c})$.
$$\begin{array}{lll} R_3 &=& B_lR_2B_l^t+D_l\\ &=&
B_l(B_kR_1B_k^t+D_k)B_l^t+D_l\\ &=&
B_lB_kR_1B_k^tB_l^t+B_lD_kB_l^t+D_l\\ &=&
(B_lB_k)R_1(B_lB_k)^t+D_k+D_l
\end{array}$$
and $$\begin{array}{lll} R'_3 &=& B_kR'_2B_k^t+D_k\\ &=&
B_k(B_{l-2}R_1B_{l-2}^t+D_{l-2})B_k^t+D_k\\ &=&
B_kB_{l-2}R_1B_{l-2}^tB_k^t+B_kD_{l-2}B_k^t+D_k\\ &=&
(B_kB_{l-2})R_1(B_kB_{l-2})^t+D_l+D_k
\end{array}$$

Thus $R'_3=R_3$. We leave it to the reader to check
$B_lD_kB_l^t=D_k$, $B_kD_{l-2}B_k^t=D_l$ and $B_lB_k=B_kB_{l-2}$.
$$\left\{\begin{array}{l}
\overrightarrow{c}=B_l*\overrightarrow{b}
=B_l*(B_k*\overrightarrow{a}) =(B_lB_k)*\overrightarrow{a}\\
\overrightarrow{c}'=B_k*\overrightarrow{b}'
=B_k*(B_{l-2}*\overrightarrow{a}) =(B_kB_{l-2})*\overrightarrow{a}
\end{array}\right. \Rightarrow \overrightarrow{c}'=\overrightarrow{c}.$$

\item[3.] $\check{T}_k\circ\hat{T}_l=
\hat{T}_{l-2}\circ\check{T}_k$ and $\check{T}_l\circ\hat{T}_k=
\hat{T}_k\circ\check{T}_{l-2}$ for $l\geq k+2$.

Let $(R_1,\overrightarrow{a})\in \mathcal{O}_n$,
$(R_2,\overrightarrow{b})=\hat{T}_l(R_1,\overrightarrow{a})$ and
$(R_3,\overrightarrow{c})=\check{T}_k(R_2,\overrightarrow{b})=
\check{T}_k\circ\hat{T}_l(R_1,\overrightarrow{a})$, and let
$(R'_2,\overrightarrow{b}')=\check{T}_k(R_1,\overrightarrow{a})$
and
$(R'_3,\overrightarrow{c}')=\hat{T}_{l-2}(R'_2,\overrightarrow{b}')=
\hat{T}_{l-2}\circ\check{T}_k(R_1,\overrightarrow{a})$.

We want to check that
$(R'_3,\overrightarrow{c}')=(R_3,\overrightarrow{c})$.

First the case $l>k+2$ (where we have the identity
$B_k^tB_l=B_{l-2}B_k^t$).

$$\begin{array}{lll} R_3 &=& (B_k^tR_2B_k)^2\\ &=&
[B_k^t(B_lR_1B_l^t+D_l)B_k]^2\\ &=&
(B_k^tB_lR_1B_l^tB_k+B_k^tD_lB_k)^2\\ &=&
(B_{l-2}B_k^tR_1B_kB_{l-2}^t+D_{l-2})^2\\ &=&
B_{l-2}B_k^tR_1B_kB_{l-2}^tB_{l-2}B_k^tR_1B_kB_{l-2}^t
+B_{l-2}B_k^tR_1B_kB_{l-2}^tD_{l-2}\\ & &
+D_{l-2}B_{l-2}B_k^tR_1B_kB_{l-2}^t+D_{l-2}^2\\ &=&
B_{l-2}B_k^tR_1B_kIB_k^tR_1B_kB_{l-2}^t +B_{l-2}B_k^tR_1B_kO
+OB_k^tR_1B_kB_{l-2}^t+D_{l-2}\\ &=&
B_{l-2}B_k^tR_1B_kB_k^tR_1B_kB_{l-2}^t +D_{l-2}\\ &=&
B_{l-2}(B_k^tR_1B_k)^2B_{l-2}^t +D_{l-2}\\ &=&
B_{l-2}R'_2B_{l-2}^t +D_{l-2}\\ &=& R'_3
\end{array}$$

We leave it to the reader to check $B_k^tD_lB_k=D_{l-2}$.

$\overrightarrow{c}=
R_3*[(B_k^t*\overrightarrow{b})\oplus(e_{k-1}*x_k)]$
 where $x_k=[r_{k-1,k+1}^{(2)}*\varphi(b_k)]\oplus[(\neg
r_{k-1,k+1}^{(2)})*(b_{k-1}\wedge b_{k+1})]$ with $
r_{k-1,k+1}^{(2)}=e_{k-1}^tR_2e_{k+1}$ and
$\overrightarrow{b}=B_l*\overrightarrow{a}$.

In this way: $$\begin{array}{lll} \overrightarrow{c} &=&
R_3*[(B_k^tB_l*\overrightarrow{a})\oplus(e_{k-1}*x_k)]\\ &=&
R'_3*[(B_{l-2}B_k^t*\overrightarrow{a})\oplus(e_{k-1}*x_k)]\\ &=&
(B_{l-2}R'_2B_{l-2}^t +D_{l-2})*
[(B_{l-2}B_k^t*\overrightarrow{a})\oplus(e_{k-1}*x_k)]\\ &=&
\{(B_{l-2}R'_2B_{l-2}^t)*
[(B_{l-2}B_k^t*\overrightarrow{a})\oplus(e_{k-1}*x_k)]\}\\ & &\vee
\{D_{l-2}*[(B_{l-2}B_k^t*\overrightarrow{a})\oplus(e_{k-1}*x_k)]\}
\end{array}$$

\begin{lemma}
Let $M=[\mu_{i,j}]$ be a Boolean matrix and let
$\overrightarrow{v}$ and $\overrightarrow{w}$ be two arrays. If,
for each index $i$, we have one of the following situations:
\begin{description}
\item[(1)] $\mu_{i,j}*w_j=\emptyset$ for all $j$
(or $\mu_{i,j}*v_j=\emptyset$ for all $j$);
\item[(2)] $\mu_{i,j}=1$ for, at most, a single index $j$;
\end{description}
then $M*(\overrightarrow{v}\oplus\overrightarrow{w})=
(M*\overrightarrow{v})\oplus(M*\overrightarrow{w})$.
\label{distrib}
\end{lemma}
\TeXButton{Proof}{\proof} Let
$\overrightarrow{x}:=M*(\overrightarrow{v}\oplus\overrightarrow{w})$
and $\overrightarrow{y}:=
(M*\overrightarrow{v})\oplus(M*\overrightarrow{w})$. We have
$$x_i=\bigvee_{j=1}^{n}\mu_{i,j}*(v_j\oplus w_j)
=\bigvee_{j=1}^{n}[(\mu_{i,j}*v_j)\oplus (\mu_{i,j}*w_j)]$$ and
$$y_i=(\bigvee_{j=1}^{n}\mu_{i,j}*v_j)\oplus
(\bigvee_{j=1}^{n}\mu_{i,j}*w_j)$$

\begin{description}

\item[(1)] If $\mu_{i,j}*w_j=\emptyset$ for all $j$ then
$$x_i=\bigvee_{j}[(\mu_{i,j}*v_j)\oplus
\emptyset]=\bigvee_{j}\mu_{i,j}*v_j$$ and
$$y_i=(\bigvee_{j}\mu_{i,j}*v_j)\oplus (\bigvee_{j}\emptyset)=
\bigvee_{j}\mu_{i,j}*v_j$$

\item[(2)] If there exists $k$ such that $\mu_{i,j}=1 \Rightarrow j=k$ then
$$x_i=\bigvee_{j}[(\mu_{i,j}*v_j)\oplus (\mu_{i,j}*w_j)]=
(\mu_{i,k}*v_k)\oplus (\mu_{i,k}*w_k)$$ and
$$y_i=(\bigvee_{j}\mu_{i,j}*v_j)\oplus
(\bigvee_{j=1}^{n}\mu_{i,j}*w_j)= (\mu_{i,k}*v_k)\oplus
(\mu_{i,k}*w_k)$$\end{description} \TeXButton{End
Proof}{\endproof}

With this lemma we have that:
$$B_{l-2}^t*[(B_{l-2}B_k^t*\overrightarrow{a})\oplus(e_{k-1}*x_k)]
=(B_{l-2}^tB_{l-2}B_k^t*\overrightarrow{a})\oplus(B_{l-2}^te_{k-1}*x_k)$$
by the condition (2) of the lemma for $i\not=l-3$, and by the
condition (1) for $i=l-3$; and
$$D_{l-2}*[(B_{l-2}B_k^t*\overrightarrow{a})\oplus(e_{k-1}*x_k)]
=(D_{l-2}B_{l-2}B_k^t*\overrightarrow{a})\oplus(D_{l-2}e_{k-1}*x_k)$$
by the condition (2) of the lemma.

Thus $$\begin{array}{lll} \overrightarrow{c} &=&
\{(B_{l-2}R'_2B_{l-2}^t)*
[(B_{l-2}B_k^t*\overrightarrow{a})\oplus(e_{k-1}*x_k)]\}\\ & &\vee
\{D_{l-2}*[(B_{l-2}B_k^t*\overrightarrow{a})\oplus(e_{k-1}*x_k)]\}\\
&=& \{(B_{l-2}R'_2)*[(B_{l-2}^tB_{l-2}B_k^t*\overrightarrow{a})
\oplus(B_{l-2}^te_{k-1}*x_k)]\}\\ & &\vee
[(D_{l-2}B_{l-2}B_k^t*\overrightarrow{a})\oplus(D_{l-2}e_{k-1}*x_k)]\\
&=& \{(B_{l-2}R'_2)*[(B_k^t*\overrightarrow{a})
\oplus(B_{l-2}^te_{k-1}*x_k)]\}\vee
[(OB_k^t*\overrightarrow{a})\oplus(O*x_k)]\\ &=&
(B_{l-2}R'_2)*[(B_k^t*\overrightarrow{a}) \oplus(e_{k-1}*x_k)]
\end{array}$$
On the other hand $$\begin{array}{lll} \overrightarrow{c}' &=&
B_{l-2}*\overrightarrow{b}'\\ &=&
B_{l-2}*\{R'_2*[(B_k^t*\overrightarrow{a})
\oplus(e_{k-1}*x'_k)]\}\\ &=&
(B_{l-2}R'_2)*[(B_k^t*\overrightarrow{a}) \oplus(e_{k-1}*x'_k)]
\end{array}$$
where $x'_k=[r_{k-1,k+1}^{(1)}*\varphi(a_k)]\oplus[(\neg
r_{k-1,k+1}^{(1)})*(a_{k-1}\wedge a_{k+1})]$ with
$r_{k-1,k+1}^{(1)}=e_{k-1}^tR_1e_{k+1}$.

To check $\overrightarrow{c}'=\overrightarrow{c}$ we only need to
prove that $x'_k=x_k$. $$\begin{array}{lll} r_{k-1,k+1}^{(2)} &=&
e_{k-1}^tR_2e_{k+1}\\ &=& e_{k-1}^t(B_lR_1B_l^t+D_l)e_{k+1}\\ &=&
e_{k-1}^tB_lR_1B_l^te_{k+1}+e_{k-1}^tD_le_{k+1}\\ &=&
e_{k-1}^tR_1e_{k+1}\\ &=& r_{k-1,k+1}^{(1)}
\end{array}$$
$\overrightarrow{b}=B_l*\overrightarrow{a}$ and $l> k+2$ implies
that $b_{k-1}=a_{k-1}$, $b_{k}=a_{k}$ and $b_{k+1}=a_{k+1}$.

Then $$\begin{array}{lll} x_k &=&
[r_{k-1,k+1}^{(2)}*\varphi(b_k)]\oplus[(\neg
r_{k-1,k+1}^{(2)})*(b_{k-1}\wedge b_{k+1})]\\ &=&
[r_{k-1,k+1}^{(1)}*\varphi(a_k)]\oplus[(\neg
r_{k-1,k+1}^{(1)})*(a_{k-1}\wedge a_{k+1})]\\ &=& x'_k
\end{array}$$

Now, let us study the case $l=k+2$.

In this case, in contrast with the case $l>k+2$, we don't have
$B_k^tB_l=B_{l-2}B_k^t$. In fact,
$B_k^tB_{k+2}=(I-D_k)+Q_{k-1,k+1}$ and
$B_kB_{k}^t=(I-D_k)+Q_{k-1,k+1}+Q_{k+1,k-1}$ where
$Q_{\alpha,\beta}=e_{\alpha}e_{\beta}^t$ (i.e. all entries of
$Q_{\alpha,\beta}$ are zero except the entry $(\alpha,\beta)$).
Thus $B_k^tB_{k+2}<B_kB_{k}^t$. $$\begin{array}{lll} R_3 &=&
(B_k^tR_2B_k)^2\\ &=& [B_k^t(B_{k+2}R_1B_{k+2}^t+D_{k+2})B_k]^2\\
&=& (B_k^tB_{k+2}R_1B_{k+2}^tB_k+B_k^tD_{k+2}B_k)^2\\ &=&
(B_k^tB_{k+2}R_1B_{k+2}^tB_k+D_k)^2\\ &=&
B_k^tB_{k+2}R_1B_{k+2}^tB_kB_k^tB_{k+2}R_1B_{k+2}^tB_k
+B_k^tB_{k+2}R_1B_{k+2}^tB_kD_k\\ & &
+D_kB_k^tB_{k+2}R_1B_{k+2}^tB_k+D_k^2\\ &=&
B_k^tB_{k+2}R_1B_{k+2}^tB_kB_k^tB_{k+2}R_1B_{k+2}^tB_k+D_k
\end{array}$$ 
because $$B_k^tB_{k+2}R_1B_{k+2}^tB_kD_k \leq
B_k^tB_{k+2}R_1B_kB_k^tD_k =B_k^tB_{k+2}R_1B_kO=O$$ and
$$D_kB_k^tB_{k+2}R_1B_{k+2}^tB_k
=(B_k^tB_{k+2}R_1B_{k+2}^tB_kD_k)^t=O$$

$$\begin{array}{lll} R'_3 &=& B_kR'_2B_k^t+D_k\\ &=&
B_k(B_k^tR_1B_k)^2B_k^t+D_k\\ &=&
B_kB_k^tR_1B_kB_k^tR_1B_kB_k^t+D_k
\end{array}$$

Since $B_{k+2}^tB_k\leq B_kB_k^t$ and $B_k^tB_{k+2}\leq B_kB_k^t$,
we have $$\begin{array}{lll} R_3 &=&
B_k^tB_{k+2}R_1B_{k+2}^tB_kB_k^tB_{k+2}R_1B_{k+2}^tB_k+D_k\\
&\leq& B_kB_k^tR_1B_kB_k^tB_kB_k^tR_1B_kB_k^t+D_k\\ &=& R'_3
\end{array}$$

On the other hand, since $B_{k+2}^tB_kB_k^tB_{k+2}=B_kB_k^t$, we
have $$\begin{array}{lll} R'_3 &=&
B_kB_k^tR_1B_kB_k^tR_1B_kB_k^t+D_k
\\ &=& B_{k+2}^tB_kB_k^tB_{k+2}R_1B_{k+2}^tB_kB_k^tB_{k+2}
R_1B_{k+2}^tB_kB_k^tB_{k+2}+D_k
\\ &=& B_{k+2}^tB_k(B_k^tB_{k+2}R_1B_{k+2}^tB_kB_k^tB_{k+2}
R_1B_{k+2}^tB_k+D_k)B_k^tB_{k+2}+D_k
\\ &=& B_{k+2}^tB_kR_3B_k^tB_{k+2}+D_k
\\ &\leq& B_kB_k^tR_3B_kB_k^t+D_k \\ &\leq& R_3^3+D_k \\ &=& R_3
\end{array}$$

Thus $R_3=R'_3$.

Here we make use of the following inequalities:
$$B_{k+2}^tB_kD_kB_k^tB_{k+2}\leq B_kB_k^tD_kB_k^tB_{k+2}=O$$ and
$$\begin{array}{lll} R_3 & = & (B_k^tR_2B_k)^2 \\
 & = & [B_k^t(B_{k+2}R_1B_{k+2}^t+D_{k+2})B_k]^2 \\
 & \geq & (B_k^tB_{k+2}R_1B_{k+2}^tB_k)^2 \\
 & \geq & (B_k^tB_{k+2}B_{k+2}^tB_k)^2 \\
 & = & (B_kB_k^t)^2 \\
 & = & B_kB_k^t
 \end{array}$$

$$\begin{array}{lll} \overrightarrow{c} &=&
R_3*[(B_k^t*\overrightarrow{b})\oplus(e_{k-1}*x_k)]\\ &=&
R'_3*[(B_k^tB_{k+2}*\overrightarrow{a})\oplus(e_{k-1}*x_k)]\\ &=&
(B_kR'_2B_k^t +D_k)*
[(B_k^tB_{k+2}*\overrightarrow{a})\oplus(e_{k-1}*x_k)]\\ &=&
\{(B_kR'_2B_k^t)*
[(B_k^tB_{k+2}*\overrightarrow{a})\oplus(e_{k-1}*x_k)]\}\\ & &\vee
\{D_k*[(B_k^tB_{k+2}*\overrightarrow{a})\oplus(e_{k-1}*x_k)]\} \\
&=& (B_kR'_2B_k^t)*
[(B_k^tB_{k+2}*\overrightarrow{a})\oplus(e_{k-1}*x_k)]
\end{array}$$
because, using lemma \ref{distrib}, we have 
$$D_k*[(B_k^tB_{k+2}*\overrightarrow{a})\oplus(e_{k-1}*x_k)]
=(D_kB_k^tB_{k+2}*\overrightarrow{a})\oplus(D_ke_{k-1}*x_k)=
\emptyset$$ since $D_kB_k^tB_{k+2}\leq D_kB_{k}B_k^t=O$ and
$D_ke_{k-1}=O$. $$\begin{array}{lll} \overrightarrow{c}' &=&
B_k*\overrightarrow{b}'\\ &=&
B_k*\{R'_2*[(B_k^t*\overrightarrow{a}) \oplus(e_{k-1}*x'_k)]\}\\
&=& (B_{k}R'_2)*[(B_k^t*\overrightarrow{a}) \oplus(e_{k-1}*x'_k)]
\end{array}$$

\begin{lemma}
Let $\overrightarrow{v}$ be an array with values in a monoid. If
$v_{k+1}\leq v_{k-1}$ then
$B_k^t*\overrightarrow{v}=X_k*\overrightarrow{v}$ where
$X_k=B_k^t(I-D_{k+1})$.
\end{lemma}
\TeXButton{Proof}{\proof} $X_k$ differs from $B_k^t$ only in the
entry $(k-1,k+1)$ which is $0$ in $X_k$ and $1$ in $B_k^t$. Thus
we only need to check the equality
$B_k^t*\overrightarrow{v}=X_k*\overrightarrow{v}$ for the $k-1$
coordinate which is $v_{k-1}\vee v_{k+1}$ in
$B_k^t*\overrightarrow{v}$ and $v_{k-1}$ in
$X_k*\overrightarrow{v}$. Since, by hypothesis, $v_{k+1}\leq
v_{k-1}$ we have the equality. \TeXButton{End Proof}{\endproof}

The $k+1$ coordinate of
$(B_k^tB_{k+2}*\overrightarrow{a})\oplus(e_{k-1}*x_k)$ is
$$\begin{array}{lll}
e_{k+1}^t*[(B_k^tB_{k+2}*\overrightarrow{a})\oplus(e_{k-1}*x_k)]
&=&
(e_{k+1}^tB_k^tB_{k+2}*\overrightarrow{a})\oplus(e_{k+1}^te_{k-1}*x_k)\\
&=& (e_{k+1}^t*\overrightarrow{a})\oplus(0*x_k)\\ &=& a_{k+1}
\end{array}$$
and the $k-1$ coordinate of
$(B_k^tB_{k+2}*\overrightarrow{a})\oplus(e_{k-1}*x_k)$ is
$$\begin{array}{lll}
e_{k-1}^t*[(B_k^tB_{k+2}*\overrightarrow{a})\oplus(e_{k-1}*x_k)]
&=&
(e_{k-1}^tB_k^tB_{k+2}*\overrightarrow{a})\oplus(e_{k-1}^te_{k-1}*x_k)\\
&=& [(e_{k-1}^t+e_{k+1}^t)*\overrightarrow{a}]\oplus(1*x_k)\\ &=&
(a_{k-1}\vee a_{k+1})\oplus x_k
\end{array}$$

Thus, 
we may apply the previous lemma $$\begin{array}{lll}
B_k^t*[(B_k^tB_{k+2}*\overrightarrow{a})\oplus(e_{k-1}*x_k)] &=&
X_k*[(B_k^tB_{k+2}*\overrightarrow{a})\oplus(e_{k-1}*x_k)]\\ &=&
(X_kB_k^tB_{k+2}*\overrightarrow{a})\oplus(X_ke_{k-1}*x_k)\\ &=&
(B_k^t*\overrightarrow{a})\oplus(e_{k-1}*x_k)
\end{array}$$

We leave it to the reader to check that $X_kB_k^tB_{k+2}=B_k^t$.
Therefore $$\begin{array}{lll} \overrightarrow{c} &=&
(B_kR'_2B_k^t)*
[(B_k^tB_{k+2}*\overrightarrow{a})\oplus(e_{k-1}*x_k)]\\ &=&
(B_{k}R'_2)*[(B_k^t*\overrightarrow{a}) \oplus(e_{k-1}*x_k)]
\end{array}$$
and thus $\overrightarrow{c}=\overrightarrow{c}'$ if $x_k=x'_k$
and we can check this using the same proof as was used in the case
$l>k+2$.

Now, let us study the other relation $\check{T}_l\circ\hat{T}_k=
\hat{T}_k\circ\check{T}_{l-2}$ for $l\geq k+2$.

We know that $\check{T}_{l'}\circ\hat{T}_{k'}=
\hat{T}_{k'-2}\circ\check{T}_{l'}$ for $k'\geq l'+2$
 (substituting
 $k$ by $l'$ and $l$ by $k'$ in the relation
 we have already proved).

We will use a mirror symmetry between these two relations. For
that purpose, let us introduce the function {\it mirror symmetry}
defined as follow:
$$\begin{array}{cccc}
  M_n: & \mathcal{O}_n & \longrightarrow & \mathcal{O}_n \\
   & (R,\overrightarrow{v}) & \longmapsto &
   (S_nRS_n,S_n*\overrightarrow{v})
\end{array}$$
where $S_n$ is the $n$-dimensional square matrix defined by
$$S_n:=[s_{i,j}] \mbox{ with } s_{i,j}=1 \mbox{ iff } i+j=n+1.$$
We have that $M_n\circ M_n$ is the identity function on
$\mathcal{O}_n$ since $S_n^2=I$. Also, hearing in mind the
relations $S_{n+2}B_{n,k}S_n =B_{n,n+3-k}$ and
$S_{n+2}D_{n+2,k}S_{n+2} =D_{n+2,n+3-k}$, it is not too
hard\footnote{It is only necessary to check the commutativity of
the two first squares since the other two are obtained from these
by a change of variables. For the second square it is useful to
check first the identity
$x_{n+3-l}(S_{n+2}RS_{n+2},S_{n+2}*\overrightarrow{v})
=x_{l}(R,\overrightarrow{v})$.} to check the commutativity of the
following squares: $$\xymatrix{\mathcal{O}_n \ar @{<->}[d]_{M_n}
\ar[r]^{\hat{T}_k} & \mathcal{O}_{n+2} \ar @{<->}[d]^{M_{n+2}} \\
\mathcal{O}_n \ar[r]_{\hat{T}_{n+3-k}} & \mathcal{O}_{n+2}}
\xymatrix{\mathcal{O}_{n+2} \ar @{<->}[d]_{M_{n+2}}
\ar[r]^{\check{T}_l} & \mathcal{O}_n \ar @{<->}[d]^{M_n}
\\ \mathcal{O}_{n+2} \ar[r]_{\check{T}_{n+3-l}} & \mathcal{O}_n}$$

$$\xymatrix{\mathcal{O}_n \ar @{<->}[d]_{M_n}
\ar[r]^{\check{T}_{l-2}} & \mathcal{O}_{n-2} \ar
@{<->}[d]^{M_{n-2}} \\ \mathcal{O}_n \ar[r]_{\check{T}_{n+3-l}} &
\mathcal{O}_{n-2}} \xymatrix{\mathcal{O}_{n-2} \ar
@{<->}[d]_{M_{n-2}} \ar[r]^{\hat{T}_k} & \mathcal{O}_n \ar
@{<->}[d]^{M_n}
\\ \mathcal{O}_{n-2} \ar[r]_{\hat{T}_{n+1-k}} & \mathcal{O}_n}$$

Thus the outside square in the following diagram:

$$\xymatrix{\mathcal{O}_n \ar[ddd]_{\check{T}_{l-2}} \ar
@{<->}[dr]|{M_n} \ar[rrr]^{\hat{T}_k} & & & \mathcal{O}_{n+2} \ar
@{<->}[dl]|{M_{n+2}} \ar[ddd]^{\check{T}_l} \\ & \mathcal{O}_n
\ar[d]_{\check{T}_{n+3-l}} \ar[r]^{\hat{T}_{n+3-k}} &
\mathcal{O}_{n+2} \ar[d]^{\check{T}_{n+3-l}} & \\ &
\mathcal{O}_{n-2} \ar @{<->}[dl]|{M_{n-2}}
\ar[r]_{\hat{T}_{n+1-k}} & \mathcal{O}_n \ar @{<->}[dr]|{M_n} & \\
\mathcal{O}_{n-2} \ar[rrr]^{\hat{T}_k} & & & \mathcal{O}_{n+2} }$$

commutes (i.e.
$\check{T}_l\circ\hat{T}_k=\hat{T}_k\circ\check{T}_{l-2}$) if and
only if the inside square commutes (i.e.
$\check{T}_{n+3-l}\circ\hat{T}_{n+3-k}=
\hat{T}_{n+1-k}\circ\check{T}_{n+3-l} \Leftrightarrow
\check{T}_{l'}\circ\hat{T}_{k'}=
\hat{T}_{k'-2}\circ\check{T}_{l'}$ by making the change of
variables: $l'=n+3-l$ and $k'=n+3-k$) which we know to be true
since $l\geq k+2 \Rightarrow k'\geq l'+2$.

\item[4.] $\check{T}_k\circ\check{T}_l=
\check{T}_{l-2}\circ\check{T}_k$ for $l\geq k+2$.

Let $(R_1,\overrightarrow{a})\in \mathcal{O}_n$,
$(R_2,\overrightarrow{b})=\check{T}_l(R_1,\overrightarrow{a})$ and
$(R_3,\overrightarrow{c})=\check{T}_k(R_2,\overrightarrow{b})=
\check{T}_k\circ\check{T}_l(R_1,\overrightarrow{a})$, and let
$(R'_2,\overrightarrow{b}')=\check{T}_k(R_1,\overrightarrow{a})$
and
$(R'_3,\overrightarrow{c}')=\check{T}_{l-2}(R'_2,\overrightarrow{b}')=
\check{T}_{l-2}\circ\check{T}_k(R_1,\overrightarrow{a})$.

We want to see that
$(R'_3,\overrightarrow{c}')=(R_3,\overrightarrow{c})$.


$$R_3=(B_k^tR_2B_k)^2=\overline{B_k^tR_2B_k}
=\overline{B_k^t\overline{B_l^tR_1B_l}B_k}
\stackrel{(*)}{=}\overline{B_k^tB_l^tR_1B_lB_k}$$

$$R'_3=\overline{B_{l-2}^tR'_2B_{l-2}}
=\overline{B_{l-2}^t\overline{B_k^tR_1B_k}B_{l-2}}
\stackrel{(*)}{=}\overline{B_{l-2}^tB_k^tR_1B_kB_{l-2}}$$ and
since $B_lB_k=B_kB_{l-2}$ we have $R_3=R'_3$.

Note that $\overline{A}$ means the transitive closure of $A$.

(*) Now let us prove that
$\overline{B_k^t\overline{B_l^tR_2B_l}B_k}
=\overline{B_k^tB_l^tR_2B_lB_k}$ and
$\overline{B_{l-2}^t\overline{B_k^tR_2B_k}B_{l-2}}
=\overline{B_{l-2}^tB_k^tR_2B_kB_{l-2}}$.

\begin{lemma}
Let $A$ be a square matrix, then $$\overline{B_k^tAB_k}=
\overline{B_k^t\overline{(I-D_k)A(I-D_k)}B_k}$$
\end{lemma}
\TeXButton{Proof}{\proof} For any natural number $n$, we have
$$\begin{array}{rcl}
  \overline{B_k^tAB_k} & \geq & (B_k^tAB_k)^n \\
   & = & B_k^tAB_kB_k^tAB_k...B_k^tAB_kB_k^tAB_k \\
   & \geq & B_k^tA(I-D_k)A(I-D_k)...(I-D_k)A(I-D_k)AB_k \\
   & = & B_k^t(I-D_k)A(I-D_k)A(I-D_k)...(I-D_k)A(I-D_k)A(I-D_k)B_k \\
   & = & B_k^t[(I-D_k)A(I-D_k)]^nB_k
\end{array}$$
therefore $\overline{B_k^tAB_k}\geq
B_k^t\overline{(I-D_k)A(I-D_k)}B_k$ and thus
$$\overline{B_k^tAB_k}\geq
\overline{B_k^t\overline{(I-D_k)A(I-D_k)}B_k}.$$

On the other hand, $$ B_k^tAB_k=B_k^t(I-D_k)A(I-D_k)B_k \leq
B_k^t\overline{(I-D_k)A(I-D_k)}B_k$$ and then
$\overline{B_k^tAB_k}\leq
\overline{B_k^t\overline{(I-D_k)A(I-D_k)}B_k}$. \TeXButton{End
Proof}{\endproof}

\begin{corollary}
If $(I-D_k)\overline{A}(I-D_k)\leq \overline{(I-D_k)A(I-D_k)}$
then $\overline{B_k^t\overline{A}B_k}=\overline{B_k^tAB_k}$.
\end{corollary}
\TeXButton{Proof}{\proof}
$$\begin{array}{rcl}
  \overline{B_k^t\overline{A}B_k} & = &
  \overline{B_k^t(I-D_k)\overline{A}(I-D_k)B_k} \\
   & \leq & \overline{B_k^t\overline{(I-D_k)A(I-D_k)}B_k} \\
   & = & \overline{B_k^tAB_k}
\end{array}$$

On the other hand, $$\overline{A}\geq A \Rightarrow
B_k^t\overline{A}B_k\geq B_k^tAB_k \Rightarrow
\overline{B_k^t\overline{A}B_k}\geq \overline{B_k^tAB_k}$$
\TeXButton{End Proof}{\endproof}

\begin{quote}
{\bf Claim:} $(I-D_k)\overline{B_l^tR_1B_l}(I-D_k) \leq
\overline{(I-D_k)B_l^tR_1B_l(I-D_k)}$.
\end{quote}

\TeXButton{Proof}{\proof} $$\begin{array}{rcl}
  \overline{(I-D_k)B_l^tR_1B_l(I-D_k)} & \geq &
  [(I-D_k)B_l^tR_1B_l(I-D_k)]^2 \\
   & = & (I-D_k)B_l^tR_1B_l(I-D_k)B_l^tR_1B_l(I-D_k)
\end{array}$$

$$\begin{array}{rcl}
  (I-D_k)\overline{B_l^tR_1B_l}(I-D_k) & = &
  (I-D_k)B_l^tR_1B_lB_l^tR_1B_l(I-D_k) \\ & = &
  (I-D_k)B_l^tR_1B_l[D_k+(I-D_k)]B_l^tR_1B_l(I-D_k) \\ & = &
  (I-D_k)B_l^tR_1B_lD_kB_l^tR_1B_l(I-D_k)\\ & & +
  (I-D_k)B_l^tR_1B_l(I-D_k)B_l^tR_1B_l(I-D_k)
\end{array}$$

$$\begin{array}{rcl}
  (I-D_k)B_l^tR_1B_lD_kB_l^tR_1B_l(I-D_k) & = &
  (I-D_k)B_l^tR_1D_kR_1B_l(I-D_k) \\
   & \leq & (I-D_k)B_l^tR_1B_l(I-D_k)
\end{array}$$

$$\begin{array}{rcl}
  (I-D_k)B_l^tR_1B_l(I-D_k)B_l^tR_1B_l(I-D_k) & \geq &
  (I-D_k)B_l^tR_1B_l(I-D_k)B_l^tB_l(I-D_k) \\
   & = & (I-D_k)B_l^tR_1B_l(I-D_k)
\end{array}$$

thus $$(I-D_k)B_l^tR_1B_lD_kB_l^tR_1B_l(I-D_k)\leq
(I-D_k)B_l^tR_1B_l(I-D_k)B_l^tR_1B_l(I-D_k).$$

And therefore, $$\begin{array}{ccl}
  (I-D_k)\overline{B_l^tR_1B_l}(I-D_k) & = &
  (I-D_k)B_l^tR_1B_l(I-D_k)B_l^tR_1B_l(I-D_k) \\
   & \leq & \overline{(I-D_k)B_l^tR_1B_l(I-D_k)}
\end{array}$$
\TeXButton{End Proof}{\endproof}

Using the same argument we can also prove the following claim.

\begin{quote}
{\bf Claim:} $(I-D_{l-2})\overline{B_k^tR_1B_k}(I-D_{l-2}) \leq
\overline{(I-D_{l-2})B_k^tR_1B_k(I-D_{l-2})}$.
\end{quote}

Therefore we have $\overline{B_k^t\overline{B_l^tR_2B_l}B_k}
=\overline{B_k^tB_l^tR_2B_lB_k}$ and
$\overline{B_{l-2}^t\overline{B_k^tR_2B_k}B_{l-2}}
=\overline{B_{l-2}^tB_k^tR_2B_kB_{l-2}}$.

Now let us see that $\overrightarrow{c}'=\overrightarrow{c}$. We
will check $c'_i=c_i$ for each index $i$.

\begin{lemma}
Let $R=[r_{i,j}]$ be a matrix with the properties E1, E2 and E3
(i.e. an equivalence relation) and $\overrightarrow{v}$ an array
fixed by the action of $R$.
\begin{description}
\item[1.] If $r_{i,j}=1$ then
$e_i^t*\overrightarrow{v}=e_j^t*\overrightarrow{v}$ (i.e.
$v_i=v_j$).

Now let $\dot{R}=[\dot{r}_{i,j}]=(B_{\alpha}^tRB_{\alpha})^2$.

\item[2.] If $\dot{r}_{i,\alpha-1}=0$ then
$$e_i^t*\{\dot{R}*[(B_{\alpha}^t*\overrightarrow{v})
\oplus(e_{\alpha-1}*x)]\}=e_i^tB_{\alpha}^t*\overrightarrow{v}
=\left\{\begin{array}{ccc}
  v_i & \mbox{if} & i<\alpha-1 \\
  v_{i-2} & \mbox{if} & i>\alpha-1
\end{array}\right. $$

\item[3.] $e_{\alpha-1}^t*\{\dot{R}*[(B_{\alpha}^t*\overrightarrow{v})
\oplus(e_{\alpha-1}*x)]\}=(v_{\alpha-1}\vee v_{\alpha+1})\oplus
x$.
\end{description}
\label{11v}
\end{lemma}
\TeXButton{Proof}{\proof}
\begin{description}
\item[1.] $$\begin{array}{ccll}
  r_{i,j}=1 \Rightarrow e_i^t*\overrightarrow{v} & = &
  e_i^tR*\overrightarrow{v} &
  (\mbox{since } \overrightarrow{v}=R*\overrightarrow{v}) \\
    & \geq & e_i^tRe_je_j^t*\overrightarrow{v} &
  (\mbox{since } e_je_j^t=D_j\leq I) \\
    & = & e_j^t*\overrightarrow{v} &
  (\mbox{since } e_iRe_j^t=r_{i,j}=1)
\end{array}$$

Since $R$ is symmetric, we have also
$e_j^tR*\overrightarrow{v}\geq e_i^tR*\overrightarrow{v}$.

\item[2.] $$\begin{array}{ccl}
  \dot{r}_{i,\alpha-1}=0 \Rightarrow e_i^t\dot{R} & = &
  e_i^t\dot{R}[(I-D_{\alpha-1})+D_{\alpha-1}] \\
    & = & e_i^t\dot{R}(I-D_{\alpha-1})
    +e_i^t\dot{R}D_{\alpha-1} \\
    & = & e_i^t\dot{R}(I-D_{\alpha-1})
    +e_i^t\dot{R}e_{\alpha-1}e_{\alpha-1}^t \\
    & = & e_i^t\dot{R}(I-D_{\alpha-1})
\end{array}$$

Thus $$\begin{array}{ccl}
  e_i^t\dot{R}*[(B_{\alpha}^t*\overrightarrow{v})
\oplus(e_{\alpha-1}*x)] & = &
  e_i^t\dot{R}(I-D_{\alpha-1})*[(B_{\alpha}^t*\overrightarrow{v})
\oplus(e_{\alpha-1}*x)] \\
    & = &
    e_i^t\dot{R}*\{[(I-D_{\alpha-1})B_{\alpha}^t*\overrightarrow{v}]\\
    & & \oplus[(I-D_{\alpha-1})e_{\alpha-1}*x]\} \\
    & = & e_i^t\dot{R}(I-D_{\alpha-1})B_{\alpha}^t*\overrightarrow{v} \\
    & = & e_i^tB_{\alpha}^tRB_{\alpha}B_{\alpha}^tR
    B_{\alpha}(I-D_{\alpha-1})B_{\alpha}^t*\overrightarrow{v} \\
    & \leq & e_i^tB_{\alpha}^tRB_{\alpha}B_{\alpha}^tR
    *\overrightarrow{v} \\ & = & e_i^tB_{\alpha}^tRB_{\alpha}B_{\alpha}^t
    *\overrightarrow{v}
\end{array}$$

Since $$e_i^tB_{\alpha}^tRB_{\alpha}e_{\alpha-1}\leq
e_i^t\dot{R}e_{\alpha-1}=0 \Rightarrow
e_i^tB_{\alpha}^tRB_{\alpha}=
e_i^tB_{\alpha}^tRB_{\alpha}(I-D_{\alpha-1})$$ we have that
$$\begin{array}{ccl}
  e_i^t\dot{R}*[(B_{\alpha}^t*\overrightarrow{v})
\oplus(e_{\alpha-1}*x)] & \leq &
  e_i^tB_{\alpha}^tRB_{\alpha}B_{\alpha}^t
    *\overrightarrow{v} \\
    & = & e_i^tB_{\alpha}^tRB_{\alpha}(I-D_{\alpha-1})B_{\alpha}^t
    *\overrightarrow{v} \\
    & \leq & e_i^tB_{\alpha}^tR*\overrightarrow{v} \\
    & = & e_i^tB_{\alpha}^t*\overrightarrow{v}
\end{array}$$

On the other hand, we have
$$e_i^t\dot{R}*[(B_{\alpha}^t*\overrightarrow{v})
\oplus(e_{\alpha-1}*x)] \geq
e_i^t\dot{R}*(B_{\alpha}^t*\overrightarrow{v}) \geq
e_i^tB_{\alpha}^t*\overrightarrow{v}$$

\item[3.] $$\begin{array}{ccl}
  e_{\alpha-1}^t\dot{R}*[(B_{\alpha}^t*\overrightarrow{v})
\oplus(e_{\alpha-1}*x)] & \geq &
  e_{\alpha-1}^t*[(B_{\alpha}^t*\overrightarrow{v})
\oplus(e_{\alpha-1}*x)] \\
    & = & (e_{\alpha-1}^tB_{\alpha}^t*\overrightarrow{v})
\oplus(e_{\alpha-1}^te_{\alpha-1}*x) \\
    & = & (v_{\alpha-1}\vee v_{\alpha+1})\oplus x
\end{array}$$

On the other hand $$\begin{array}{ccl}
  e_{\alpha-1}^t\dot{R}*[(B_{\alpha}^t*\overrightarrow{v})
\oplus(e_{\alpha-1}*x)] & \leq &
  (e_{\alpha-1}^t\dot{R}B_{\alpha}^t*\overrightarrow{v})
\oplus(e_{\alpha-1}^t\dot{R}e_{\alpha-1}*x) \\
    & = & (e_{\alpha-1}^t\dot{R}B_{\alpha}^t*\overrightarrow{v})
\oplus x \\
    & \stackrel{(*)}{=} & (v_{\alpha-1}\vee v_{\alpha+1})\oplus x
\end{array}$$
\end{description}

\begin{description}
\item[(*)] $e_{\alpha-1}^t\dot{R}B_{\alpha}^t*\overrightarrow{v}
\geq e_{\alpha-1}^tB_{\alpha}^t*\overrightarrow{v}
=v_{\alpha-1}\vee v_{\alpha+1}$, and on the other hand
\end{description}

$$\begin{array}{ccl}
  e_{\alpha-1}^t\dot{R}B_{\alpha}^t*\overrightarrow{v} & = &
  e_{\alpha-1}^tB_{\alpha}^tRB_{\alpha}B_{\alpha}^tR
  B_{\alpha}B_{\alpha}^t*\overrightarrow{v} \\
    & \leq & e_{\alpha-1}^tB_{\alpha}^tRB_{\alpha}B_{\alpha}^tR
  (I+B_{\alpha}D_{\alpha-1}B_{\alpha}^t)*\overrightarrow{v}  \\
    & = & (e_{\alpha-1}^tB_{\alpha}^tRB_{\alpha}B_{\alpha}^tR
  *\overrightarrow{v})
  \vee (e_{\alpha-1}^tB_{\alpha}^tRB_{\alpha}B_{\alpha}^tR
  B_{\alpha}e_{\alpha-1}e_{\alpha-1}^tB_{\alpha}^t*\overrightarrow{v})\\
    & = & (e_{\alpha-1}^tB_{\alpha}^tRB_{\alpha}B_{\alpha}^t
  *\overrightarrow{v}) \vee (e_{\alpha-1}^tB_{\alpha}^t*\overrightarrow{v})\\
    & \leq & [e_{\alpha-1}^tB_{\alpha}^tR
    (I+B_{\alpha}D_{\alpha-1}B_{\alpha}^t)*\overrightarrow{v}]
  \vee (e_{\alpha-1}^tB_{\alpha}^t*\overrightarrow{v})\\
    & = & (e_{\alpha-1}^tB_{\alpha}^tR*\overrightarrow{v}) \vee
    (e_{\alpha-1}^tB_{\alpha}^tRB_{\alpha}e_{\alpha-1}e_{\alpha-1}^t
    B_{\alpha}^t*\overrightarrow{v})
  \vee (e_{\alpha-1}^tB_{\alpha}^t*\overrightarrow{v})\\
  & = & e_{\alpha-1}^tB_{\alpha}^t*\overrightarrow{v}
\end{array}$$

\TeXButton{End Proof}{\endproof}

\begin{lemma}
Let $R=[r_{i,j}]$ be a matrix with the properties E1, E2 and E3
(i.e. an equivalence relation) and let
$\dot{R}=[\dot{r}_{i,j}]=(B_{\alpha}^tRB_{\alpha})^2$. Then:

\begin{description}
\item[(a)] For any $j\not=\alpha-1$,
$\dot{r}_{\alpha-1,j}=r_{\alpha-1,\check{\alpha}(j)}+
r_{\alpha+1,\check{\alpha}(j)}$ where
$\check{\alpha}(j)=\left\{\begin{array}{ccc}
  j & \mbox{if} & j<\alpha-1 \\
  j+2 & \mbox{if} & j>\alpha-1
\end{array}\right.$;

\item[(b)] If $i\not=\alpha-1$ and $j\not=\alpha-1$ then
$\dot{r}_{i,j}\geq r_{\check{\alpha}(i),\check{\alpha}(j)}$;

\item[(c)] If $i\not=\alpha-1$, $j\not=\alpha-1$ and
$\dot{r}_{i,\alpha-1}=0$ (or $\dot{r}_{j,\alpha-1}=0$) then
$\dot{r}_{i,j}= r_{\check{\alpha}(i),\check{\alpha}(j)}$.
\end{description}
\label{12R}
\end{lemma}
\TeXButton{Proof}{\proof}
\begin{description}
\item[(a)] $$\begin{array}{ccl}
  \dot{r}_{\alpha-1,j} & = & e_{\alpha-1}^t\dot{R}e_j \\
    & = & e_{\alpha-1}^t(B_{\alpha}^tRB_{\alpha})^2e_j \\
    & \geq & e_{\alpha-1}^tB_{\alpha}^tRB_{\alpha}e_j \\
    & = & (e_{\alpha-1}^t+e_{\alpha+1}^t)RB_{\alpha}e_j \\
    & = & r_{\alpha-1,\check{\alpha}(j)}+
r_{\alpha+1,\check{\alpha}(j)}
\end{array}$$

On the other hand $$\begin{array}{ccl}
  \dot{r}_{\alpha-1,j} & = & e_{\alpha-1}^tB_{\alpha}^tRB_{\alpha}
  B_{\alpha}^tRB_{\alpha}e_j \\
    & \leq & e_{\alpha-1}^tB_{\alpha}^tR(I+B_{\alpha}D_{\alpha-1}
  B_{\alpha}^t)RB_{\alpha}e_j \\
    & = & e_{\alpha-1}^tB_{\alpha}^tRB_{\alpha}e_j+
    e_{\alpha-1}^tB_{\alpha}^tRB_{\alpha}e_{\alpha-1}e_{\alpha-1}^t
  B_{\alpha}^tRB_{\alpha}e_j \\
    & = & e_{\alpha-1}^tB_{\alpha}^tRB_{\alpha}e_j \\
    & = & r_{\alpha-1,\check{\alpha}(j)}+
r_{\alpha+1,\check{\alpha}(j)}
\end{array}$$

\item[(b)] $$\begin{array}{ccl}
  \dot{r}_{i,j} & = & e_i^t(B_{\alpha}^tRB_{\alpha})^2e_j \\
    & \geq & e_i^tB_{\alpha}^tRB_{\alpha}e_j \\
    & = & e_{\check{\alpha}(i)}^tRe_{\check{\alpha}(j)} \\
    & = & r_{\check{\alpha}(i),\check{\alpha}(j)}
\end{array}$$

\item[(c)] $$\begin{array}{ccl}
  \dot{r}_{i,j} & = & e_i^t(B_{\alpha}^tRB_{\alpha})^2e_j \\
    & = & e_i^tB_{\alpha}^tRB_{\alpha}B_{\alpha}^tRB_{\alpha}e_j \\
    & \leq & e_i^tB_{\alpha}^tR(I+B_{\alpha}D_{\alpha-1}
  B_{\alpha}^t)RB_{\alpha}e_j \\
    & = & e_i^tB_{\alpha}^tRB_{\alpha}e_j+
    e_i^tB_{\alpha}^tRB_{\alpha}e_{\alpha-1}e_{\alpha-1}^t
  B_{\alpha}^tRB_{\alpha}e_j \\
    & \stackrel{(*)}{=} & e_i^tB_{\alpha}^tRB_{\alpha}e_j \\
    & = & r_{\check{\alpha}(i),\check{\alpha}(j)}
\end{array}$$ (*) $e_i^tB_{\alpha}^tRB_{\alpha}e_{\alpha-1}
\leq e_i^t(B_{\alpha}^tRB_{\alpha})^2e_{\alpha-1}
=\dot{r}_{i,\alpha-1}=0$ (or
$e_{\alpha-1}^tB_{\alpha}^tRB_{\alpha}e_j \leq
\dot{r}_{\alpha-1,j}=0$).

And thus, by (b), we have $\dot{r}_{i,j}=
r_{\check{\alpha}(i),\check{\alpha}(j)}$.
\end{description}
\TeXButton{End Proof}{\endproof}

Now, we are going to prove that $c_{k-1}=c'_{k-1}$ for the case
$\ddot{r}_{k-1,l-3}=0$.

{\it Convention:} $R_1=[r_{i,j}]$,
$R_2=[\dot{r}_{i,j}]=(B_l^tR_1B_l)^2$,
$R_3=[\ddot{r}_{i,j}]=(B_k^tR_2B_k)^2$,
$R'_2=[\dot{r}'_{i,j}]=(B_k^tR_1B_k)^2$ and
$R'_3=[\ddot{r}'_{i,j}]=(B_{l-2}^tR'_2B_{l-2})^2$.

By lemma \ref{11v}, we have:

$$\begin{array}{ccl}
  c_{k-1} & = & e_{k-1}^tR_3*
  [(B_k^t*\overrightarrow{b})\oplus(e_{k-1}*x_k)] \\
          & = & (b_{k-1}\vee b_{k+1})\oplus x_k
\end{array}$$
where $x_k=[\dot{r}_{k-1,k+1}*\varphi(b_k)]\oplus[(\neg
\dot{r}_{k-1,k+1})*(b_{k-1}\wedge b_{k+1})]$ with $
\dot{r}_{k-1,k+1}=e_{k-1}^tR_2e_{k+1}$.

By lemma \ref{12R}, $\ddot{r}_{k-1,l-3}=0 \Rightarrow
\dot{r}_{k-1,l-1}=0 \mbox{ and } \dot{r}_{k+1,l-1}=0$. And then,
by lemma \ref{11v}, $$b_{k-1}=e_{k-1}^tR_2*
  [(B_l^t*\overrightarrow{a})\oplus(e_{l-1}*x_l)]=a_{k-1}$$ and
$$b_{k+1}=e_{k+1}^tR_2*
  [(B_l^t*\overrightarrow{a})\oplus(e_{l-1}*x_l)]=a_{k+1}$$

Thus $$c_{k-1}=(a_{k-1}\vee a_{k+1})\oplus x_k$$

Since $R'_3=R_3$, $\ddot{r}'_{k-1,l-3}=\ddot{r}_{k-1,l-3}=0$.
Then, by lemma \ref{11v}, $$
\begin{array}{ccl}
  c'_{k-1} & = & e_{k-1}^tR'_3*
  [(B_{l-2}^t*\overrightarrow{b})\oplus(e_{l-3}*x'_{l-2})] \\
    & = & b'_{k-1} \\
    & = & e_{k-1}^tR'_2*
  [(B_k^t*\overrightarrow{a})\oplus(e_{k-1}*x'_k)] \\
    & = & (a_{k-1}\vee a_{k+1})\oplus x'_k
\end{array}$$
where $$x'_k=[r_{k-1,k+1}*\varphi(a_k)]\oplus[(\neg
r_{k-1,k+1})*(a_{k-1}\wedge a_{k+1})]$$ with $
r_{k-1,k+1}=e_{k-1}^tR_1e_{k+1}$.

Thus $c_{k-1}=c'_{k-1}$ if $x_k=x'_k$. $$
\begin{array}{ccl}
  x_k & = & [\dot{r}_{k-1,k+1}*\varphi(b_k)]\oplus[(\neg
\dot{r}_{k-1,k+1})*(b_{k-1}\wedge b_{k+1})] \\
    & = & [\dot{r}_{k-1,k+1}*\varphi(b_k)]\oplus[(\neg
\dot{r}_{k-1,k+1})*(a_{k-1}\wedge a_{k+1})] \\
    & \stackrel{(*)}{=} & [\dot{r}_{k-1,k+1}*\varphi(a_k)]\oplus[(\neg
\dot{r}_{k-1,k+1})*(a_{k-1}\wedge a_{k+1})] \\
    & \stackrel{(\dag)}{=} & [r_{k-1,k+1}*\varphi(a_k)]\oplus[(\neg
r_{k-1,k+1})*(a_{k-1}\wedge a_{k+1})] \\
    & = & x'_k
\end{array}$$

\begin{description}
\item[(*)] If $\dot{r}_{k-1,k+1}=0$ then
$\dot{r}_{k-1,k+1}*\varphi(b_k)=
\emptyset=\dot{r}_{k-1,k+1}*\varphi(a_k)$.

If $\dot{r}_{k-1,k+1}=1$ then $\dot{r}_{k,l-1}=0$ by the
properties T1 and T2 of $R_2$. Thus $b_{k}=e_{k}^tR_2*
  [(B_l^t*\overrightarrow{a})\oplus(e_{l-1}*x_l)]=a_{k}$.

\item[($\dag$)] $\ddot{r}_{k-1,l-3}=0 \Rightarrow
\dot{r}_{k-1,l-1}=0 \mbox{ and } \dot{r}_{k+1,l-1}=0$ by lemma 11.

$\dot{r}_{k-1,l-1}=0 \Rightarrow \dot{r}_{k-1,k+1}=r_{k-1,k+1}$ by
lemma 11.
\end{description}

We now prove $c_{k-1}=c'_{k-1}$ for the case
$\ddot{r}_{k-1,l-3}=1$.

By lemma 12, $$\begin{array}{ccc}
  \ddot{r}_{k-1,l-3}=1 & \Rightarrow & \dot{r}_{k-1,l-1}=1 \mbox{ or }
   \dot{r}_{k+1,l-1}=1 \\
       & \Rightarrow & r_{k-1,l-1}=1 \mbox{ or }
   r_{k-1,l+1}=1 \mbox{ or } r_{k+1,l-1}=1 \mbox{ or }
   r_{k+1,l+1}=1
\end{array}$$

Now we have, in theory, 15 cases to study: $$(r_{k-1,l-1},
r_{k-1,l+1},r_{k+1,l-1},r_{k+1,l+1})\in \{0,1\}^4\setminus
\{(0,0,0,0)\}$$ but we can exclude the cases $(r_{k-1,l-1},
r_{k-1,l+1},r_{k+1,l-1},r_{k+1,l+1})=(1,1,1,0)$, $(1,1,0,1)$,
$(1,0,1,1)$ and $(0,1,1,1)$, because $R_1=[r_{i,j}]$ satisfies the
properties E1, E2 and E3, and the case $(r_{k-1,l-1},
r_{k-1,l+1},r_{k+1,l-1},r_{k+1,l+1})=(1,0,0,1)$, because the
inequalities $k-1<k+1\leq l-1<l+1$ and property T2 of $R_1$ imply
that if $r_{k-1,l-1}=r_{k+1,l+1}=1$ then $r_{k+1,l-1}=1$ and
$r_{k-1,l+1}=1$.

Thus, we have the following ten cases to study:

\begin{tabular}{|c|c|c|c|c|}  \hline \hline
  $r_{k-1,l-1}$ & $r_{k-1,l+1}$ & $r_{k+1,l-1}$ & $r_{k+1,l+1}$ &
   $\mbox{one geometric realization}$ \\ \hline \hline
  1 & 1 & 1 & 1 & $\psdiag{5}{10}{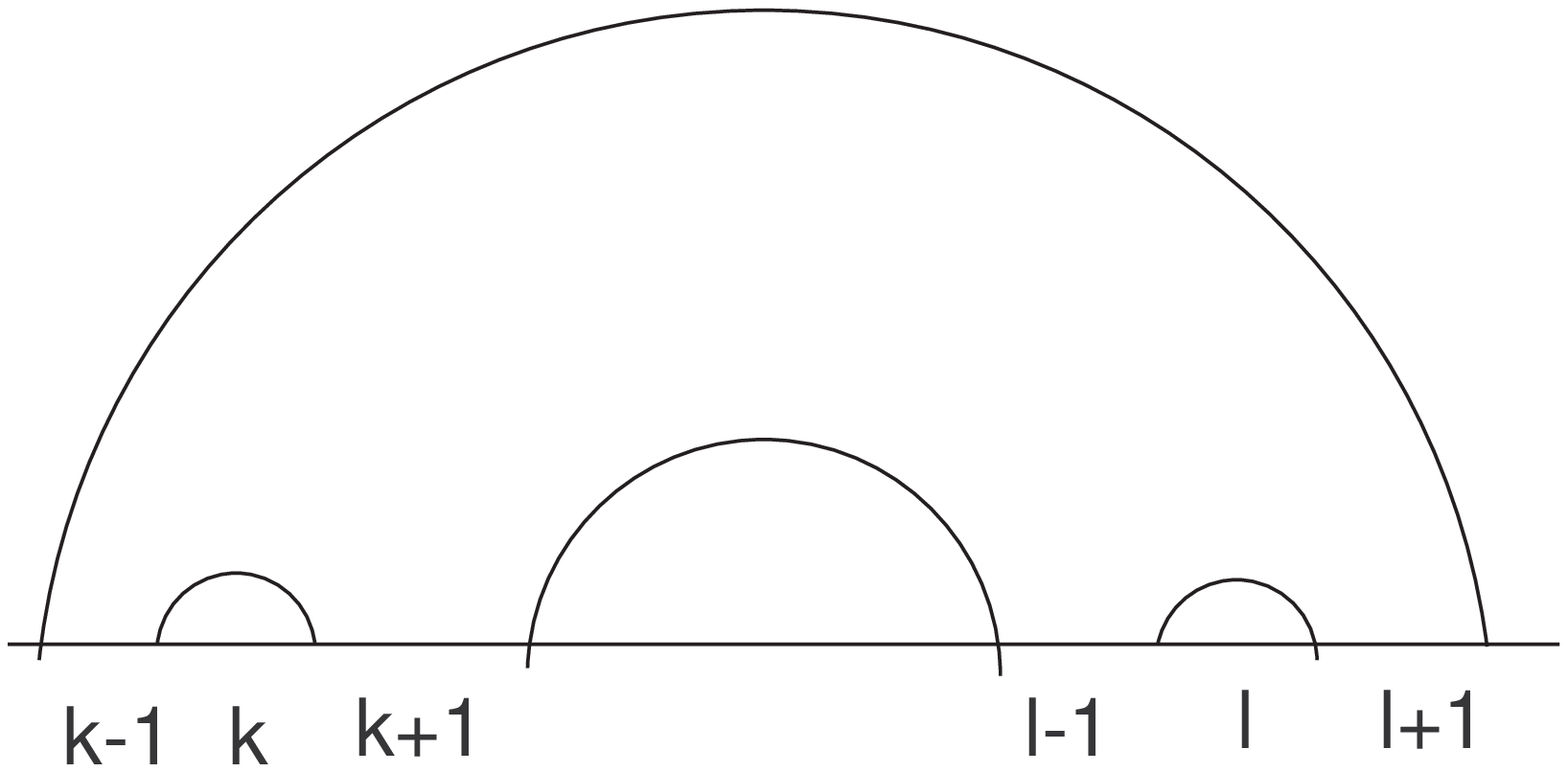}$  \\ \hline
  1 & 0 & 1 & 0 & $\psdiag{5}{10}{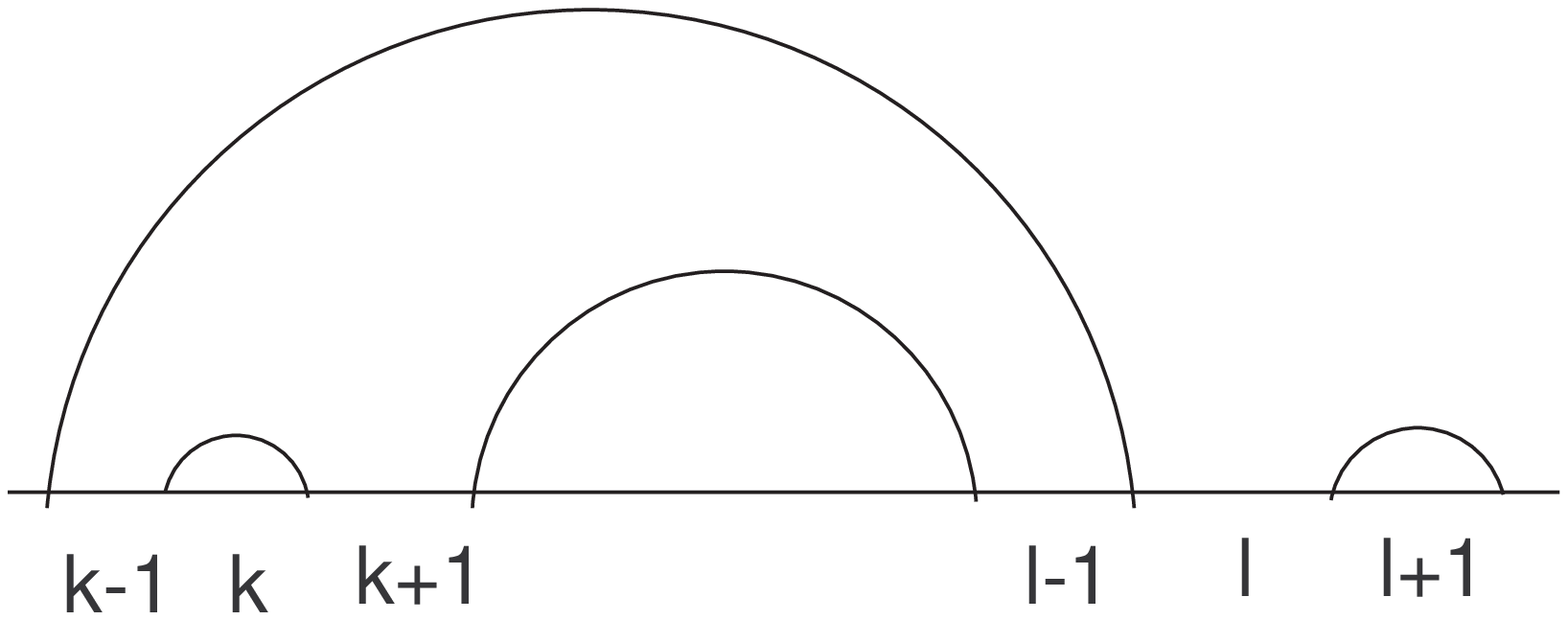}$  \\ \hline
  0 & 1 & 0 & 1 & $\psdiag{5}{10}{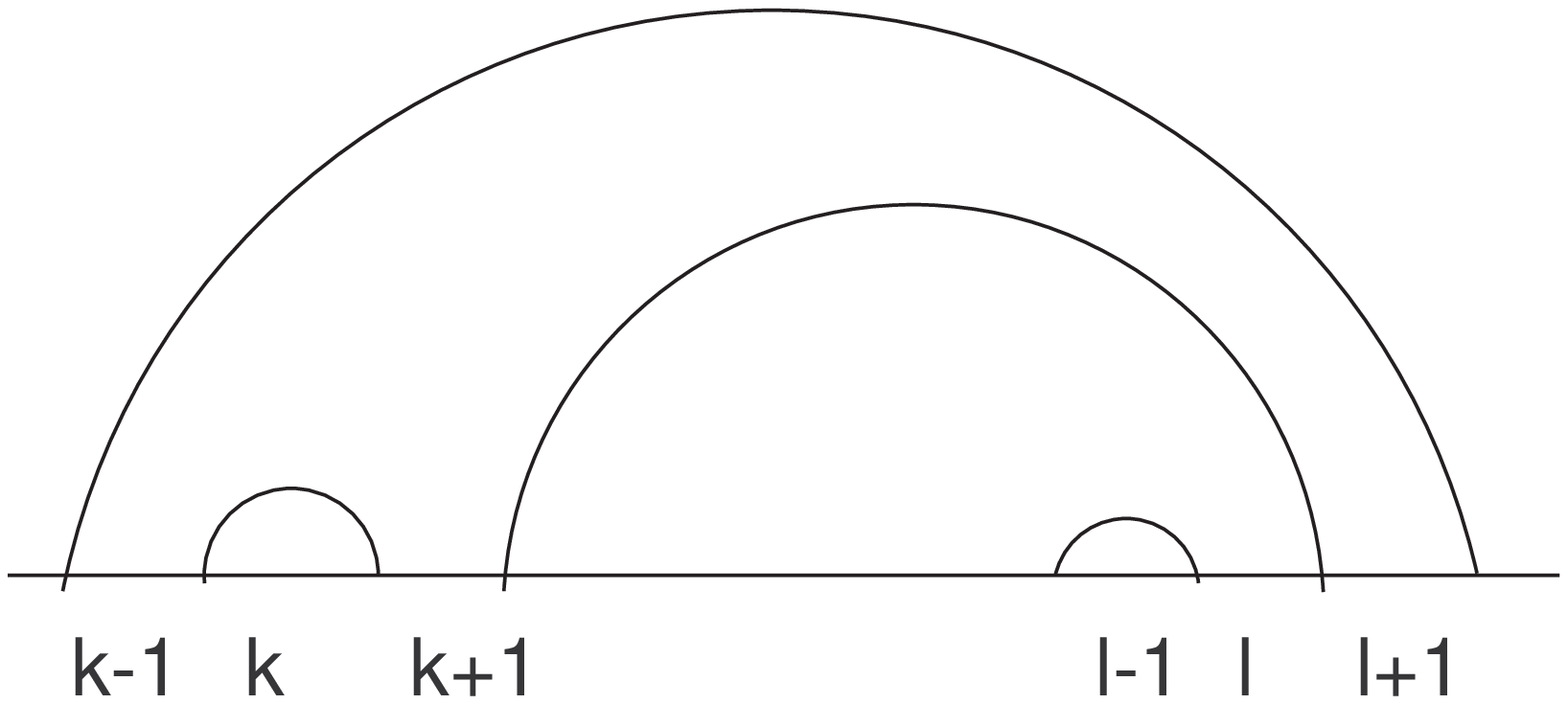}$  \\ \hline
  1 & 1 & 0 & 0 & $\psdiag{5}{10}{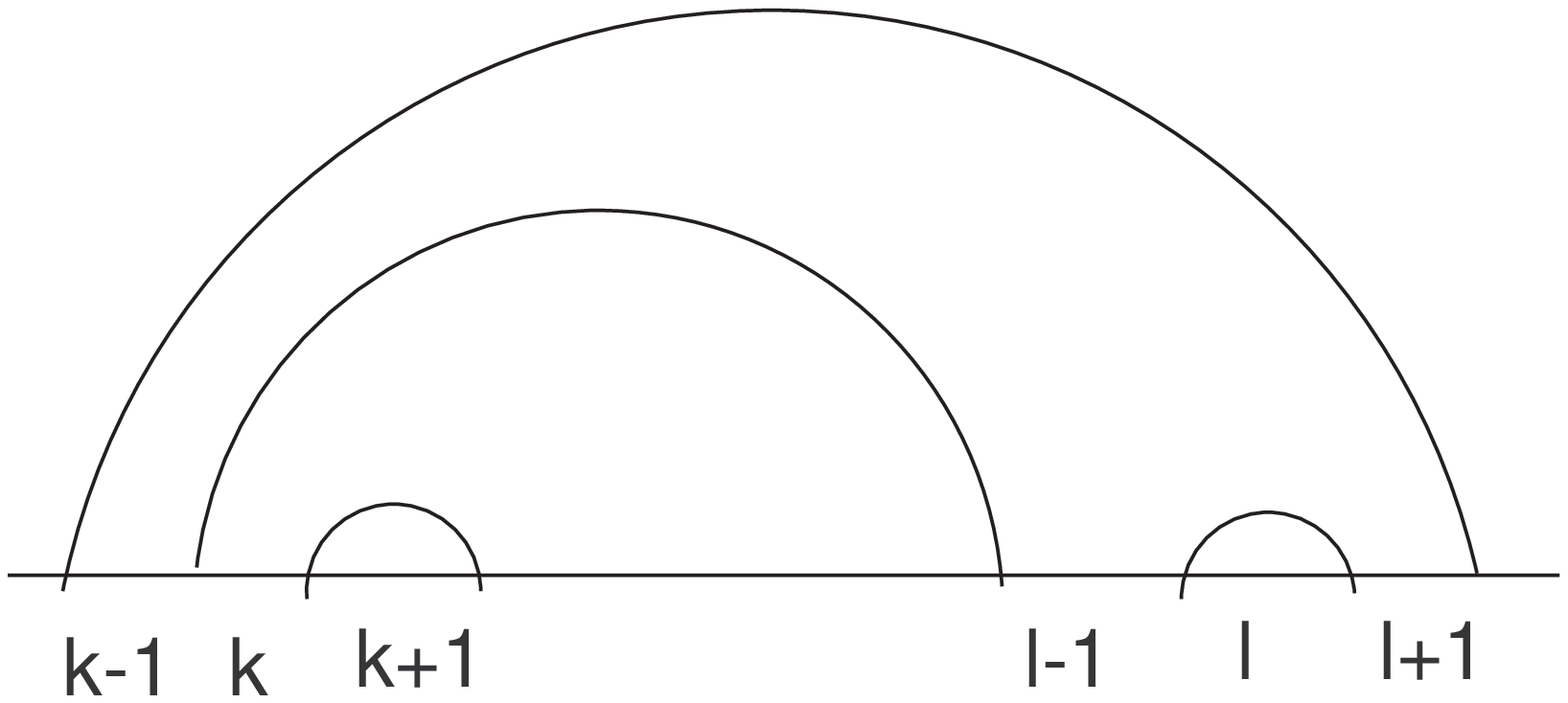}$  \\ \hline
  0 & 0 & 1 & 1 & $\psdiag{5}{10}{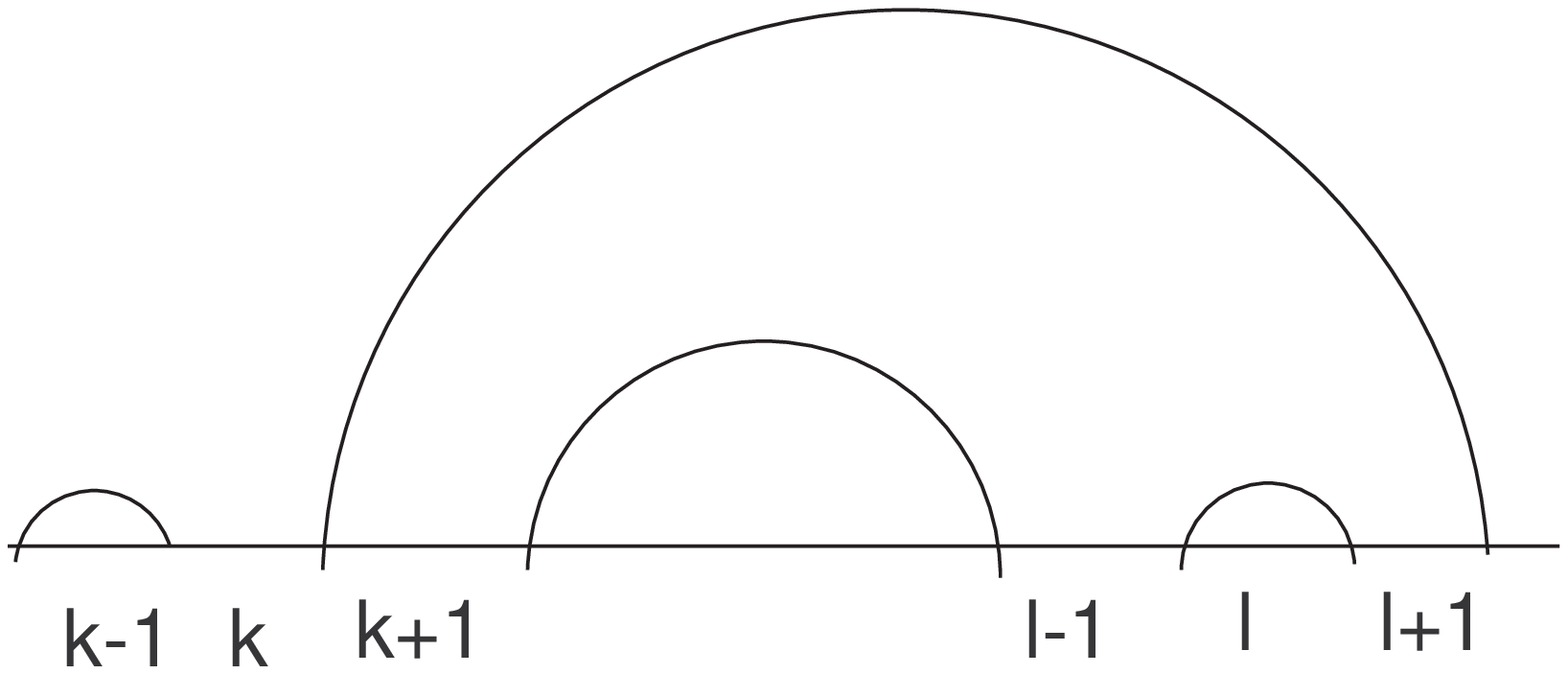}$  \\ \hline
  0 & 1 & 1 & 0 & $\psdiag{5}{10}{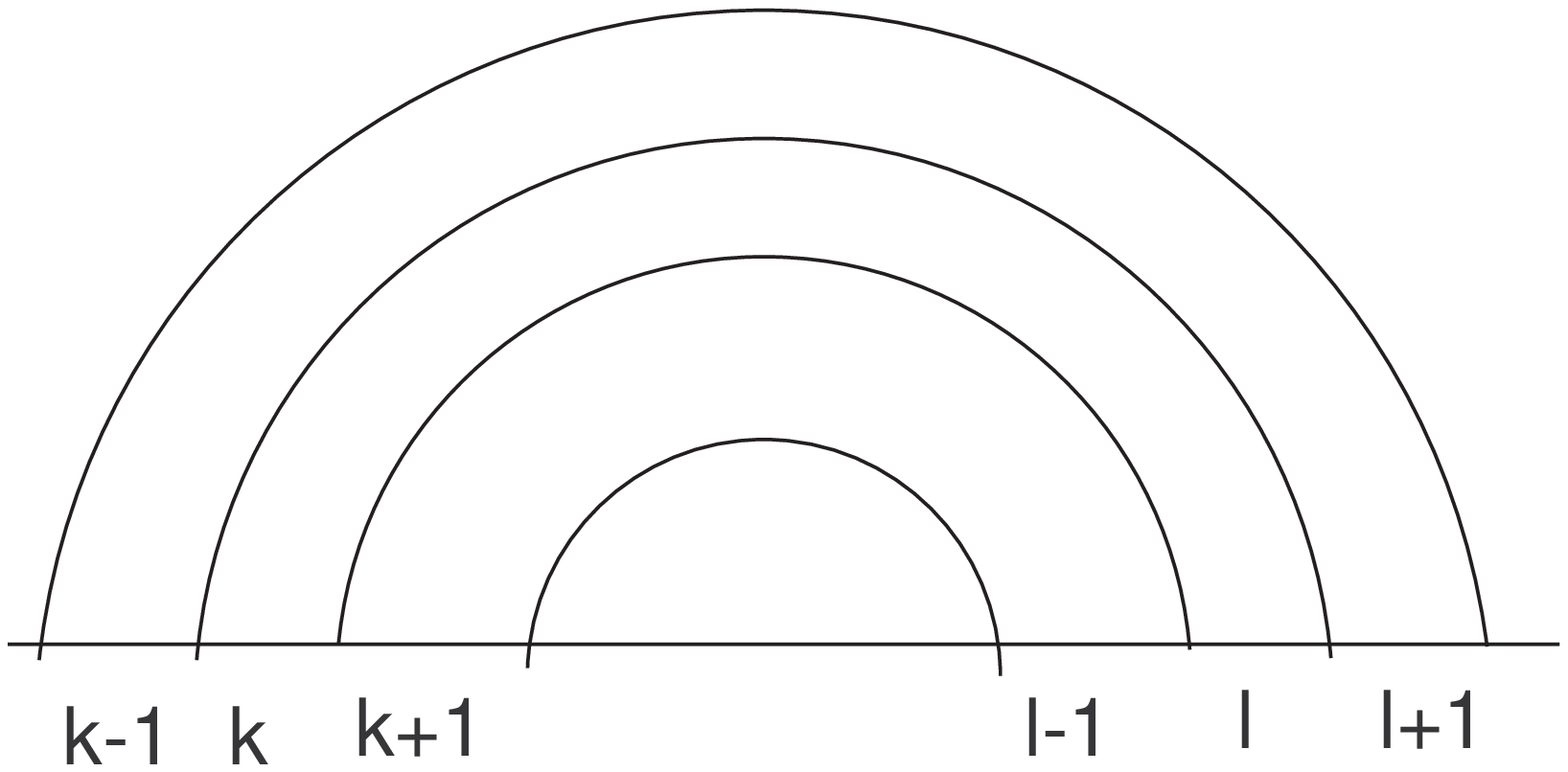}$  \\ \hline
  1 & 0 & 0 & 0 & $\psdiag{5}{10}{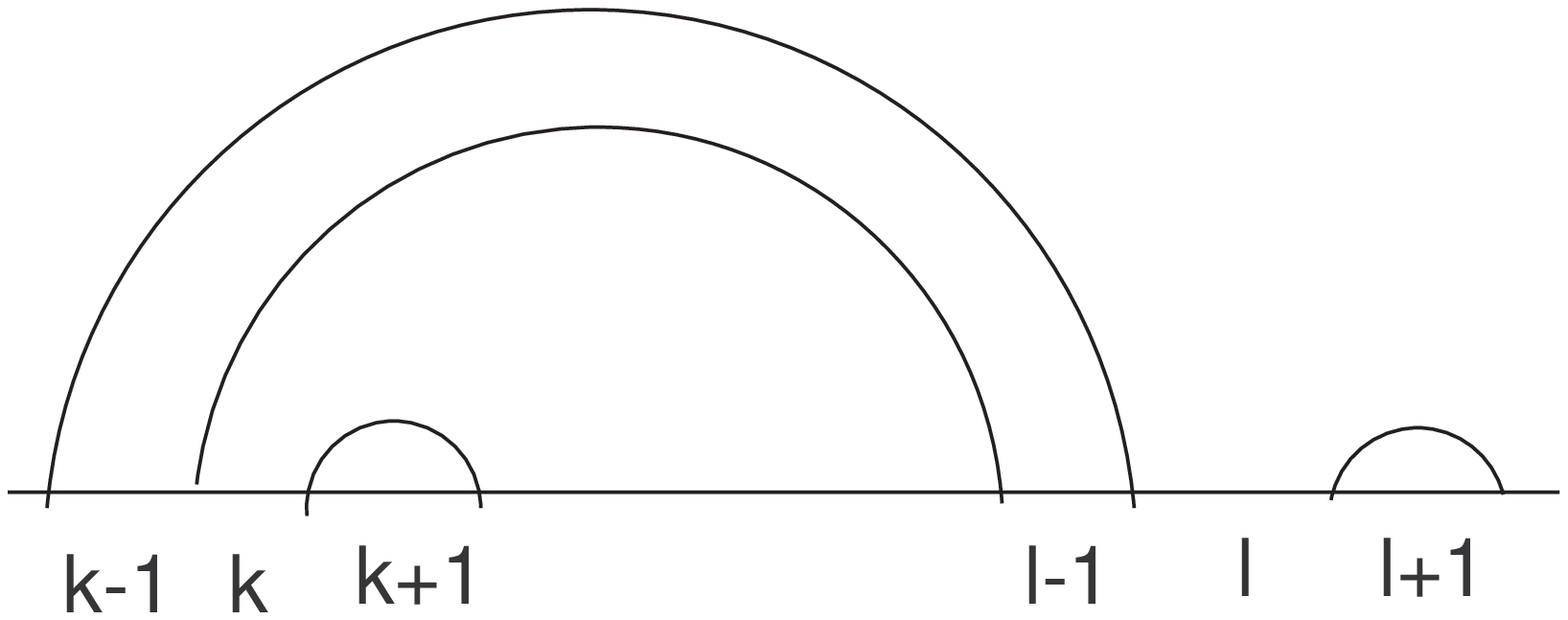}$  \\ \hline
  0 & 1 & 0 & 0 & $\psdiag{5}{10}{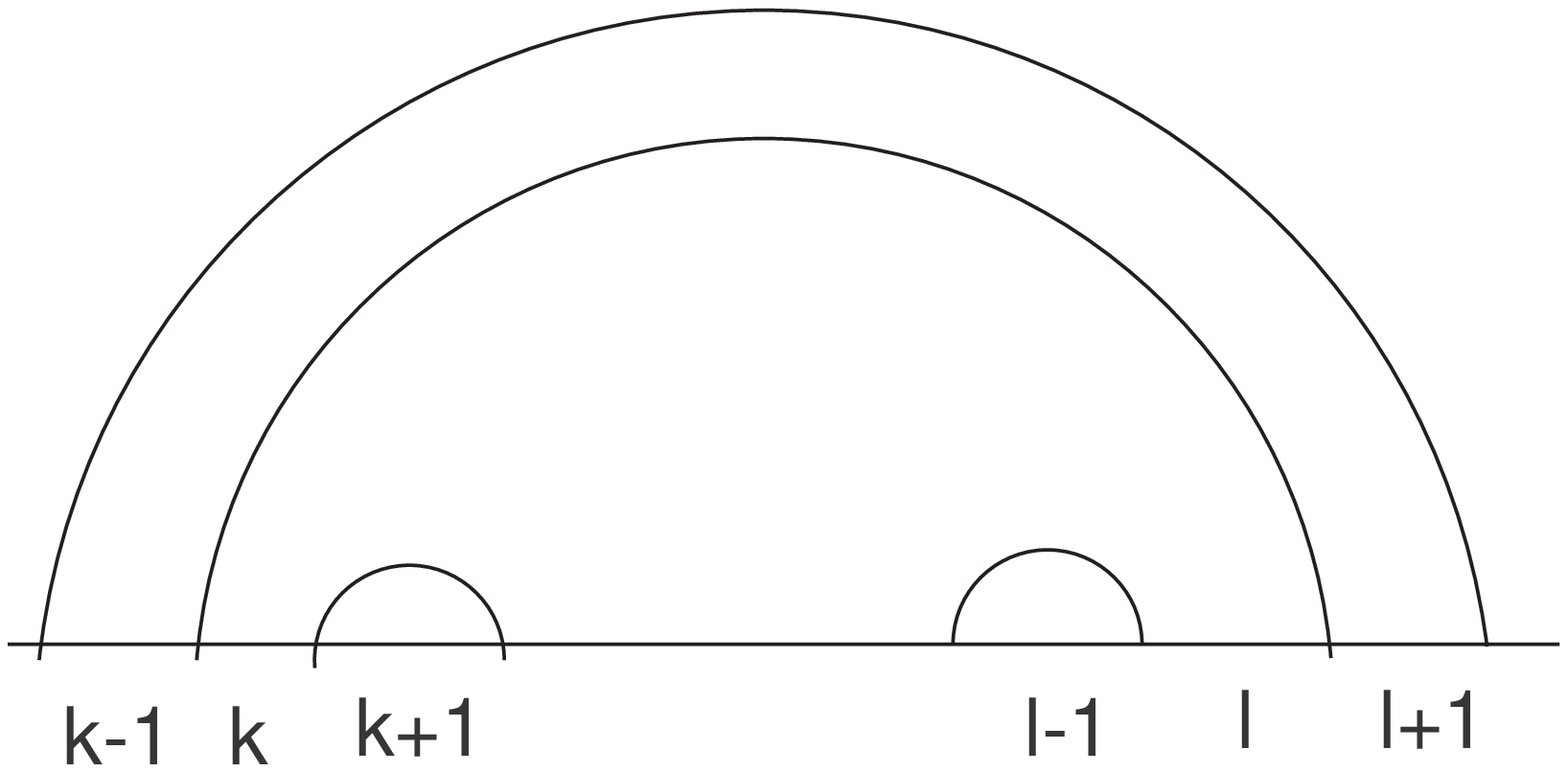}$  \\ \hline
  0 & 0 & 1 & 0 & $\psdiag{4}{8}{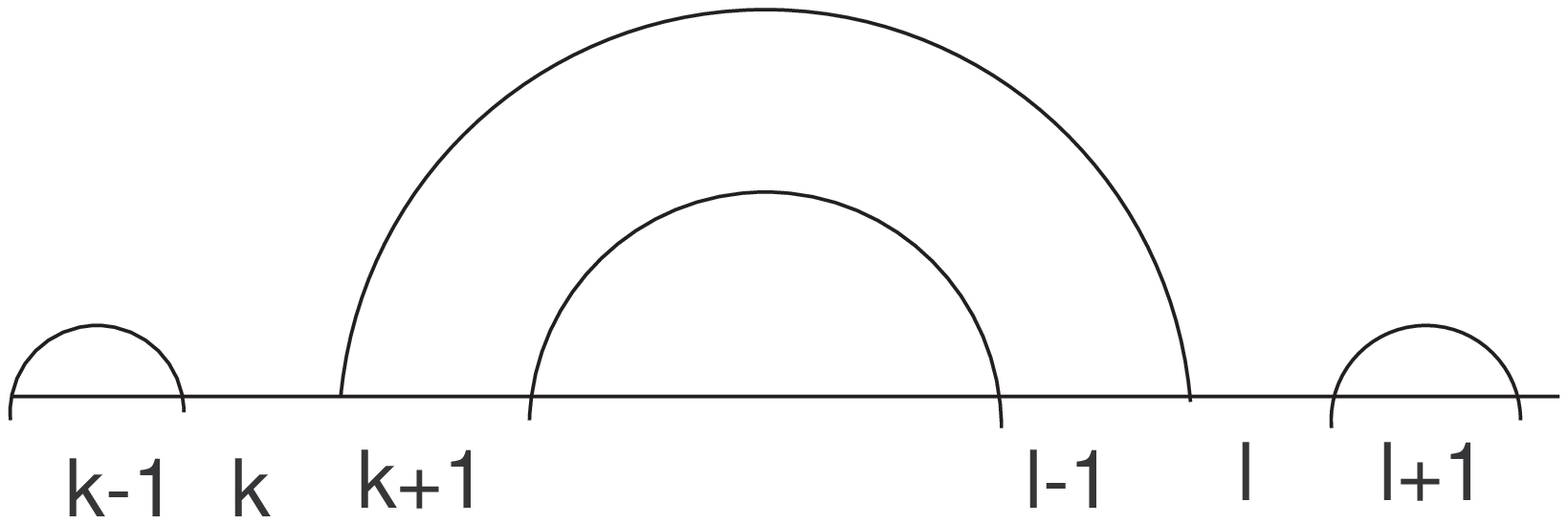}$  \\ \hline
  0 & 0 & 0 & 1 & $\psdiag{5}{10}{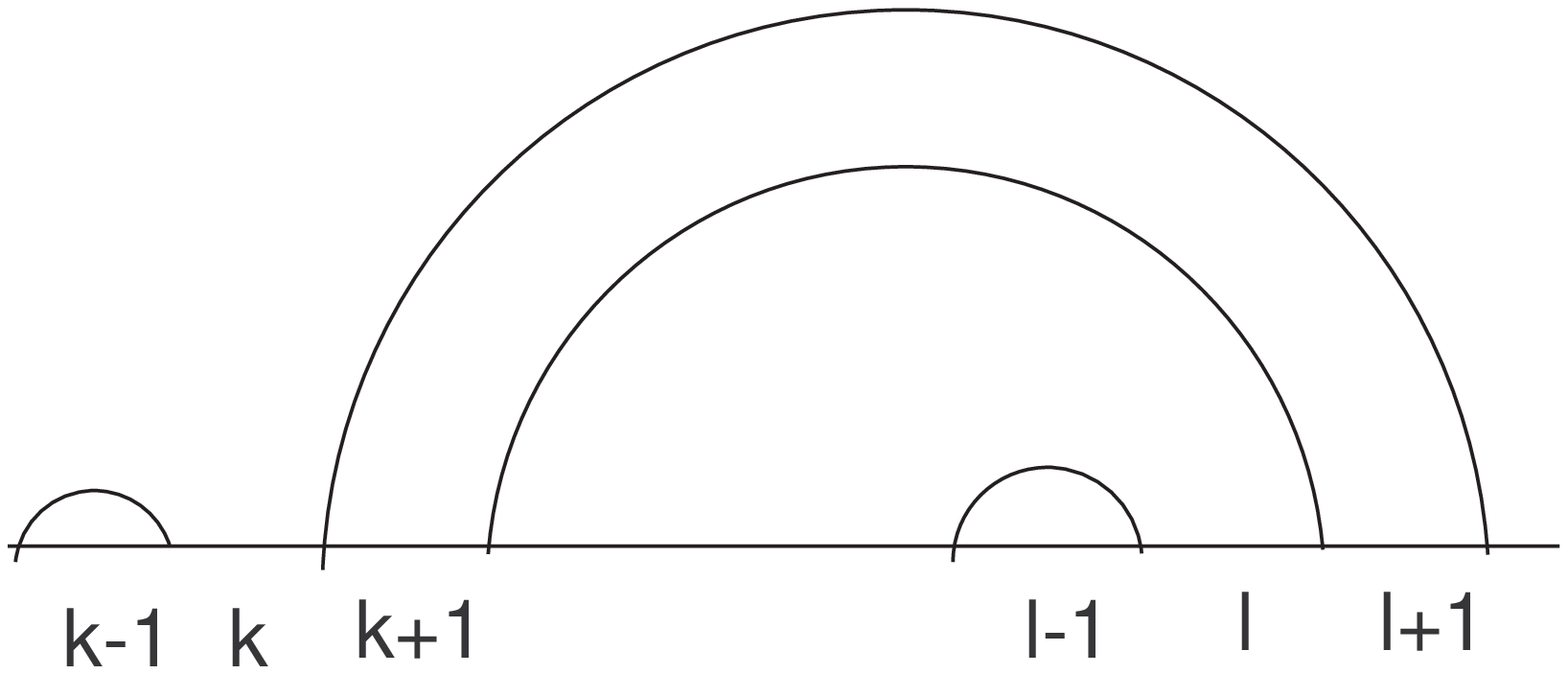}$  \\ \hline \hline
\end{tabular}

Using the properties of $R_1=[r_{i,j}]$ and $R_2=[\dot{r}_{i,j}]$,
and the relations between them ($R_2=(B_l^tR_1B_l)^2$), it is easy
to prove the following statements:

\begin{description}
\item[1.] If $r_{k-1,l-1}=1$ or $r_{k-1,l+1}=1$ or $r_{k+1,l-1}=1$
or $r_{k+1,l+1}=1$ then: $$r_{k-1,k+1}=1 \Leftrightarrow
(r_{k-1,l-1}=r_{k+1,l-1} \mbox{ and } r_{k-1,l+1}=r_{k+1,l+1});$$

\item[2.] If $r_{k-1,l-1}=1$ or $r_{k-1,l+1}=1$ or $r_{k+1,l-1}=1$
or $r_{k+1,l+1}=1$ then: $$r_{l-1,l+1}=1 \Leftrightarrow
(r_{k-1,l-1}=r_{k-1,l+1} \mbox{ and } r_{k+1,l-1}=r_{k+1,l+1});$$

\item[3.] $\dot{r}_{k-1,k+1}=r_{k-1,k+1}+r_{k-1,l+1}r_{k+1,l-1};$

\item[4.] $\dot{r}_{k-1,l-1}=r_{k-1,l-1}+r_{k-1,l+1};$

\item[5.] $\dot{r}_{k+1,l-1}=r_{k+1,l-1}+r_{k+1,l+1}.$
\end{description}

Also, using lemma \ref{11v}, we can easily check:

$$c_{k-1}=(b_{k-1}\vee b_{k+1})\oplus x_k= \left\{\begin{array}{ccc}
  b_{k-1}\oplus b_{k+1} & \mbox{if} & \dot{r}_{k-1,k+1}=0 \\
  b_{k-1}\oplus \varphi(b_k) & \mbox{if} & \dot{r}_{k-1,k+1}=1
\end{array}\right.$$

$$b_{k-1}= \left\{\begin{array}{ccc}
  a_{k-1} & \mbox{if} & \dot{r}_{k-1,l-1}=0 \\
  b_{l-1} & \mbox{if} & \dot{r}_{k-1,l-1}=1
\end{array}\right.$$

$$b_{k+1}= \left\{\begin{array}{ccc}
  a_{k+1} & \mbox{if} & \dot{r}_{k+1,l-1}=0 \\
  b_{l-1} & \mbox{if} & \dot{r}_{k+1,l-1}=1
\end{array}\right.$$

$$b_{l-1}=(a_{l-1}\vee a_{l+1})\oplus x_l=
\left\{\begin{array}{ccc}
  a_{l-1}\oplus a_{l+1} & \mbox{if} & r_{l-1,l+1}=0 \\
  a_{l-1}\oplus \varphi(a_l) & \mbox{if} & r_{l-1,l+1}=1
\end{array}\right.$$

Since $\ddot{r}_{k-1,l-3}=1 \Rightarrow k-1-(l-3)\in 2\mathbb{Z}$
we have $\dot{r}_{k,l-1}=0$ and then $b_k=a_k$.

With this, we can construct the following table:

{\tiny
\begin{tabular}{|c|c|c|c||c|c|c|}
\hline \hline
  $r_{k-1,l-1}$ & $r_{k-1,l+1}$ & $r_{k+1,l-1}$ & $r_{k+1,l+1}$ &
  $r_{k-1,k+1}$ & $\dot{r}_{k-1,k+1}$ & $c_{k-1}$
  \\ \hline \hline
  1 & 1 & 1 & 1 & 1 & 1 & $b_{k-1}\oplus \varphi(b_k)$
  \\ \hline
  1 & 0 & 1 & 0 & 1 & 1 & $b_{k-1}\oplus \varphi(b_k)$
  \\ \hline
  0 & 1 & 0 & 1 & 1 & 1 & $b_{k-1}\oplus \varphi(b_k)$
  \\ \hline
  1 & 1 & 0 & 0 & 0 & 0 & $b_{k-1}\oplus b_{k+1}$
  \\ \hline
  0 & 0 & 1 & 1 & 0 & 0 & $b_{k-1}\oplus b_{k+1}$
  \\ \hline
  0 & 1 & 1 & 0 & 0 & 1 & $b_{k-1}\oplus \varphi(b_k)$
  \\ \hline
  1 & 0 & 0 & 0 & 0 & 0 & $b_{k-1}\oplus b_{k+1}$
  \\ \hline
  0 & 1 & 0 & 0 & 0 & 0 & $b_{k-1}\oplus b_{k+1}$ 
  \\ \hline
  0 & 0 & 1 & 0 & 0 & 0 & $b_{k-1}\oplus b_{k+1}$ 
  \\ \hline
  0 & 0 & 0 & 1 & 0 & 0 & $b_{k-1}\oplus b_{k+1}$
  \\ \hline \hline
\end{tabular}
}

{\tiny
\begin{tabular}{|c|c|c|c|c|c|c|}
\hline \hline
  $r_{l-1,l+1}$ &
  $b_{l-1}$ & $\dot{r}_{k-1,l-1}$ & $b_{k-1}$ & $\dot{r}_{k+1,l-1}$ &
  $b_{k+1}$ & $c_{k-1}$ \\ \hline \hline
  1 &
  $a_{l-1}\oplus \varphi(a_l)$ & 1 & $a_{l-1}\oplus \varphi(a_l)$ &
  1 & $a_{l-1}\oplus \varphi(a_l)$ &
  $[a_{l-1}\oplus \varphi(a_l)]\oplus\varphi(a_k)$ \\ \hline
  0 &
  $a_{l-1}\oplus a_{l+1}$ & 1 & $a_{l-1}\oplus a_{l+1}$ &
  1 & $a_{l-1}\oplus a_{l+1}$ &
  $[a_{l-1}\oplus a_{l+1}]\oplus\varphi(a_k)$ \\ \hline
  0 &
  $a_{l-1}\oplus a_{l+1}$ & 1 & $a_{l-1}\oplus a_{l+1}$ &
  1 & $a_{l-1}\oplus a_{l+1}$ &
  $[a_{l-1}\oplus a_{l+1}]\oplus\varphi(a_k)$ \\ \hline
  1 &
  $a_{l-1}\oplus \varphi(a_l)$ & 1 & $a_{l-1}\oplus \varphi(a_l)$ &
  0 & $a_{k+1}$ &
  $[a_{l-1}\oplus \varphi(a_l)]\oplus a_{k+1}$ \\ \hline
  1 &
  $a_{l-1}\oplus \varphi(a_l)$ & 0 & $a_{k-1}$ &
  1 & $a_{l-1}\oplus \varphi(a_l)$ &
  $a_{k-1}\oplus[a_{l-1}\oplus \varphi(a_l)]$ \\ \hline
  0 &
  $a_{l-1}\oplus a_{l+1}$ & 1 & $a_{l-1}\oplus a_{l+1}$ &
  1 & $a_{l-1}\oplus a_{l+1}$ &
  $[a_{l-1}\oplus a_{l+1}]\oplus\varphi(a_k)$ \\ \hline
  0 &
  $a_{l-1}\oplus a_{l+1}$ & 1 & $a_{l-1}\oplus a_{l+1}$ &
  0 & $a_{k+1}$ &
  $[a_{l-1}\oplus a_{l+1}]\oplus a_{k+1}$ \\ \hline
  0 &
  $a_{l-1}\oplus a_{l+1}$ & 1 & $a_{l-1}\oplus a_{l+1}$ &
  0 & $a_{k+1}$ &
  $[a_{l-1}\oplus a_{l+1}]\oplus a_{k+1}$ \\ \hline
  0 &
  $a_{l-1}\oplus a_{l+1}$ & 0 & $a_{k-1}$ &
  1 & $a_{l-1}\oplus a_{l+1}$ &
  $a_{k-1}\oplus[a_{l-1}\oplus a_{l+1}]$ \\ \hline
  0 &
  $a_{l-1}\oplus a_{l+1}$ & 0 & $a_{k-1}$ &
  1 & $a_{l-1}\oplus a_{l+1}$ &
  $a_{k-1}\oplus[a_{l-1}\oplus a_{l+1}]$ \\ \hline \hline
\end{tabular}
}

We can construct an analogous table for the value $c'_{k-1}$. All
we need to know is that:

\begin{description}
\item[1.] If $r_{k-1,l-1}=1$ or $r_{k-1,l+1}=1$ or $r_{k+1,l-1}=1$
or $r_{k+1,l+1}=1$ then: $$r_{k-1,k+1}=1 \Leftrightarrow
(r_{k-1,l-1}=r_{k+1,l-1} \mbox{ and } r_{k-1,l+1}=r_{k+1,l+1});$$

\item[2.] If $r_{k-1,l-1}=1$ or $r_{k-1,l+1}=1$ or $r_{k+1,l-1}=1$
or $r_{k+1,l+1}=1$ then: $$r_{l-1,l+1}=1 \Leftrightarrow
(r_{k-1,l-1}=r_{k-1,l+1} \mbox{ and } r_{k+1,l-1}=r_{k+1,l+1});$$

\item[3.] $\dot{r}'_{l-3,l-1}=r_{l-1,l+1}+r_{k-1,l+1}r_{k+1,l-1};$

\item[4.] $\dot{r}'_{k-1,l-3}=r_{k-1,l-1}+r_{k+1,l-1};$

\item[5.] $\dot{r}'_{k-1,l-1}=r_{k-1,l+1}+r_{k+1,l+1}.$
\end{description}

$$\ddot{r}'_{k-1,l-3}=\ddot{r}_{k-1,l-3}=1 \Rightarrow$$
$$c'_{k-1}=c'_{l-3}=(b'_{l-3}\vee b'_{l-1})\oplus x'_{l-2}=
\left\{\begin{array}{ccc}
  b'_{l-3}\oplus b'_{l-1} & \mbox{if} & \dot{r}'_{l-3,l-1}=0 \\
  b'_{l-3}\oplus \varphi(b'_{l-2}) & \mbox{if} & \dot{r}'_{l-3,l-1}=1
\end{array}\right.$$

$$b'_{l-3}= \left\{\begin{array}{ccc}
  a_{l-1} & \mbox{if} & \dot{r}'_{k-1,l-3}=0 \\
  b'_{k-1} & \mbox{if} & \dot{r}'_{k-1,l-3}=1
\end{array}\right.$$

$$b'_{l-1}= \left\{\begin{array}{ccc}
  a_{l+1} & \mbox{if} & \dot{r}'_{k-1,l-1}=0 \\
  b'_{k-1} & \mbox{if} & \dot{r}'_{k-1,l-1}=1
\end{array}\right.$$

$$b'_{k-1}=(a_{k-1}\vee a_{k+1})\oplus x'_k=
\left\{\begin{array}{ccc}
  a_{k-1}\oplus a_{k+1} & \mbox{if} & r_{k-1,k+1}=0 \\
  a_{k-1}\oplus \varphi(a_k) & \mbox{if} & r_{k-1,k+1}=1
\end{array}\right.$$
and $b'_{l-2}=a_l$ since $\ddot{r}'_{k-1,l-3}=1 \Rightarrow
\dot{r}'_{k-1,l-2}=0$.

{\tiny
\begin{tabular}{|c|c|c|c||c|c|c|}
\hline \hline
  $r_{k-1,l-1}$ & $r_{k-1,l+1}$ & $r_{k+1,l-1}$ & $r_{k+1,l+1}$ &
  $r_{l-1,l+1}$ & $\dot{r}'_{l-3,l-1}$ & $c'_{k-1}$ 
  \\ \hline \hline
  1 & 1 & 1 & 1 & 1 & 1 & $b'_{l-3}\oplus \varphi(b'_{l-2})$ 
  \\ \hline
  1 & 0 & 1 & 0 & 0 & 0 & $b'_{l-3}\oplus b'_{l-1}$
  \\ \hline
  0 & 1 & 0 & 1 & 0 & 0 & $b'_{l-3}\oplus b'_{l-1}$
  \\ \hline
  1 & 1 & 0 & 0 & 1 & 1 & $b'_{l-3}\oplus \varphi(b'_{l-2})$
  \\ \hline
  0 & 0 & 1 & 1 & 1 & 1 & $b'_{l-3}\oplus \varphi(b'_{l-2})$
  \\ \hline
  0 & 1 & 1 & 0 & 0 & 1 & $b'_{l-3}\oplus \varphi(b'_{l-2})$
  \\ \hline
  1 & 0 & 0 & 0 & 0 & 0 & $b'_{l-3}\oplus b'_{l-1}$
  \\ \hline
  0 & 1 & 0 & 0 & 0 & 0 & $b'_{l-3}\oplus b'_{l-1}$
  \\ \hline
  0 & 0 & 1 & 0 & 0 & 0 & $b'_{l-3}\oplus b'_{l-1}$
  \\ \hline
  0 & 0 & 0 & 1 & 0 & 0 & $b'_{l-3}\oplus b'_{l-1}$
  \\ \hline
 \hline
\end{tabular}
}

{\tiny
\begin{tabular}{|c|c|c|c|c|c|c|}
\hline \hline
  $r_{k-1,k+1}$ &
  $b'_{k-1}$ & $\dot{r}'_{l-3,k-1}$ & $b'_{l-3}$ & $\dot{r}'_{k-1,l-1}$ &
  $b'_{l-1}$ & $c'_{k-1}$ \\ \hline \hline
  1 &
  $a_{k-1}\oplus \varphi(a_k)$ & 1 & $a_{k-1}\oplus \varphi(a_k)$ &
  1 & $a_{k-1}\oplus \varphi(a_k)$ &
  $[a_{k-1}\oplus \varphi(a_k)]\oplus\varphi(a_l)$ \\ \hline
  1 &
  $a_{k-1}\oplus \varphi(a_k)$ & 1 & $a_{k-1}\oplus \varphi(a_k)$ &
  0 & $a_{l+1}$ &
  $[a_{k-1}\oplus \varphi(a_k)]\oplus a_{l+1}$ \\ \hline
  1 &
  $a_{k-1}\oplus \varphi(a_k)$ & 0 & $a_{l-1}$ &
  1 & $a_{k-1}\oplus \varphi(a_k)$ &
  $a_{l-1}\oplus[a_{k-1}\oplus \varphi(a_k)]$ \\ \hline
  0 &
  $a_{k-1}\oplus a_{k+1}$ & 1 & $a_{k-1}\oplus a_{k+1}$ &
  1 & $a_{k-1}\oplus a_{k+1}$ &
  $[a_{k-1}\oplus a_{k+1}]\oplus\varphi(a_l)$ \\ \hline
  0 &
  $a_{k-1}\oplus a_{k+1}$ & 1 & $a_{k-1}\oplus a_{k+1}$ &
  1 & $a_{k-1}\oplus a_{k+1}$ &
  $[a_{k-1}\oplus a_{k+1}]\oplus\varphi(a_l)$ \\ \hline
  0 &
  $a_{k-1}\oplus a_{k+1}$ & 1 & $a_{k-1}\oplus a_{k+1}$ &
  1 & $a_{k-1}\oplus a_{k+1}$ &
  $[a_{k-1}\oplus a_{k+1}]\oplus\varphi(a_l)$ \\ \hline
  0 &
  $a_{k-1}\oplus a_{k+1}$ & 1 & $a_{k-1}\oplus a_{k+1}$ &
  0 & $a_{l+1}$ &
  $[a_{k-1}\oplus a_{k+1}]\oplus a_{l+1}$ \\ \hline
  0 &
  $a_{k-1}\oplus a_{k+1}$ & 0 & $a_{l-1}$ &
  1 & $a_{k-1}\oplus a_{k+1}$ &
  $a_{l-1}\oplus [a_{k-1}\oplus a_{k+1}]$ \\ \hline
  0 &
  $a_{k-1}\oplus a_{k+1}$ & 1 & $a_{k-1}\oplus a_{k+1}$ &
  0 & $a_{l+1}$ &
  $[a_{k-1}\oplus a_{k+1}]\oplus a_{l+1}$ \\ \hline
  0 &
  $a_{k-1}\oplus a_{k+1}$ & 0 & $a_{l-1}$ &
  1 & $a_{k-1}\oplus a_{k+1}$ &
  $a_{l-1}\oplus [a_{k-1}\oplus a_{k+1}]$ \\ \hline
 \hline
\end{tabular}
}

In order to facilitate the comparison of $c_{k-1}$ and $c'_{k-1}$
we will substitute each $a_i$ (with $i\in \{k+1,l-1,l+1\}$)
appearing in the expressions of $c_{k-1}$ and $c'_{k-1}$ by $a_j$
where $j\in \{k-1,k+1,l-1,l+1\}$ is the smallest index such that
$r_{i,j}=1$.

The results can be seen in the following table:


\begin{tabular}{|c|c|c|c||c|c|}
\hline \hline
  $r_{k-1,l-1}$ & $r_{k-1,l+1}$ & $r_{k+1,l-1}$ & $r_{k+1,l+1}$ &
  $r_{k-1,k+1}$ & $r_{l-1,l+1}$ 
  \\ \hline \hline
  1 & 1 & 1 & 1 & 1 & 1 
  \\ \hline
  1 & 0 & 1 & 0 & 1 & 0 
  \\ \hline
  0 & 1 & 0 & 1 & 1 & 0 
  \\ \hline
  1 & 1 & 0 & 0 & 0 & 1 
  \\ \hline
  0 & 0 & 1 & 1 & 0 & 1 
  \\ \hline
  0 & 1 & 1 & 0 & 0 & 0 
  \\ \hline
  1 & 0 & 0 & 0 & 0 & 0 
  \\ \hline
  0 & 1 & 0 & 0 & 0 & 0 
  \\ \hline
  0 & 0 & 1 & 0 & 0 & 0 
  \\ \hline
  0 & 0 & 0 & 1 & 0 & 0 
  \\ \hline \hline
\end{tabular}

\begin{tabular}{|c|c|c|c|c|}
\hline \hline
  $a_{k+1}$ & $a_{l-1}$ & $a_{l+1}$ &
  $c_{k-1}$ & $c'_{k-1}$ \\ \hline \hline
  $a_{k-1}$ & $a_{k-1}$ & $a_{k-1}$ &
  $a_{k-1}\oplus \varphi(a_l) \oplus \varphi(a_k)$ &
  $a_{k-1}\oplus \varphi(a_k) \oplus \varphi(a_l)$ \\ \hline
  $a_{k-1}$ & $a_{k-1}$ & $a_{l+1}$ &
  $a_{k-1}\oplus a_{l+1} \oplus \varphi(a_k)$ &
  $a_{k-1}\oplus \varphi(a_k) \oplus a_{l+1}$  \\ \hline
  $a_{k-1}$ & $a_{l-1}$ & $a_{k-1}$ &
  $a_{l-1}\oplus a_{k-1} \oplus \varphi(a_k)$ &
  $a_{l-1}\oplus a_{k-1} \oplus \varphi(a_k)$  \\ \hline
  $a_{k+1}$ & $a_{k-1}$ & $a_{k-1}$ &
  $a_{k-1}\oplus \varphi(a_l) \oplus a_{k+1}$ &
  $a_{k-1}\oplus a_{k+1} \oplus \varphi(a_l)$  \\ \hline
  $a_{k+1}$ & $a_{k+1}$ & $a_{k+1}$ &
  $a_{k-1}\oplus a_{k+1} \oplus \varphi(a_l)$ &
  $a_{k-1}\oplus a_{k+1} \oplus \varphi(a_l)$  \\ \hline
  $a_{k+1}$ & $a_{k+1}$ & $a_{k-1}$ &
  $a_{k+1}\oplus a_{k-1} \oplus \varphi(a_k)$ &
  $a_{k-1}\oplus a_{k+1} \oplus \varphi(a_l)$  \\ \hline
  $a_{k+1}$ & $a_{k-1}$ & $a_{l+1}$ &
  $a_{k-1}\oplus a_{l+1} \oplus a_{k+1}$ &
  $a_{k-1}\oplus a_{k+1} \oplus a_{l+1}$  \\ \hline
  $a_{k+1}$ & $a_{l-1}$ & $a_{k-1}$ &
  $a_{l-1}\oplus a_{k-1} \oplus a_{k+1}$ &
  $a_{l-1}\oplus a_{k-1} \oplus a_{k+1}$  \\ \hline
  $a_{k+1}$ & $a_{k+1}$ & $a_{l+1}$ &
  $a_{k-1}\oplus a_{k+1} \oplus a_{l+1}$ &
  $a_{k-1}\oplus a_{k+1} \oplus a_{l+1}$  \\ \hline
  $a_{k+1}$ & $a_{l-1}$ & $a_{k+1}$ &
  $a_{k-1}\oplus a_{l-1} \oplus a_{k+1}$ &
  $a_{l-1}\oplus a_{k-1} \oplus a_{k+1}$  \\ \hline \hline
\end{tabular}


We can easily see that $c_{k-1}=c'_{k-1}$ for all rows except the
sixth row where $c_{k-1}=c'_{k-1}$ if $a_k=a_l$. But, since
$r_{k-1,l+1}=r_{k+1,l-1}=1$ and $r_{k-1,k+1}=0$ in this case, we
have, by the topological properties T2 and T3, that $r_{k,l}=1$
and then $a_k=a_l$.

Now, to see that $c_{l-3}=c'_{l-3}$ we proceed in the same way as
we did to show $c_{k-1}=c'_{k-1}$ for the case
$\ddot{r}_{k-1,l-3}=0$. For the case $\ddot{r}_{k-1,l-3}=1$, we
have $c_{l-3}=c_{k-1}=c'_{k-1}=c'_{l-3}$ (by lemma 11).

For another generic index $i$, we have:

\begin{quote}
If $\ddot{r}_{i,k-1}=1$ then $c_i=c_{k-1}=c'_{k-1}=c'_i$.
\end{quote}

\begin{quote}
If $\ddot{r}_{i,l-3}=1$ then $c_i=c_{l-3}=c'_{l-3}=c'_i$.
\end{quote}

\begin{quote}
If $\ddot{r}_{i,k-1}=\ddot{r}_{i,l-3}=0$ then we have:
\end{quote}

$$c_i=e_i^tR_3*[(B_k^t*\overrightarrow{b})
\oplus(e_{k-1}*x_k)]=e_i^tB_k^t*\overrightarrow{b}=
b_{\check{k}(i)}$$ $$\mbox{where} \quad
\check{k}(i)=\left\{\begin{array}{ccc}
  i & \mbox{if} & i<k-1 \\
  i-2 & \mbox{if} & i>k-1
\end{array}\right. $$

$\ddot{r}_{i,l-3}=0 \Rightarrow \dot{r}_{\check{k}(i),l-1}=0$, by
lemma 11. Thus $$c_i=b_{\check{k}(i)}=
e_{\check{k}(i)}^tR_2*[(B_l^t*\overrightarrow{a})
\oplus(e_{l-1}*x_l)]\}=e_{\check{k}(i)}^tB_l^t*\overrightarrow{a}=
e_i^tB_k^tB_l^t*\overrightarrow{a}$$

$$c'_i=e_i^tR'_3*[(B_{l-2}^t*\overrightarrow{b}')
\oplus(e_{l-3}*x'_{l-2})]\}=e_i^tB_{l-2}^t*\overrightarrow{b}'=
b'_{\check{l-2}(i)}$$

$\ddot{r}'_{i,k-1}=\ddot{r}_{i,k-1}=0 \Rightarrow
\dot{r}'_{\check{l-2}(i),k-1}=0$, by lemma \ref{11v}. Thus
$$c'_i=b'_{\check{l-2}(i)}=
e_{\check{l-2}(i)}^tR'_2*[(B_k^t*\overrightarrow{a})
\oplus(e_{k-1}*x'_k)]\}=e_{\check{l-2}(i)}^tB_k^t*\overrightarrow{a}=
e_i^tB_{l-2}^tB_k^t*\overrightarrow{a}$$

Since $B_{l-2}^tB_k^t=B_k^tB_l^t$ we have $c_i=c'_i$.

\end{description}

\section{Systems of non-singular planar curves}


Now, we are going to study this representation in the particular
case of systems of non-singular planar curves. We mean by a system
of non-singular planar curves an immersion in the plane of a
finite number of disjoint circles which may be regarded as a
morphism from the empty partition of the line to itself (i.e. an
element of $\hom(\emptyset,\emptyset)$).

Any morphism $t\in\hom(\emptyset,\emptyset)$ is a word made of
generators $\hat{t}_{n,k}$ and $\check{t}_{n,k}$. Let us make the
substitution: $\hat{t}_{n,k}\mapsto(2,2k-n-3)$ and
$\check{t}_{n,k}\mapsto(-2,2k-n-3)$ for each word. The number
$2k-n-3$ counts the number of strings on the left of the local
maximum (minimum) minus the number of strings on the right.

The following lemma shows that no information is lost after this
substitution. But first, let us introduce a notation for the
following sets: $$S_+:=\{(2,n):n\in 2\mathbb{Z}\}$$ and
$$S_-:=\{(-2,n):n\in 2\mathbb{Z}\}$$

\begin{lemma}
After the substitution: $\hat{t}_{n,k}\mapsto(2,2k-n-3)$ and
$\check{t}_{n,k}\mapsto(-2,2k-n-3)$ of a morphism
$t\in\hom(\emptyset,\emptyset)$ we get a word
$(c_1,d_1)(c_2,d_2)...(c_n,d_n)$ in $S_+\cup S_-$ satisfying the
following condition:
\begin{description}
\item[C.] For any index $i$:
\begin{description}
\item[] if $c_i=2$ then $|d_i|\leq
-(\sum_{j<i}c_j)-2 =\sum_{j>i}c_j$,
\item[] if $c_i=-2$ then $|d_i|\leq
-\sum_{j<i}c_j =(\sum_{j>i}c_j)-2$.
\end{description}
\end{description}
Also, if a word in $S_+\cup S_-$ satisfies this condition then
there exists a unique morphism $t\in\hom(\emptyset,\emptyset)$
that becomes this word after the substitution:
$\hat{t}_{n,k}\mapsto(2,2k-n-3)$ and
$\check{t}_{n,k}\mapsto(-2,2k-n-3)$.
\end{lemma}
\TeXButton{Proof}{\proof} Observing that each $(c_i,d_i)$
substitutes one generator $\hat{t}_{n_i,k_i}$ or
$\check{t}_{n_i,k_i}$ it easy to see that $-\sum_{j<i}c_j$ is the
number of points at the bottom of the generator
($\hat{t}_{n_i,k_i}$ or $\check{t}_{n_i,k_i}$) and that
$\sum_{j>i}c_j$ is number of points at the top of the generator.
The condition follows naturally from this fact.

For the second part of the lemma, we make the inverse substitution
$(c_i,d_i)=(2,k)\mapsto \hat{t}_{n',k'}$ with
$n'=(\sum_{j>i}c_j)+1$ and $k'=\frac{k+n'+3}{2}$, and
$(c_i,d_i)=(-2,k)\mapsto \check{t}_{n',k'}$ with
$n'=(\sum_{j>i}c_j)-1$ and $k'=\frac{k+n'+3}{2}$. \TeXButton{End
Proof}{\endproof}

So we have a new language for morphisms in
$\hom(\emptyset,\emptyset)$ and in this language the local
relations become:
\begin{description}
\item[1] $...(c_i,d_i)(-2,k)(2,k+2)(c_{i+3},d_{i+3})...
=...(c_i,d_i) (c_{i+3},d_{i+3})...
=...(c_i,d_i)(-2,k)(2,k-2)(c_{i+3},d_{i+3})...$ (i.e.
$\check{t}_{k'+1}\circ\hat{t}_{k'}= id=
\check{t}_{k'-1}\circ\hat{t}_{k'}$);
\item[2] $...(2,k)(2,l)...=...(2,l+2)(2,k+2)...$
for $k\leq l-2$ (i.e. $\hat{t}_{k'}\circ\hat{t}_{l'}=
\hat{t}_{l'+2}\circ\hat{t}_{k'}$ for $l'\geq k'+2$);
\item[3.1]
$...(-2,k)(2,l)...=...(2,l-2)(-2,k+2)...$ for $k\leq l-4$ (i.e.
$\check{t}_{k'}\circ\hat{t}_{l'}=
\hat{t}_{l'-2}\circ\check{t}_{k'}$ for $l'\geq k'+2$)
\item[3.2] $...(2,k)(-2,l)...=...(-2,l+2)(2,k-2)...$ for
$k\leq l$ (i.e. $\hat{t}_{k'}\circ\check{t}_{l'}=
\check{t}_{l'+2}\circ\hat{t}_{k'}$ for $l'\geq k'+2$)
\item[4] $...(-2,k)(-2,l)...=...(-2,l-2)(-2,k-2)$ for $k\leq l-2$
(i.e. $\check{t}_{l'-2}\circ\check{t}_{k'}=
\check{T}_{k'}\circ\check{T}_{l'}$ for $l'\geq k'+2$).
\end{description}

Note that, by the previous lemma, these relations preserve the
condition C.


\begin{proposition}
Any word in $S_+\cup S_-$ satisfying the condition C is
equivalent, by the previous relations, to a word of symbols
$(2,0)$ and $(-2,0)$.
\end{proposition}
\TeXButton{Proof}{\proof} We use the following algorithm to
transform any word in $S_+\cup S_-$ satisfying the condition C
into a word of symbols $(2,0)$ and $(-2,0)$.

ALGORITHM:
\begin{description}
\item[input:] Take a word in $S_+\cup S_-$ satisfying the
condition C;
\item[step 1.] Apply the relations 3.1 ($...(2,l-2)(-2,k+2)...=
...(-2,k)(2,l)...$ for $k\leq l-4$) and 3.2 ($...(2,k)(-2,l)...=
...(-2,l+2)(2,k-2)...$ for $k\leq l$) so as to put all the symbols
of the form $(-2,k)$ on the left side of the word and all the
symbols of the form $(2,k)$ on the right;
\item[step 2.] Use the relations 2 ($...(2,k)(2,l)...=...(2,l+2)(2,k+2)...$
for $k\leq l-2$) and 4 ($...(-2,k)(-2,l)...=...(-2,l-2)(-2,k-2)$
for $k\leq l-2$) to order in semi-decreasing order (according to
$d_i$) each block consisting of $(c_i,d_i)$ of the form $(2,k)$ or
$(-2,k)$;
\item[step 3.] If there exists a sequence of the type
$(-2,k)(2,k+2)$ in the word we apply the relation 1
($...(c_i,d_i)(-2,k)(2,k+2)(c_{i+3},d_{i+3})... =...(c_i,d_i)
(c_{i+3},d_{i+3})...$) and go back to step 1;
\item[step 4.] If there exists a sequence of the type
$(-2,k)(2,l)$ with $k\leq l-4$ then we apply the relation 3.1
($...(-2,k)(2,l)...=...(2,l-2)(-2,k+2)...$) and go back to the
step 3;
\item[step 5.] If there doesn't exist a sequence of the type
$(-2,k)(2,l)$ with $k\leq l-2$ in the word we {\bf output} this.
\end{description}

\begin{description}
\item[Claim 1.] The algorithm always terminates in a finite number of
steps.
\end{description}

We will prove this by induction on the number of symbols since the
algorithm doesn't increase this.

It easy to see that any word satisfying the condition C has the
same number of symbols of the form $(2,k)$ as of the form
$(-2,k)$, and always begins with $(-2,0)$ and ends with $(2,0)$.
So $(-2,0)(2,0)$ is the unique word of two symbols and this passes
through the algorithm without any changes.

Now, assuming that the algorithm terminates for any word with $2n$
symbols, we take a word with $2n+2$ symbols. We consider for any
word $w=(c_1,d_1)(c_2,d_2)...(c_n,d_n)$ the "potential"
$$E(w)=(\sum_{c_i=2}i,\sum_{c_i=-2}d_i -\sum_{c_i=2}d_i)$$ which
takes values in $\mathbb{Z}^2$ with lexicographic order (i.e.
$(a,b)\leq (c,d)$ iff $a\leq c$ or $a=c$ and $b\leq d$). We can
see that, after being increased in step 1, the potential of the
word is always decreased until (if it occurs) the algorithm
returns to step 1 (after step 3), but in this case the word has
$2n$ symbols and then by the induction hypothesis the algorithm
terminates. The condition C implies that there are a finite number
of potentials for a word with $2n+2$ symbols (or less), and thus
the algorithm terminates in a finite number of steps.

\begin{description}
\item[Claim 2.] The output word has only $(2,0)$ and $(-2,0)$ as
symbols.
\end{description}

For a word $w=(c_1,d_1)...(c_{2n},d_{2n})$ let
$\alpha_1<...<\alpha_n$ be the indices such that $c_{\alpha_i}=2$
and $\beta_1<...<\beta_n$ be the indices such that
$c_{\beta_i}=-2$. We consider a new "potential" $E_2$ defined by
the formula: $$E_2(w)=\max\left(\{d_{\alpha_{i+1}}-d_{\alpha_i}
-2(\alpha_{i+1}-\alpha_i -1) : 1\leq i\leq n-1\}\right. \cup$$
$$\quad \left. \{d_{\beta_{i+1}}-d_{\beta_i}
-2(\beta_{i+1}-\beta_i -1) : 1\leq i\leq n-1\}\right)$$

After step 2 the potential of the word becomes non-positive and
step 4 doesn't change this. This means that the output word
contains no sequence of the type $(2,k)(2,l)$ or $(-2,k)(-2,l)$
with $k<l$. It is a condition of the algorithm that the output
word contains no sequence of the type $(-2,k)(2,l)$ with $k<l$.

Thus, if it does not contain a sequence of the type $(2,k)(-2,l)$
with $k<l$, then the word is ordered by semi-decreasing order and,
since it has to begin with $(-2,0)$ and to end with $(2,0)$, the
word has only $(2,0)$ and $(-2,0)$ as symbols.

To show that the output word contains no sequence of the type
$(2,k)(-2,l)$ with $k<l$, it is enough to observe that after step
1 there is no sequence of the type $(2,k)(-2,l)$ and that steps 2
and 4 do not produce any new sequences of the type $(2,k)(-2,l)$
with $k<l$. \TeXButton{End Proof}{\endproof}

\begin{quote}
NOTE: This algorithm was not conceived to be the most efficient
but to guarantee an easy proof that it terminates. The author
conjectures that there exist more efficient algorithms.
\end{quote}

Next, we observe that the family of non-singular planar curves
(i.e. $\hom(\emptyset,\emptyset)$) together with the composition
has a structure of a commutative monoid. We will see that any
irreducible element of this monoid (i.e. a system of curves that
is not a composition of other systems of curves) is a system of
curves encircled by another curve. Note that to encircle a system
of curves by another curve is, in the $(\pm 2,k)$ words language,
the same as adding a $(-2,0)$ at the beginning and a $(2,0)$ at
the end of the corresponding word. Thus if a word
$(c_1,0)(c_2,0)...(c_n,0)$ of symbols $(-2,0)$ and
$(2,0)$\footnote{By the previous proposition, any system of
non-singular planar curves can be represented in this way.}
satisfing the condition C (and thus $(c_1,0)=(-2,0)$ and
$(c_n,0)=(2,0)$) doesn't come from a system of curves encircled by
another curve (this means that the word $(c_2,0)...(c_{n-1},0)$
doesn't satisfy the condition C) then there exists $2<k<n$ such
that $\sum_{j<k}c_j =0$, which implies that $(c_{k-1},0)=(2,0)$,
$(c_k,0)=(-2,0)$ and the words $(c_1,0)...(c_{k-1},0)$ and
$(c_{k},0)...(c_{n},0)$ satisfy the condition C. Therefore the
morphism in $\hom(\emptyset,\emptyset)$ (corresponding to the word
$(c_1,0)...(c_n,0)$) is the composition of two morphisms in
$\hom(\emptyset,\emptyset)$ (corresponding to the words
$(c_1,0)...(c_{k-1},0)$ and $(c_{k},0)...(c_{n},0)$).

Note that we have just proved that an irreducible morphism is
another morphism encircled by an exterior curve, but we haven't
proved yet that a morphism encircled by an exterior curve is an
irreducible morphism. This is because the last proposition says
nothing about when two words of symbols $(-2,0)$ and $(2,0)$
represent equivalent words (by the relations 1, 2, 3 and 4).

However we do know that, using the composition of morphisms and
the operation of encircling, we can generate all morphisms in
$\hom(\emptyset,\emptyset)$ from the identity.

With this in mind, we are going to see what happens in the
representation when we compose two morphisms (see corollary
\ref{additivity}) or encircle one by a circle (see corollary
\ref{varphi}).


Consider, for a fixed value $\mathbf{m}\in \mathbb{M}$ and for
each natural number $n$, the following morphism:
$$\begin{array}{cccc} \Psi_{n,\mathbf{m}}: & \mathcal{O}_n &
\longrightarrow & \mathcal{O}_{n}\\ & (R,\overrightarrow{v}) &
\longmapsto &(R,R*(e_1 *\mathbf{m}\oplus\overrightarrow{v}))
\end{array}$$

The purpose of this morphism is to add the value $\mathbf{m}$ to
the region associated to the first interval.


\begin{lemma}
For any $2\leq k\leq n+1$, $\hat{T}_{n,k}\circ\Psi_{n,\mathbf{m}}
=\Psi_{n+2,\mathbf{m}}\circ\hat{T}_{n,k}$ and
$\check{T}_{n,k}\circ\Psi_{n+2,\mathbf{m}}
=\Psi_{n,\mathbf{m}}\circ\check{T}_{n,k}$
\end{lemma}
\TeXButton{Proof}{\proof}
\begin{description}
\item[] First equation: $\hat{T}_{n,k}\circ\Psi_{n,\mathbf{m}}
=\Psi_{n+2,\mathbf{m}}\circ\hat{T}_{n,k}$.

For an arbitrary $(R,\overrightarrow{v})\in \mathcal{O}_n$, let
$(R_1,\overrightarrow{v}_1)=\Psi_{n,\mathbf{m}}(R,\overrightarrow{v})$,
$(R_2,\overrightarrow{v}_2)=\hat{T}_{n,k}(R_1,\overrightarrow{v}_1)$,
$(R'_1,\overrightarrow{v}'_1)=\hat{T}_{n,k}(R,\overrightarrow{v})$
and
$(R'_2,\overrightarrow{v}'_2)=\Psi_{n+2,\mathbf{m}}(R'_1,\overrightarrow{v}'_1)$.
We want to check
$(R_2,\overrightarrow{v}_2)=(R'_2,\overrightarrow{v}'_2)$.

Since $\Psi$ does not change the matrices and the changes of the
matrices under the morphisms $\hat{T}$ do not depend on the choice
of the array of values, it is obvious that $R_2=R'_2$.

$$\begin{array}{lllll} \overrightarrow{v}_2 &=&
B_{n,k}*\overrightarrow{v}_1 & &\\ &=&
B_{n,k}*[R*(e_1*\mathbf{m}\oplus\overrightarrow{v}_1)] & &\\ &=&
B_{n,k}R*(e_1*\mathbf{m}\oplus\overrightarrow{v}_1)& &\\ &=&
B_{n,k}RB_{n,k}^tB_{n,k}*(e_1*\mathbf{m}\oplus\overrightarrow{v}_1)&
&
\\ &=&
[B_{n,k}RB_{n,k}^tB_{n,k}*(e_1*\mathbf{m}\oplus\overrightarrow{v}_1)]
\vee
[D_{n+2,k}B_{n,k}*(e_1*\mathbf{m}\oplus\overrightarrow{v}_1)]&
&
\\ &=&
(B_{n,k}RB_{n,k}^t+D_{n+2,k})*[B_{n,k}*(e_1*\mathbf{m}\oplus\overrightarrow{v}_1)]&
& \\ &=&
R'_1*[(B_{n,k}e_1*\mathbf{m})\oplus(B_{n,k}*\overrightarrow{v}_1)]&
& \\ &=&
R'_1*[(B_{n,k}e_1*\mathbf{m})\oplus\overrightarrow{v}'_1]& &
\end{array}$$

On the other hand, $$\begin{array}{ccc}
  \overrightarrow{v}'_2 & = & R'_1*(e_1*\mathbf{m}\oplus\overrightarrow{v}'_1)
\end{array}$$

If $k>2$ then $B_{n,k}e_1=e_1$ and therefore
$\overrightarrow{v}_2=\overrightarrow{v}'_2$.

If $k=2$ then $$\begin{array}{ccl}
  \overrightarrow{v}_2 & = & R'_1*[(B_{n,k}e_1*\mathbf{m})\oplus\overrightarrow{v}'_1] \\
   & = & R'_1*[(e_1*\mathbf{m}\vee e_3*\mathbf{m})\oplus\overrightarrow{v}'_1] \\
   & = & R'_1*(e_1*\mathbf{m}\oplus\overrightarrow{v}'_1)\vee
   R'_1*(e_3*\mathbf{m}\oplus\overrightarrow{v}'_1)\\
   & = & R'_1*(e_1*\mathbf{m}\oplus\overrightarrow{v}'_1)\\
   & = & \overrightarrow{v}'_2
\end{array}$$
because $R'_1*(e_3*\mathbf{m}\oplus\overrightarrow{v}'_1)=
R'_1*(e_1*\mathbf{m}\oplus\overrightarrow{v}'_1)$ as we will prove
next:

$$\begin{array}{ccccccc}
  R'_1 e_1 & = & (B_{n,2}RB_{n,2}^t+D_{n+2,2})e_1 & = &
  B_{n,2}RB_{n,2}^t e_1\vee D_{n+2,2}e_1 & = & B_{n,2}R e_1 \\
  R'_1 e_3 & = & (B_{n,2}RB_{n,2}^t+D_{n+2,2})e_3 & = &
  B_{n,2}RB_{n,2}^t e_3\vee D_{n+2,2}e_3 & = & B_{n,2}R e_1
\end{array}$$
thus $R'_1 e_1*\mathbf{m}=R'_1 e_3*\mathbf{m}$, and using the next
lemma we conclude
$R'_1*(e_3*\mathbf{m}\oplus\overrightarrow{v}'_1)=
R'_1*(e_1*\mathbf{m}\oplus\overrightarrow{v}'_1)$.

\begin{lemma}
If $R=R^t$ and $R*\overrightarrow{v}=\overrightarrow{v}$ then
$R*(\overrightarrow{u}\oplus\overrightarrow{v})
=R*\overrightarrow{u}\oplus \overrightarrow{v}$.
\end{lemma}
\TeXButton{Proof}{\proof} Let
$\overrightarrow{x}:=R*(\overrightarrow{u}\oplus\overrightarrow{v})$
and $\overrightarrow{y}:=R*\overrightarrow{u}\oplus
\overrightarrow{v}$, then $$x_i=\bigvee_j r_{i,j}*(u_j\oplus v_j)
\quad \mbox{and} \quad y_i=(\bigvee_j r_{i,j}*u_j)\oplus v_i$$

$$
\begin{array}{ccl}
  x_i & = & \bigvee_j r_{i,j}*(u_j\oplus v_j) \\
   & = & \bigvee_{j: r_{i,j}=1}(u_j\oplus v_j) \\
   & = & \bigvee_{j: r_{i,j}=1}(u_j\oplus v_i) \\
   & = & (\bigvee_{j: r_{i,j}=1}u_j)\oplus v_i \\
   & = & (\bigvee_j r_{i,j}*u_j)\oplus v_i \\
   & = & y_i
\end{array}
$$

Note that, since $R*\overrightarrow{v}=\overrightarrow{v}$, we
have that $v_j=v_i$ whenever $r_{i,j}=1$.
 \TeXButton{End Proof}{\endproof}

\item[] Second equation $\check{T}_{n,k}\circ\Psi_{n+2,\mathbf{m}}
=\Psi_{n,\mathbf{m}}\circ\check{T}_{n,k}$.

For an arbitrary $(R,\overrightarrow{v})\in \mathcal{O}_n$, let
$(R_1,\overrightarrow{a})=\Psi_{n+2,\mathbf{m}}(R,\overrightarrow{v})$,
$(R_2,\overrightarrow{b})=\check{T}_{n,k}(R_1,\overrightarrow{a})$,
$(R'_1,\overrightarrow{a}')=\check{T}_{n,k}(R,\overrightarrow{v})$
and
$(R'_2,\overrightarrow{b}')=\Psi_{n,\mathbf{m}}(R'_1,\overrightarrow{a}')$.
We want to check
$(R_2,\overrightarrow{b})=(R'_2,\overrightarrow{b}')$.

Since $\Psi$ does not change the matrices and the changes of the
matrices under the morphisms $\check{T}$ do not depend on the
choice of the array of values, it is obvious that $R_2=R'_2$.

$$\overrightarrow{b}= R_2*[(B_k^t*\overrightarrow{a})\oplus
e_{k-1}*\tilde{x}_k]$$ with
$$\overrightarrow{a}=R*(e_1*\mathbf{m}\oplus\overrightarrow{v})
=Re_1*\mathbf{m}\oplus\overrightarrow{v}$$ and
$$\tilde{x}_k=\left\{\begin{array}{ccc}
  a_{k-1}\wedge a_{k+1} & \mbox{if} & r_{k-1,k+1}=0 \\
  \varphi(a_k) & \mbox{if} & r_{k-1,k+1}=1
\end{array}\right.$$

On the other hand
$$\overrightarrow{b}'=R'_2*(e_1*\mathbf{m}\oplus\overrightarrow{a}')
=R'_2e_1*\mathbf{m}\oplus\overrightarrow{a}'$$

We want to check $\overrightarrow{b}=\overrightarrow{b}'$.

\begin{description}
\item[A)] $b_{k-1}=b'_{k-1}$.

By lemma \ref{11v}, we have $$b_{k-1}=(a_{k-1}\vee a_{k+1})\oplus
\tilde{x}_k =\left\{\begin{array}{ccc}
   a_{k-1}\oplus a_{k+1} & \mbox{if} & r_{k-1,k+1}=0 \\
  a_{k-1}\oplus\varphi(a_k) & \mbox{if} & r_{k-1,k+1}=1
\end{array} \right.$$

Bearing in mind that

$$a_{k-1}=e_{k-1}^t*(Re_1*\mathbf{m}\oplus\overrightarrow{v})
=r_{k-1,1}*\mathbf{m}\oplus v_{k-1}$$
$$a_{k+1}=e_{k+1}^t*(Re_1*\mathbf{m}\oplus\overrightarrow{v})
=r_{k+1,1}*\mathbf{m}\oplus v_{k+1}$$
$$a_{k}=e_{k}^t*(Re_1*\mathbf{m}\oplus\overrightarrow{v})
=r_{k,1}*\mathbf{m}\oplus v_{k}$$ we have that

$$b_{k-1}=\left\{\begin{array}{ccc}
   (r_{k-1,1}*\mathbf{m})\oplus v_{k-1}\oplus
   (r_{k+1,1}*\mathbf{m})\oplus v_{k+1}
   & \mbox{if} & r_{k-1,k+1}=0 \\
  (r_{k-1,1}*\mathbf{m})\oplus v_{k-1}\oplus
  \varphi(r_{k,1}*\mathbf{m}\oplus v_{k}) & \mbox{if} & r_{k-1,k+1}=1
\end{array} \right.$$

Now we note that $r_{k-1,k+1}=0$ implies that $r_{k-1,1}=0$ or
$r_{k+1,1}=0$ and therefore
$(r_{k-1,1}*\mathbf{m})\oplus(r_{k+1,1}*\mathbf{m})=
(r_{k-1,1}+r_{k+1,1})*\mathbf{m}$. On the other hand
$r_{k-1,k+1}=1$ implies that $r_{k,1}=0$ and
$r_{k-1,1}=r_{k+1,1}=r_{k-1,1}+r_{k+1,1}$.

Thus

$$b_{k-1}=\left\{\begin{array}{ccc}
   (r_{k-1,1}+r_{k+1,1})*\mathbf{m}\oplus v_{k-1}\oplus v_{k+1}
   & \mbox{if} & r_{k-1,k+1}=0 \\
  (r_{k-1,1}+r_{k+1,1})*\mathbf{m}\oplus v_{k-1}\oplus
  \varphi(v_{k}) & \mbox{if} & r_{k-1,k+1}=1
\end{array} \right.$$

That means $$b_{k-1}=(r_{k-1,1}+r_{k+1,1})*\mathbf{m}\oplus
a'_{k-1}$$ by lemma \ref{11v}. Now

$$\overrightarrow{b}'=R'_1*(e_1*\mathbf{m}\oplus\overrightarrow{a}')
=R'_1e_1*\mathbf{m}\oplus\overrightarrow{a}'$$ i.e.

$$\begin{array}{ccl}
  b'_{k-1} & = & e_{k-1}^t*(R'_1e_1*\mathbf{m}\oplus\overrightarrow{a}') \\
   & = & e_{k-1}^tR'_1e_1*\mathbf{m}\oplus a'_{k-1} \\
   & = & r'_{k-1,1}*\mathbf{m}\oplus a'_{k-1}
\end{array}$$ and using lemma \ref{12R} we have

$$\begin{array}{ccl}
  b'_{k-1} & = & (r_{k-1,1}+r_{k+1,1})*\mathbf{m}\oplus a'_{k-1} \\
   & = & b_{k-1}
\end{array}$$

\item[B)] $b_{i}=b'_{i}$ with $r'_{i,k-1}=0$.

By lemma \ref{11v}, $b_i=e_i^tB_k^t*\overrightarrow{a}$.

$$\begin{array}{ccl}
  b_i & = & e_i^tB_k^t*(Re_1*\mathbf{m}\oplus\overrightarrow{v}) \\
   & = & e_i^tB_k^tRe_1*\mathbf{m}\oplus e_i^tB_k^t\overrightarrow{v} \\
   & = & r_{\check{k}(i),1}*\mathbf{m}\oplus e_i^tB_k^t\overrightarrow{v}
\end{array}$$

On the other hand,

$$\begin{array}{ccl}
  b'_i & = & e_i^t(R'_1e_1*\mathbf{m}\oplus\overrightarrow{a}') \\
   & = & e_i^tR'_1e_1*\mathbf{m}\oplus e_1^t\overrightarrow{a}') \\
   & = & r'_{i,1}*\mathbf{m}\oplus e_i^tB_k^t\overrightarrow{v}
\end{array}$$

Thus we have $b_{i}=b'_{i}$ if $r'_{i,1}=r_{\check{k}(i),1}$.

If $k-1\not= 1$ then, by lemma \ref{12R},
$r'_{i,1}=r_{\check{k}(i),\check{k}(1)}=r_{\check{k}(i),1}$.

If $k-1=1$ then, by lemma \ref{12R}, $r'_{i,1}=r'_{i,k-1}\geq
r_{\check{k}(i),k-1}=r_{\check{k}(i),1}$ and since
$r'_{i,1}=r'_{i,k-1}=0$ we have $r'_{i,1}=r_{\check{k}(i),1}$.

\item[C)] $b_{i}=b'_{i}$ with $r'_{i,k-1}=1$.

Since, by lemma \ref{11v}, $b_{i}=b_{k-1}$ and $b'_{i}=b'_{k-1}$
this shows that $b_{i}=b'_{i}$.

\end{description}

\end{description}
 \TeXButton{End Proof}{\endproof}

\begin{corollary}
\begin{description}
\item[1.] For any $T\in \hom(\mathcal{O}_{m},\mathcal{O}_{n})$,
$T\circ\Psi_{m,\mathbf{m}}=\Psi_{n,\mathbf{m}}\circ T;$
\item[2.] For any $T\in \hom(\mathcal{O}_{1},\mathcal{O}_{1})$
(corresponding to a system of non-singular planar curves)
$T([1],\mathbf{\emptyset})=([1],\mathbf{x}) \quad \Rightarrow
\quad T([1],\mathbf{m})=([1],\mathbf{x}\oplus\mathbf{m})$;
\item[3.] For any $T_1 , T_2\in \hom(\mathcal{O}_{1},\mathcal{O}_{1})$
$T_2\circ
T_1([1],\mathbf{\emptyset})=([1],\mathbf{m_2}\oplus\mathbf{m_1})$
where $([1],\mathbf{m_i})= T_i([1],\mathbf{\emptyset})$, $i=1,2$.
\end{description}
\label{additivity}
\end{corollary}
\TeXButton{Proof}{\proof}
\begin{description}
\item[1.] It is obvious that if
$\{\Psi_{n,\mathbf{m}}\}_{n\in\mathbb{N}}$ commute with the
generators $\hat{T}_{n,k}$ and $\check{T}_{n,k}$ then they commute
with any morphism $T$.

\item[2.] $T([1],\mathbf{m})=T\circ\Psi_{1,\mathbf{m}}([1],\mathbf{\emptyset})
=\Psi_{1,\mathbf{m}}\circ T([1],\mathbf{\emptyset})=
\Psi_{1,\mathbf{m}}([1],\mathbf{x})=
([1],\mathbf{x}\oplus\mathbf{m})$.

\item[3.] $T_2\circ T_1([1],\mathbf{\emptyset})=
T_2([1],\mathbf{m_1})=([1],\mathbf{m_2}\oplus\mathbf{m_1})$.
\end{description}
 \TeXButton{End Proof}{\endproof}

This result shows that all information about a morphism $T\in
\hom(\mathcal{O}_{1},\mathcal{O}_{1})$ is contained in a single
value $\mathbf{v}\in \mathbb{M}$ (which is obtained evaluating the
morphism $T$ on $([1],\emptyset)$). Thus, we can associate a value
in $\mathbb{M}$ to each system of non-singular planar curves.
Moreover this association is a monoid homomorphism.

Next we will look at what value is obtained when a collection of
planar curves is encircled by another curve.


Consider, for each $n$, the following subset of $\mathcal{O}_n$
$$\mathcal{U}_n =\{(R,\overrightarrow{v})\in \mathcal{O}_n :
e_1^tRe_n=1\}$$

By the property T1 we have that $\mathcal{U}_n =\emptyset$ for any
even number $n$. It is clear that $\mathcal{U}_1=\mathcal{O}_1$
and therefore it easy to see that $\mathcal{U}_n \not=\emptyset$
for any odd number $n$, using the following result.

\begin{proposition}
For any $n$ and $k$, if $(R,\overrightarrow{v})\in\mathcal{U}_n$
then $\hat{T}_k(R,\overrightarrow{v})\in\mathcal{U}_{n+2}$ and
$\check{T}_k(R,\overrightarrow{v})\in\mathcal{U}_{n-2}$.
\end{proposition}
\TeXButton{Proof}{\proof}
Let $(R',\overrightarrow{v}')=\hat{T}_k(R,\overrightarrow{v})$. We
want to see that if $e_1^tRe_n=1$ then $e_1^tR'e_{n+2}=1$.
$$\begin{array}{ccl}
  e_1^tR'e_{n+2} & = & e_1^t(B_kRB_k^t+D_k)e_{n+2} \\
   & = & e_1^tB_kRB_k^te_{n+2} \\
   & = & e_1^tRe_n=1
\end{array}$$

Let $(R',\overrightarrow{v}')=\check{T}_k(R,\overrightarrow{v})$.
We want to see that if $e_1^tRe_n=1$ then $e_1^tR'e_{n-2}=1$.
$$\begin{array}{ccl}
  e_1^tR'e_{n-2} & = & e_1^t(B_k^tRB_k)^2e_{n-2} \\
   & \geq & e_1^tB_k^tRB_ke_{n-2} \\
   & \geq & e_1^tRe_n=1
\end{array}$$
 \TeXButton{End Proof}{\endproof}

Now we consider, for each odd natural number $n$, the following
morphism:

$$\begin{array}{cccc}
  \varepsilon_n : & \mathcal{U}_n & \longrightarrow & \mathcal{U}_{n+2} \\
   & (R,\overrightarrow{v}) & \longmapsto & (\tilde{R},\overrightarrow{\tilde{v}})
\end{array}$$
with $$\tilde{R}=E_n RE_n^t+F_{n+2}$$ and
$$\overrightarrow{\tilde{v}}=E_n*\overrightarrow{v}$$ where
$E_n=[e_{i,j}]$ is an $(n+2)\times n$ matrix defined by $e_{i,j}=1
\Leftrightarrow i-1=j$ and $F_{n+2}=[f_{i,j}]$ is an $(n+2)\times
(n+2)$ matrix defined by $f_{i,j}=1 \Leftrightarrow
i,j\in\{1,n+2\}$.

\begin{proposition}
For any odd natural number $n$, $\varepsilon_n$ is well defined.
\end{proposition}
\TeXButton{Proof}{\proof} We need to see that if
$(R,\overrightarrow{v})\in\mathcal{U}_{n}$ then $\varepsilon_n
(R,\overrightarrow{v})\in \mathcal{U}_{n+2}$.

Let $(\tilde{R},\overrightarrow{\tilde{v}})=\varepsilon_n
(R,\overrightarrow{v})$. We are going to check that
$R*\overrightarrow{v}=\overrightarrow{v}$ implies
$\tilde{R}*\overrightarrow{\tilde{v}}=\overrightarrow{\tilde{v}}$,
$\tilde{R}$ satisfies the properties E1, E2, E3, T1, T2 and T3 and
the $(1,n+2)$ entry of $\tilde{R}$ is $1$.

For this purpose, we are going to use the identities:
$E_n^tE_n=I$, $E_n^tF_{n+2}=O$ and $F_{n+2}E_n=O$, which we will
leave to the reader to check.

$$\begin{array}{ccl}
  R*\overrightarrow{v}=\overrightarrow{v} \Rightarrow \tilde{R}*\overrightarrow{\tilde{v}}
   & = & (E_nRE_n^t+F_{n+2})*(E_n*\overrightarrow{v}) \\
   & = & (E_nRE_n^tE_n)*\overrightarrow{v} \vee (F_{n+2}E_n)*\overrightarrow{v} \\
   & = & E_nR*\overrightarrow{v} \\
   & = & E_n*\overrightarrow{v} \\
   & = & \overrightarrow{\tilde{v}}
\end{array}$$

The properties E1, E2 and E3 of $\tilde{R}$ are very easy to
verify, so we leave them as an exercise.

For the properties T1, T2 and T3, we observe that the $(i,j)$
entry of $\tilde{R}$ is the $(i-1,j-1)$ entry of $R$ if $2\leq
i,j\leq n+1$, $1$ if $i,j\in\{1,n+2\}$ and $0$ in all other cases,
i.e.: $$\tilde{r}_{i,j}=e_i^t(E_nRE_n^t+F_{n+2})e_j
=\left\{\begin{array}{cl}
  r_{i-1,j-1} & \mbox{if} \quad 2\leq i,j\leq n+1 \\
  1 & \mbox{if} \quad i,j\in\{1,n+2\} \\
  0 & \mbox{otherwise}
\end{array}\right.$$

In particular, we have $\tilde{r}_{1,n+2}=1$.

\begin{description}

\item[] Property T1: $\tilde{r}_{i,j}=1 \Rightarrow j-i\in
2\mathbb{Z}$.

We have $\tilde{r}_{i,j}=1 \Rightarrow i,j\in \{1,n+2\}$ or
$r_{i-1,j-1}=1$. Since we have taken $n$ to be odd and $R$
satisfies T1 we have $j-i \in 2\mathbb{Z}$.

\item[] Property T2: $\forall_{\alpha\leq \beta\leq \gamma\leq \delta}
\tilde{r}_{\alpha,\gamma}=\tilde{r}_{\beta,\delta}=1 \Rightarrow
\tilde{r}_{\alpha,\beta}=\tilde{r}_{\beta,\gamma}=\tilde{r}_{\gamma,\delta}=1$.

We only need consider the case $\alpha<\beta<\gamma<\delta$, and
then $\tilde{r}_{\alpha,\gamma}=\tilde{r}_{\beta,\delta}=1
\Rightarrow 2\leq \alpha<\beta<\gamma<\delta\leq n+1$. Thus
$\tilde{r}_{\alpha,\gamma}=r_{\alpha-1,\gamma-1}$ and
$\tilde{r}_{\beta,\delta}=r_{\beta-1,\delta-1}$, hence, since $R$
satisfies T2, we have
$r_{\alpha-1,\beta-1}=r_{\beta-1,\gamma-1}=r_{\gamma-1,\delta-1}=1$,
that is
$\tilde{r}_{\alpha,\beta}=\tilde{r}_{\beta,\gamma}=\tilde{r}_{\gamma,\delta}=1$.

\item[] Property T3: $\forall_{\alpha<\beta}
\tilde{r}_{\alpha,\beta}=1 \Rightarrow
\tilde{r}_{\alpha+1,\beta-1}=1$ or
$\exists_{\alpha<\gamma<\beta}:\tilde{r}_{\alpha,\gamma}=1$.

If $\alpha=1$ then $\tilde{r}_{\alpha,\beta}=1 \Rightarrow
\beta=n+2$. Thus
$\tilde{r}_{\alpha+1,\beta-1}=\tilde{r}_{2,n+1}=r_{1,n}=1$ since
$R\in\mathcal{U}_n$\footnote{This is an abuse of notation, since
we should write $(R,\overrightarrow{v})\in\mathcal{U}_n$.}.

If $\alpha>1$ then $\tilde{r}_{\alpha,\beta}=1 \Rightarrow
\beta<n+2$. Thus $$\begin{array}{ccccc}
  \tilde{r}_{\alpha,\beta}=r_{\alpha-1,\beta-1} & \Rightarrow
  & r_{\alpha,\beta-2}=1 & \mbox{or} &
  \exists_{\alpha-1<\gamma-1<\beta-1}:r_{\alpha-1,\gamma-1}=1 \\
   & \Rightarrow & \tilde{r}_{\alpha+1,\beta-1}=1 & \mbox{or} &
   \exists_{\alpha<\gamma<\beta}:\tilde{r}_{\alpha,\gamma}=1
\end{array}$$

\end{description}

 \TeXButton{End Proof}{\endproof}

\begin{lemma}
For any $2\leq k\leq n+1$, we have the identities
$\hat{T}_{n+2,k+1}\circ\varepsilon_{n}
=\varepsilon_{n+2}\circ\hat{T}_{n,k}$ and
$\check{T}_{n,k+1}\circ\varepsilon_{n}
=\varepsilon_{n-2}\circ\check{T}_{n-2,k}$
\end{lemma}
\TeXButton{Proof}{\proof}


\begin{description}

\item[] First identity: $\hat{T}_{n+2,k+1}\circ\varepsilon_{n}
=\varepsilon_{n+2}\circ\hat{T}_{n,k}$.

For an arbitrary $(R,\overrightarrow{v})\in \mathcal{U}_n$, let
$(R_1,\overrightarrow{v}_1)=\varepsilon_n(R,\overrightarrow{v})$,
$(R_2,\overrightarrow{v}_2)=\hat{T}_{n+2,k+1}(R_1,\overrightarrow{v}_1)$,
$(R'_1,\overrightarrow{v}'_1)=\hat{T}_{n,k}(R,\overrightarrow{v})$
and
$(R'_2,\overrightarrow{v}'_2)=\varepsilon_{n+2}(R'_1,\overrightarrow{v}'_1)$.
We want to check
$(R_2,\overrightarrow{v}_2)=(R'_2,\overrightarrow{v}'_2)$.

$$\begin{array}{ccl}
  R_2 & = & B_{n+2,k+1}R_1B_{n+2,k+1}^t+D_{n+4,k+1} \\
   & = & B_{n+2,k+1}(E_nRE_n^t+F_{n+2})B_{n+2,k+1}^t+D_{n+4,k+1} \\
   & = & B_{n+2,k+1}E_nRE_n^tB_{n+2,k+1}^t+B_{n+2,k+1}F_{n+2}B_{n+2,k+1}^t+D_{n+4,k+1} \\
   & = & E_{n+2}B_{n,k}RB_{n,k}^tE_{n+2}^t+F_{n+4}+E_{n+2}D_{n+2,k}E_{n+2}^t \\
   & = & E_{n+2}(B_{n,k}RB_{n,k}^t+D_{n+2,k})E_{n+2}^t+F_{n+4} \\
   & = & E_{n+2}R'_1E_{n+2}^t+F_{n+4} \\
   & = & R'_2
\end{array}$$

$$\begin{array}{ccl}
  \overrightarrow{v}_2 & = & B_{n+2,k+1}*\overrightarrow{v}_1 \\
   & = & B_{n+2,k+1}*(E_n*\overrightarrow{v}) \\
   & = & (B_{n+2,k+1}E_n)*\overrightarrow{v} \\
   & = & (E_{n+2}B_{n,k})*\overrightarrow{v} \\
   & = & E_{n+2}*(B_{n,k}*\overrightarrow{v}) \\
   & = & E_{n+2}*\overrightarrow{v}'_1 \\
   & = & \overrightarrow{v}'_2
\end{array}$$


\item[] Second identity $\check{T}_{n,k+1}\circ\varepsilon_{n}
=\varepsilon_{n-2}\circ\check{T}_{n-2,k}$.

For an arbitrary $(R,\overrightarrow{v})\in \mathcal{U}_n$, let
$(R_1,\overrightarrow{a})=\varepsilon_n(R,\overrightarrow{v})$,
$(R_2,\overrightarrow{b})=\check{T}_{n,k+1}(R_1,\overrightarrow{a})$,
$(R'_1,\overrightarrow{a}')=\check{T}_{n-2,k}(R,\overrightarrow{v})$
and
$(R'_2,\overrightarrow{b}')=\varepsilon_{n-2}(R'_1,\overrightarrow{a}')$.
We want to check
$(R_2,\overrightarrow{b})=(R'_2,\overrightarrow{b}')$.

$$\begin{array}{ccl}
  R_2 & = & [B_{n,k+1}^tR_1B_{n,k+1}]^2 \\
   & = & [B_{n,k+1}^t(E_nRE_n^t+F_{n+2})B_{n,k+1}]^2 \\
   & = & [B_{n,k+1}^tE_nRE_n^tB_{n,k+1}+B_{n,k+1}^tF_{n+2}B_{n,k+1}]^2 \\
   & = & (E_{n-2}B_{n-2,k}^tRB_{n-2,k}E_{n-2}^t+F_n)^2 \\
   & = & E_{n-2}B_{n-2,k}^tRB_{n-2,k}E_{n-2}^tE_{n-2}B_{n-2,k}^tRB_{n-2,k}E_{n-2}^t \\
   &  & +E_{n-2}B_{n-2,k}^tRB_{n-2,k}E_{n-2}^tF_n +F_nE_{n-2}B_{n-2,k}^tRB_{n-2,k}E_{n-2}^t+F_n^2 \\
   & = & E_{n-2}R'_1E_{n-2}^t+F_n \\
   & = & R'_2
\end{array}$$

$$\begin{array}{ccl}
  \overrightarrow{b} & = & R_2*[(B_{n,k+1}^t*\overrightarrow{a})\oplus(e_k*\tilde{x}_{k+1})] \\
   & = & (E_{n-2}R'_1E_{n-2}^t+F_n)*
   \{[B_{n,k+1}^t*(E_n*\overrightarrow{v})]\oplus(e_k*\tilde{x}_{k+1})\} \\
   & = &
   E_{n-2}R'_1E_{n-2}^t*[(B_{n,k+1}^tE_n*\overrightarrow{v})\oplus(e_k*\tilde{x}_{k+1})] \\
   &  & \vee F_n*[(B_{n,k+1}^tE_n*\overrightarrow{v})\oplus(e_k*\tilde{x}_{k+1})] \\
   & = & E_{n-2}R'_1*\{E_{n-2}^t*(E_{n-2}B_{n-2,k}^t*\overrightarrow{v})
   \oplus(E_{n-2}^te_k*\tilde{x}_{k+1})\} \\
   &  & \vee (F_nE_{n-2}B_{n-2,k}^t*\overrightarrow{v})\oplus(F_ne_k*\tilde{x}_{k+1})] \\
   & = & E_{n-2}R'_1*[(E_{n-2}^tE_{n-2}B_{n-2,k}^t*\overrightarrow{v})\oplus(e_{k-1}*\tilde{x}_{k+1})] \\
   & = & E_{n-2}R'_1*[(B_{n-2,k}^t*\overrightarrow{v})\oplus(e_{k-1}*x_k)] \\
   & = & E_{n-2}*\overrightarrow{a}' \\
   & = & \overrightarrow{b}'
\end{array}$$

Here $\tilde{x}_{k+1}=\neg\tilde{r}_{k,k+2}*(a_k\wedge
a_{k+2})\oplus \tilde{r}_{k,k+2}*\varphi(a_{k+1})$ where
$\tilde{r}_{k,k+2}=e_k^tR_1e_{k+2}$.

Since $1<k, k+2<n+2$ we have that $\tilde{r}_{k,k+2}=r_{k-1,k+1}$
and $a_k=v_{k-1}$, $a_{k+1}=v_{k}$ and $a_{k+2}=v_{k+1}$. Thus
$\tilde{x}_{k+1}=x_k$.

\end{description}


 \TeXButton{End Proof}{\endproof}

\begin{corollary}
Let $F:\mathcal{PI}_\mathbb{M}\longrightarrow
\mathcal{PI}_\mathbb{M}$ be the functor which sends
$\hat{T}_{n,k}$ and $\check{T}_{n,k}$ to $\hat{T}_{n+2,k+1}$ and
$\check{T}_{n+2,k+1}$ (respectively). Then

\begin{description}
\item[1.] For any $T\in\hom(\mathcal{U}_n,\mathcal{U}_m)$,
$\varepsilon_m\circ T=F(T)\circ\varepsilon_n$.

\item[2.] For any $T\in\hom(\mathcal{O}_1,\mathcal{O}_1)$
(note that $\mathcal{U}_1=\mathcal{O}_1$),
$T([1],\mathbf{\emptyset})=([1],\mathbf{x})\Rightarrow
\check{T}_{1,2}\circ
F(T)\circ\hat{T}_{1,2}([1],\mathbf{\emptyset})=([1],\varphi(\mathbf{x}))$.
\end{description}
\label{varphi}
\end{corollary}
\TeXButton{Proof}{\proof}

\begin{description}
\item[1.] This follows immediately from the previous lemma and
the definition of the functor $F$.

\item[2.] $$\begin{array}{ccl}
  \check{T}_{1,2}\circ F(T)\circ\hat{T}_{1,2}([1],\mathbf{\emptyset})
  & = & \check{T}_{1,2}\circ
F(T)\left(\left[
\begin{array}{ccc}
  1 & 0 & 1 \\
  0 & 1 & 0 \\
  1 & 0 & 1
\end{array}
\right],\left(
\begin{array}{c}
  \emptyset \\
  \emptyset \\
  \emptyset
\end{array}
\right)\right) \\
   & = & \check{T}_{1,2}\circ F(T)\circ\varepsilon_1([1],\mathbf{\emptyset}) \\
   & = & \check{T}_{1,2}\circ\varepsilon_1\circ T([1],\mathbf{\emptyset}) \\
   & = & \check{T}_{1,2}\circ\varepsilon_1 ([1],\mathbf{x}) \\
   & = & \check{T}_{1,2}\left(\left[\begin{array}{ccc}
  1 & 0 & 1 \\
  0 & 1 & 0 \\
  1 & 0 & 1
\end{array}\right],\left(\begin{array}{c}
  \emptyset \\
  \mathbf{x} \\
  \emptyset
\end{array}\right)\right) \\
   & = & ([1],\varphi(\mathbf{x}))
\end{array}
$$

\end{description}

 \TeXButton{End Proof}{\endproof}

Note that $\check{T}_{1,2}\circ F(T)\circ\hat{T}_{1,2}$ is the
representation of $\check{t}_{1,2}\circ
F(t)\circ\hat{t}_{1,2}$\footnote{Here the functor $F$ is defined
in the same way for ${\bf{PT}}$ as it was for
$\mathbf{PI}_\mathbb{M}$.} which corresponds to the encirclement
of the morphism $t$.


In summary we see that the representation of a system of
non-singular planar curves is a morphism in
$\hom(\mathcal{O}_1,\mathcal{O}_1)$, but this is determined by a
single value. This means that the representation gives an
application of the monoid of systems of non-singular planar curves
$\hom(\emptyset,\emptyset)$ to the monoid of the representation
$\mathbb{M}$. To a morphism $t\in\hom(\emptyset,\emptyset)$ we
associate the value $v(t)\in\mathbb{M}$ such that
$([1],v(t))=T([1],\emptyset)$ where
$T\in\hom(\mathcal{O}_1,\mathcal{O}_1)$ is the representation of
$t$. This application is a monoid morphism ($v(t_1\circ
t_2)=v(t_1)\oplus v(t_2)$) as we have seen (see corollary
\ref{additivity}), and furthermore, the function $\varphi$ is what
corresponds in $\mathbb{M}$ to the operation of encircling a
system of curves by another curve (see corollary \ref{varphi}),
i.e. $v(\langle t\rangle )=\varphi(v(t))$ where $\langle t\rangle$
denotes the encirclement of the morphism $t$.


We are going to study this application in the following two
particular cases.

{\bf First case:} the monoid is the natural numbers with the usual
multiplication as the operation of the monoid and with the
division order giving the lattice structure; the function
$\varphi$ is the function that sends a number $n$ to the $n^{th}$
prime number.

We are going to see that, for systems of non-singular planar
curves, the application is a monoid isomorphism (and thus is a
complete invariant for such systems). This means that two systems
of curves with the same value are equivalent (isotopic). Let us
prove this by induction on the value $v(s)=n$ of the system $s$.
Note that it only remains to prove that the application is
bijective.

First we observe that if a morphism
$s\in\hom(\emptyset,\emptyset)$ is irreducible (and thus is
encircled) then $v(s)$ is a prime number.

If $v(s)=1$ then $s$ is the empty system of curves (the identity
in $\hom(\emptyset,\emptyset)$) because a non-empty system of
curves is a non-empty composition of irreducible morphisms
therefore it has a non-empty product of prime numbers as value $v$
(i.e. $v(s)>1$). This implies that if a morphism has a prime
number as value then it is irreducible. Thus an encircled morphism
is irreducible.

Now let $v(s_1)=v(s_2)=n$, and suppose by the induction hypothesis
that, for $k<n$, $v(s_1)=v(s_2)=k$ implies $s_1=s_2$.

If $n$ is a prime number (say the $k^{th}$ prime) then $s_1$ and
$s_2$ are irreducible (encircled), that is $s_1=\langle s_3
\rangle$ and $s_2=\langle s_4 \rangle$ and $v(s_3)=v(s_4)=k$.
Thus, by the induction hypothesis, $s_3=s_4$ and therefore
$s_1=s_2$.

If $n$ is a composite number then $s_1$ and $s_2$ factorize into
the same irreducible morphisms because if an irreducible morphism
$s_3$ is a factor of $s_1$ then $v(s_3)$ is a prime number that
divides $n$ and, by the induction hypothesis, $s_3$ is the unique
morphism with value $v(s_3)$, therefore $s_3$ is also a factor of
$s_2$. The same argument can be used to prove that each factor
appears in $s_1$ and $s_2$ the same number of times. We only need
to do the following exercise to conclude that $s_1=s_2$.

\begin{exercise}
Show that the monoid $\hom(\emptyset,\emptyset)$ is commutative.
\end{exercise}


{\bf Second case:} the monoid is the non-negative integer numbers
with the usual sum as the operation of the monoid and with the
usual order giving the lattice structure; the function $\varphi$
is the function that sends a number $n$ to its successor $n+1$.

We are going to see that, for a system of non-singular planar
curves $s$, $v(s)$ is simply the number of curves of $s$.

This is very easy because, since the application $v$ is uniquely
determined by the relations $v(s_1\circ s_2)=v(s_1)+v(s_2)$ and
$v(\langle s\rangle )=\varphi(v(s))=v(s)+1$, we only need to
observe that the number of curves $\nu(s)$ of a system of curves
$s$ satisfies the relation $\nu(s_1\circ s_2)=\nu(s_1)+\nu(s_2)$
and $\nu(\langle s\rangle )=\nu(s)+1$.



{\bf Acknowledgment} - I wish to thank all the people who helped
me to write this paper: L. Kauffman, R. Picken, P. Resende and M.
Stosic for their useful comments and the Department of
Combinatorics of the University of Waterloo where part of this
paper was written.

This work was supported by the programme {\em Programa Operacional
``Ci\^{e}ncia, Tecnologia, Inova\c{c}\~{a}o''} (POCTI) of the
{\em Funda\c{c}\~{a}o para a Ci\^{e}ncia e a Tecnologia} (FCT),
cofinanced by the European Community fund FEDER.


\begin{thebibliography}{1}
\bibitem[1]{1} Birkhoff, G. {\it Lattice Theory, 3rd ed.},
Providence, RI: Amer. Math. Soc., 1967.

\bibitem[2]{2} Ki Hang Kim, {\it Boolean
matrix theory and application}, Marcel Decker, Inc., 1982.



\bibitem[3]{3}  V.G.Turaev, {\it Operator invariants of tangles,
and $R$-matrices.}(Russian), Izv. Akad. Nauk SSSR Ser. Mat. 53
(1989), no. 5, 1073--1107, 1135; translation in Math. USSR-Izv. 35
(1990), no. 2, 411--444.





\end{thebibliography}
\end{document}